\documentclass[a4paper,11pt,oneside,reqno]{amsart}
\usepackage[T1]{fontenc}
\usepackage[utf8]{inputenc}
\usepackage[sc]{mathpazo}
\usepackage[a4paper,margin=26mm]{geometry}
\usepackage{amsmath,amssymb,amsthm,mathtools,mathrsfs}
\usepackage{array,booktabs,longtable,graphicx,tikz,needspace}
\usepackage[section]{placeins}
\usetikzlibrary{arrows.meta,calc,positioning}
\usepackage{microtype}
\usepackage[unicode,hidelinks]{hyperref}
\usepackage{bookmark}
\usepackage[all]{xy}
\usepackage{textcase}
\newcommand{\ResearchPart}[2]{%
  \par\bigskip\Needspace{10\baselineskip}%
  \refstepcounter{part}%
  \addcontentsline{toc}{part}{Part \thepart. #1}%
  \noindent\rule{\linewidth}{0.55pt}\par\vspace{0.7\baselineskip}%
  \noindent{\footnotesize\scshape\color{black!60}Part \thepart}\par\vspace{0.24\baselineskip}%
  \noindent{\Large\bfseries #1}\par\vspace{0.55\baselineskip}%
  \noindent{\small #2}\par\vspace{0.5\baselineskip}%
  \noindent\rule{\linewidth}{0.3pt}\par\medskip
}

\newtheorem{theorem}{Theorem}[section]
\newtheorem{lemma}[theorem]{Lemma}
\newtheorem{proposition}[theorem]{Proposition}
\newtheorem{corollary}[theorem]{Corollary}
\theoremstyle{definition}
\newtheorem{definition}[theorem]{Definition}
\newtheorem{example}[theorem]{Example}
\newtheorem{remark}[theorem]{Remark}

\numberwithin{equation}{section}
\allowdisplaybreaks[2]
\hypersetup{pdftitle={Hochschild cohomology of Temperley--Lieb algebras: a Kauffman diagram perspective},pdfauthor={Viktor Lopatkin},pdfsubject={Hochschild cohomology; Anick resolution; diagrammatic cocycles; Temperley--Lieb algebras}}
\title[Hochschild cohomology and Kauffman diagrams]{Hochschild cohomology of Temperley--Lieb algebras:\\ a Kauffman diagram perspective}
\author{Viktor Lopatkin}
\date{\today}

\address{Faculty of Computer Science, HSE University,
	Pokrovsky Boulevard 11, Moscow, 109028, Russia}

\email{wickktor@gmail.com}

\begin{document}
	
	\begin{abstract}
We describe Hochschild cohomology of Temperley--Lieb algebras by explicit
minimal-resolution coordinates and multilinear operations on Kauffman
diagrams. For each doubled-line basic algebra, a finite four-term
bimodule complex is verified by a signed-shuffle comparison and an
augmented contracting homotopy. We compute cocycles and coboundaries
and recover the known cup product and Gerstenhaber bracket in explicit
representatives. Over $\mathbb C$, at nonzero root-of-unity loop
parameters, we construct a tower-compatible derivation and the
marked-circle two-cocycle. Their classes are nonzero on every
nonsemisimple block; their cup products therefore represent all positive
classes over the centre. Relative full twists give a compact expression
for the derivation. Worked planar diagrams distinguish the cup product
from insertion and strict cochain equalities from equalities of classes.
A finite-lattice section interprets the marked-circle operation as
loop-weight response, derives its planar Potts normalization, and
exhibits a six-vertex realization. It also shows why passing to a
deformable spin representation may annihilate the class with different
coefficients. No thermodynamic conclusion is claimed, and the even-strand
zero-parameter family remains separate.
\end{abstract}
\maketitle

\begingroup
\small
\setlength{\parskip}{0pt}
\tableofcontents
\endgroup

\section{Introduction: from a minimal resolution to diagrammatic operations}
\label{tl:introduction}
Temperley--Lieb algebras connect the representation theory of finite-dimensional
algebras with planar topology and transfer-matrix methods for lattice models.
Their defining planar objects are sufficiently concrete to formulate the
cohomological problem at the level of the diagrams themselves.

\begin{definition}[Kauffman diagrams and the Temperley--Lieb algebra]
Let $N\ge1$. An \emph{$N$-strand Kauffman diagram} is an isotopy class,
relative to the boundary, of $N$ pairwise disjoint, nonintersecting arcs
in a rectangle that pair its $N$ marked upper and $N$ marked lower
boundary points. There are no closed components in a basis diagram.
For $\beta\in\mathbb C$, the \emph{Temperley--Lieb algebra}
$TL_N(\beta)$ is the complex vector space with these diagrams as a
basis. The product of two basis diagrams is defined by stacking the
first above the second, deleting each internal closed circle, and
multiplying the remaining diagram by $\beta^L$, where $L$ is the
number of deleted circles. The unit is the diagram of $N$ vertical
through strands.
\end{definition}
Equivalently, $TL_N(\beta)$ is generated by the elementary diagrams
$u_i$ with a cap and a cup joining the neighbouring positions
$i,i+1$ and through strands elsewhere. They satisfy
\begin{equation}\label{tl:relations}
 u_i^2=\beta u_i,\qquad u_iu_{i\pm1}u_i=u_i,
 \qquad u_iu_j=u_ju_i\quad(|i-j|>1).
\end{equation}
This diagrammatic presentation is valid for every $\beta$,
with $\dim TL_N(\beta)=\frac1{N+1}\binom{2N}{N}$;
see \cite[Theorem~2.4]{RSA} and Section~\ref{tl:prelim}.
The loop parameter originated in the transfer-matrix description
of percolation, colouring, and ice-type models \cite{TemperleyLieb};
the graphical correspondence with planar Potts models was developed
by Baxter--Kelland--Wu \cite{BKW}.

Our result gives explicit representatives of Hochschild cohomology
classes as operations on this fixed diagram basis. Thus the first
cohomology is represented by diagrammatic derivations modulo inner
ones, the second by bilinear corrections to diagram multiplication,
and their cup products and brackets become operations with specified
input and output diagrams.

The abstract cohomology ring of the nonsemisimple basic blocks is already
known. Li--Xu--Chen \cite{LiXuChen} constructed an explicit minimal
projective bimodule resolution and computed the cup algebra, and
Lyu--Wu \cite{LyuWu} determined the Gerstenhaber structure.
The present work is not a claim of the first minimal resolution or the
first ring presentation. We give a critical-path realization, uniform
four-term coordinates, explicit comparison and contracting maps, and
ordinary diagrammatic representatives of the low-degree generators.
Our main additional statements are their nonvanishing on every
nonsemisimple block and their realization by loop marking and relative
full twists. Worked diagrams and a finite-lattice interpretation then
connect these representatives with actual gluing operations.

\paragraph{The underlying bound quiver.}
For the nonsemisimple parameters considered below, write
$\beta=q+q^{-1}\ne0$ and $e=\operatorname{ord}(q^2)\ge3$.
Each nonsemisimple block $B$ of $TL_N(\beta)$, with $m$ isomorphism
classes of simple modules, is Morita equivalent to the \emph{bound}
quiver algebra $\Lambda(m)=\mathbb C Q_m/I_m$ \cite{Westbury,RSA}.
Here $Q_m$ is the doubled line
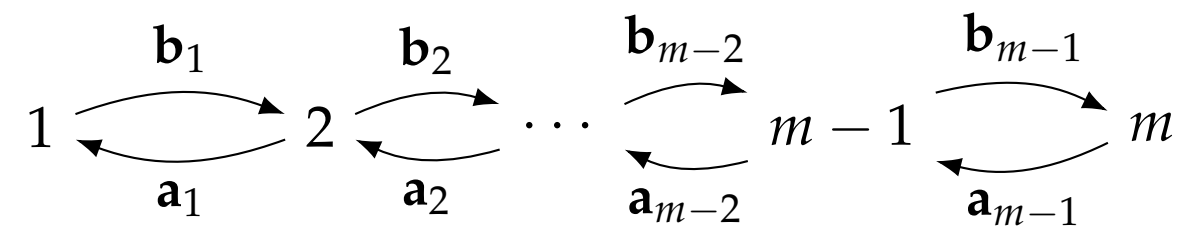
\begin{center}
\begin{tikzpicture}[>=Latex,baseline=(current bounding box.center)]
\node (q1) at (0,0) {$1$};
\node (q2) at (1.85,0) {$2$};
\node (qd) at (3.45,0) {$\cdots$};
\node (qm1) at (5.3,0) {$m-1$};
\node (qm) at (7.35,0) {$m$};
\draw[->,bend left=20] (q1) to node[above,font=\small] {$\mathbf b_1$} (q2);
\draw[->,bend left=20] (q2) to node[above,font=\small] {$\mathbf b_2$} (qd);
\draw[->,bend left=20] (qd) to node[above,font=\small] {$\mathbf b_{m-2}$} (qm1);
\draw[->,bend left=20] (qd) to node[below,font=\small] {$\mathbf a_2$} (q2);
\draw[->,bend left=20] (q2) to node[below,font=\small] {$\mathbf a_1$} (q1);
\draw[->,bend left=20] (qm1) to node[above,font=\small] {$\mathbf b_{m-1}$} (qm);
\draw[->,bend left=20] (qm1) to node[below,font=\small] {$\mathbf a_{m-2}$} (qd);
\draw[->,bend left=20] (qm) to node[below,font=\small] {$\mathbf a_{m-1}$} (qm1);
\end{tikzpicture}
\end{center}
with paths multiplied from left to right, and with relations
\begin{equation}\label{intro:quiver-relations}
\begin{gathered}
 \mathbf a_{i+1}\mathbf a_i=0,\quad
 \mathbf b_i\mathbf b_{i+1}=0,\quad
 \mathbf a_i\mathbf b_i=\mathbf b_{i+1}\mathbf a_{i+1}
 \quad(1\le i\le m-2),\\
 \mathbf a_{m-1}\mathbf b_{m-1}=0.
\end{gathered}
\end{equation}
The central block idempotent, the primitive vertex idempotents,
and the physical strands of a Kauffman diagram are different objects:
$m$ counts the simple modules of $B$, not the strands $N$.
The general presentation is given again in Section~\ref{sec:GS}
for the resolution calculations. In the first illustrative case below
there is an actual algebra isomorphism $TL_3(1)\cong\Lambda(2)$.

\subsection{A first diagrammatic glimpse: $TL_3(1)$}
\begin{figure}[!htbp]
\centering
\smallskip
\begin{tikzpicture}[line width=.75pt,scale=.93]
\begin{scope}[shift={(0,0)}]
\draw (0,2.08) -- (0,.62);
\draw (.43,2.08) .. controls (.43,1.69) and (.86,1.69) .. (.86,2.08);
\draw (.43,.62) .. controls (.43,1.01) and (.86,1.01) .. (.86,.62);
\node at (.43,.21) {$\scriptstyle v$};
\end{scope}
\node at (1.65,1.37) {$\xrightarrow{\ D\ }$};
\begin{scope}[shift={(2.75,0)}]
\draw (0,2.08) -- (0,.62);
\draw (.43,2.08) .. controls (.43,1.69) and (.86,1.69) .. (.86,2.08);
\draw (.43,.62) .. controls (.43,1.01) and (.86,1.01) .. (.86,.62);
\node at (.43,.21) {$\scriptstyle v$};
\end{scope}
\node at (4.08,1.37) {$-$};
\begin{scope}[shift={(4.74,0)}]
\draw (0,2.08) .. controls (0,1.69) and (.43,1.69) .. (.43,2.08);
\draw (.43,.62) .. controls (.43,1.01) and (.86,1.01) .. (.86,.62);
\draw (.86,2.08) .. controls (.86,1.55) and (0,1.15) .. (0,.62);
\node at (.43,.21) {$\scriptstyle uv$};
\end{scope}
\end{tikzpicture}
\caption{The explicit outer class $[D]\in HH^1(TL_3(1))$:
$D(v)=v-uv$. Each term on the right is an ordinary three-strand
Kauffman diagram; the formula is an equality of algebra elements.}
\label{fig:intro-TL3}
\end{figure}
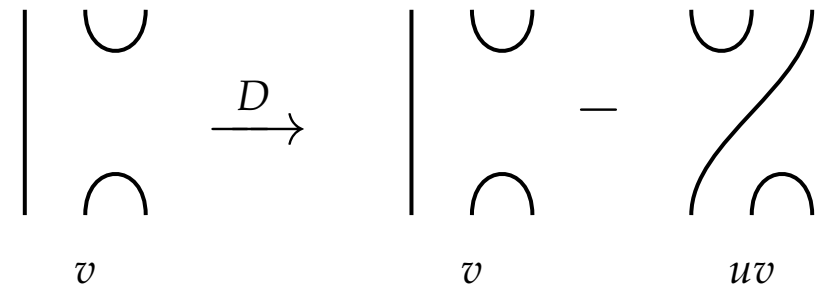

Put $A=TL_3(1)$, $u=u_1$ and $v=u_2$. Already at this level the first
Hochschild cohomology is nonzero. More precisely,
\begin{equation}\label{intro:HH1-explicit}
 \boxed{HH^1(A)=\operatorname{Der}(A)/\operatorname{Inn}(A)
       =\mathbb C[D],\qquad D(u)=0,\quad D(v)=v-uv.}
\end{equation}
Thus the \emph{outer} class $[D]$ is represented by the explicit
one-input rule in Figure~\ref{fig:intro-TL3}. The second cohomology
is also one-dimensional, with generator the class of the marked-circle
operation $\Phi$ obtained by differentiating diagram multiplication
with respect to its loop weight. These are conclusions, not merely
candidates for representatives; their verification is given below.

\begin{figure}[!htbp]
\centering
\smallskip
\begin{tikzpicture}[line width=.72pt,scale=1.12]
\draw[gray!60,line width=.45pt,rounded corners=2.3pt]
  (-.18,-3.43) rectangle (7.63,3.08);
\draw[gray!65,line width=.45pt] (-.18,-.28) -- (7.63,-.28);
\draw[gray!48,line width=.45pt] (3.78,-3.43) -- (3.78,3.08);
\begin{scope}[shift={(0.00,0)}]
\node[font=\small] at (1.58,2.73) {$\Phi(u,u)$};
\draw[gray!52,line width=.36pt] (-.10,1.17) -- (1.14,1.17);
\draw (0.080,1.170) .. controls (0.080,0.860) and (0.560,0.860) .. (0.560,1.170);
\draw (0.080,0.120) .. controls (0.080,0.430) and (0.560,0.430) .. (0.560,0.120);
\draw (1.040,0.120) -- (1.040,1.170);
\draw (0.080,2.220) .. controls (0.080,1.910) and (0.560,1.910) .. (0.560,2.220);
\draw (0.080,1.170) .. controls (0.080,1.480) and (0.560,1.480) .. (0.560,1.170);
\draw (1.040,1.170) -- (1.040,2.220);
\node[font=\scriptsize] at (1.25,1.75) {$ u $};
\node[font=\scriptsize] at (1.25,.62) {$ u $};
\node[font=\small,text=blue!65!black] at (0.320,1.17) {$\star$};
\node at (1.86,1.20) {$\xrightarrow{\ \Phi_1\ }$};
\draw (2.630,1.700) .. controls (2.630,1.390) and (3.110,1.390) .. (3.110,1.700);
\draw (2.630,0.650) .. controls (2.630,0.960) and (3.110,0.960) .. (3.110,0.650);
\draw (3.590,0.650) -- (3.590,1.700);
\end{scope}
\begin{scope}[shift={(3.85,0)}]
\node[font=\small] at (1.58,2.73) {$\Phi(u,v)$};
\draw[gray!52,line width=.36pt] (-.10,1.17) -- (1.14,1.17);
\draw (0.560,1.170) .. controls (0.560,0.860) and (1.040,0.860) .. (1.040,1.170);
\draw (0.560,0.120) .. controls (0.560,0.430) and (1.040,0.430) .. (1.040,0.120);
\draw (0.080,0.120) -- (0.080,1.170);
\draw (0.080,2.220) .. controls (0.080,1.910) and (0.560,1.910) .. (0.560,2.220);
\draw (0.080,1.170) .. controls (0.080,1.480) and (0.560,1.480) .. (0.560,1.170);
\draw (1.040,1.170) -- (1.040,2.220);
\node[font=\scriptsize] at (1.25,1.75) {$ u $};
\node[font=\scriptsize] at (1.25,.62) {$ v $};
\node at (1.86,1.20) {$\xrightarrow{\ \Phi_1\ }$};
\node at (2.90,1.18) {$0$};
\end{scope}
\begin{scope}[shift={(0.00,-3.25)}]
\node[font=\small] at (1.58,2.73) {$\Phi(v,u)$};
\draw[gray!52,line width=.36pt] (-.10,1.17) -- (1.14,1.17);
\draw (0.080,1.170) .. controls (0.080,0.860) and (0.560,0.860) .. (0.560,1.170);
\draw (0.080,0.120) .. controls (0.080,0.430) and (0.560,0.430) .. (0.560,0.120);
\draw (1.040,0.120) -- (1.040,1.170);
\draw (0.560,2.220) .. controls (0.560,1.910) and (1.040,1.910) .. (1.040,2.220);
\draw (0.560,1.170) .. controls (0.560,1.480) and (1.040,1.480) .. (1.040,1.170);
\draw (0.080,1.170) -- (0.080,2.220);
\node[font=\scriptsize] at (1.25,1.75) {$ v $};
\node[font=\scriptsize] at (1.25,.62) {$ u $};
\node at (1.86,1.20) {$\xrightarrow{\ \Phi_1\ }$};
\node at (2.90,1.18) {$0$};
\end{scope}
\begin{scope}[shift={(3.85,-3.25)}]
\node[font=\small] at (1.58,2.73) {$\Phi(v,v)$};
\draw[gray!52,line width=.36pt] (-.10,1.17) -- (1.14,1.17);
\draw (0.560,1.170) .. controls (0.560,0.860) and (1.040,0.860) .. (1.040,1.170);
\draw (0.560,0.120) .. controls (0.560,0.430) and (1.040,0.430) .. (1.040,0.120);
\draw (0.080,0.120) -- (0.080,1.170);
\draw (0.560,2.220) .. controls (0.560,1.910) and (1.040,1.910) .. (1.040,2.220);
\draw (0.560,1.170) .. controls (0.560,1.480) and (1.040,1.480) .. (1.040,1.170);
\draw (0.080,1.170) -- (0.080,2.220);
\node[font=\scriptsize] at (1.25,1.75) {$ v $};
\node[font=\scriptsize] at (1.25,.62) {$ v $};
\node[font=\small,text=blue!65!black] at (0.800,1.17) {$\star$};
\node at (1.86,1.20) {$\xrightarrow{\ \Phi_1\ }$};
\draw (3.110,1.700) .. controls (3.110,1.390) and (3.590,1.390) .. (3.590,1.700);
\draw (3.110,0.650) .. controls (3.110,0.960) and (3.590,0.960) .. (3.590,0.650);
\draw (2.630,0.650) -- (2.630,1.700);
\end{scope}
\end{tikzpicture}
\caption{Four elementary evaluations of the marked-circle cocycle
on the three-strand diagram basis, arranged as a single calculation
(top row: $u\cdot u$, $u\cdot v$; bottom row: $v\cdot u$,
$v\cdot v$). The first factor is stacked above the second.
An internal contractible circle, marked by $\star$, occurs in
$u\cdot u$ and $v\cdot v$, but not in the two mixed products.
The frame rules separate the four cases, and the thin pale line
within each cell marks the interface of the two stacked input factors;
none of the gray rules represents a strand. The values are
$\Phi_1(u,u)=u$, $\Phi_1(u,v)=0$,
$\Phi_1(v,u)=0$, $\Phi_1(v,v)=v$.}
\label{fig:intro-TL3-Phi-four}
\end{figure}
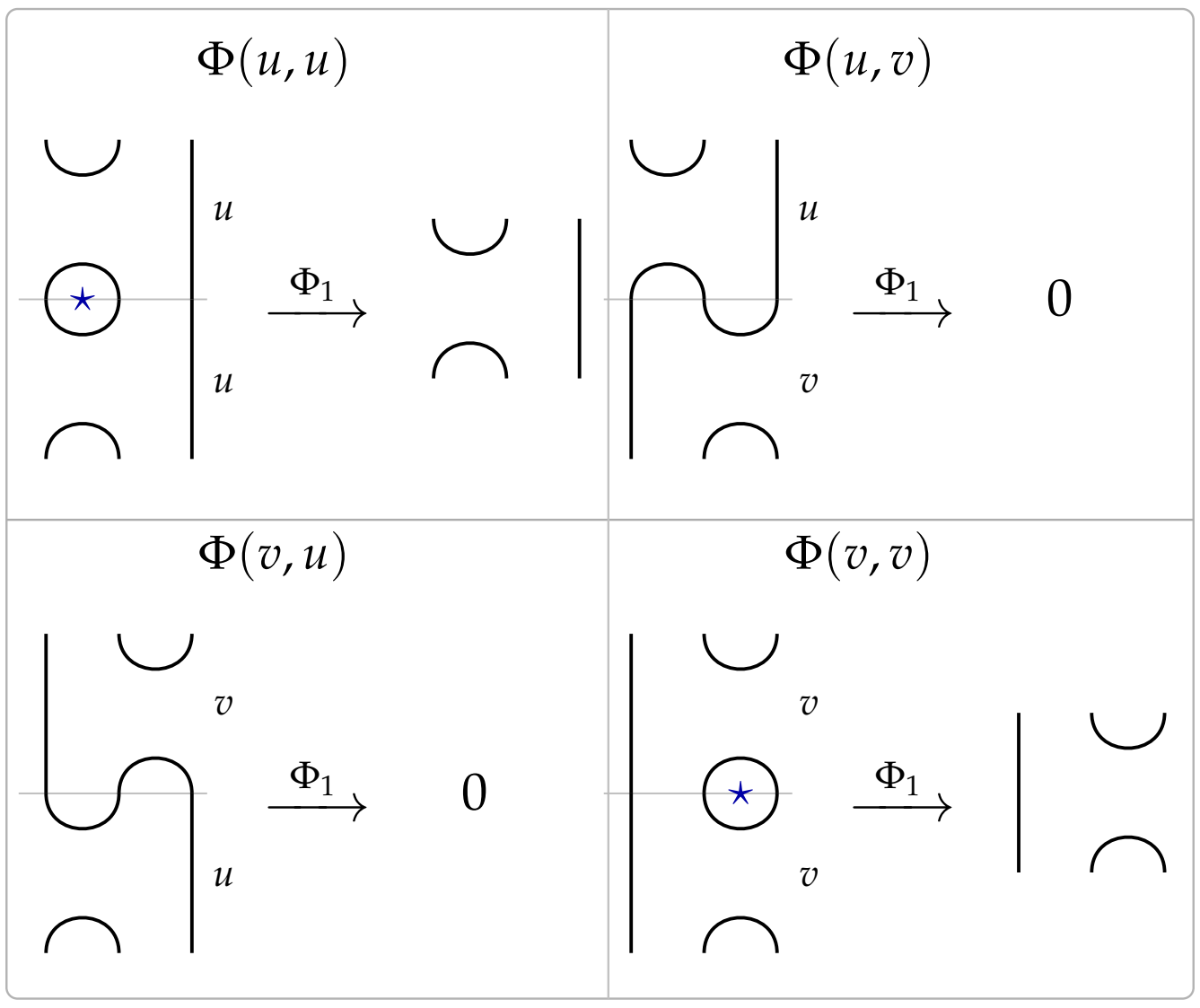

For the second Hochschild generator, the marking rule is equally
explicit. If the stacked product $d_1d_2$ creates $L$ closed circles
and leaves the loop-free diagram $d_0$, then
$\Phi_\beta(d_1,d_2)=L\beta^{L-1}d_0$ for $L>0$ and is zero for
$L=0$. Hence, with $\beta=1$,
\begin{equation}\label{intro:Phi-four}
 \boxed{\Phi(u,u)=u,\quad \Phi(u,v)=0,\quad
 \Phi(v,u)=0,\quad \Phi(v,v)=v.}
\end{equation}
The four input pairs, the internal circles and the outputs are drawn in
Figure~\ref{fig:intro-TL3-Phi-four}. The operation has \emph{two}
inputs, unlike the derivation $D$ with one input. The nonzero class
$[\Phi]\in HH^2(A)$ is established below; an individual nonzero
cochain value alone would not imply that conclusion.

The same algebra provides a complete small test case. Put
$z=u+v-uv-vu$. Then $z\in Z(A)$, $z^2=0$ and $D(z)=z\ne0$, so
$D$ is not inner. Its Hochschild cohomology is
\begin{equation}\label{eq:intro-HH3}
\boxed{\begin{gathered}
HH^0(A)=\mathbb C1\oplus\mathbb Cz,\qquad
HH^1(A)=\mathbb C[D],\qquad HH^2(A)=\mathbb C[\Phi],\\
HH^n(A)=0\quad(n\ge3).
\end{gathered}}
\end{equation}
All cup products of two positive-degree classes vanish in this
small example; nevertheless the Gerstenhaber bracket is nontrivial:
$\bigl[[D],[\Phi]\bigr]=-[\Phi]$ and $[[D],z]=z$.
This is why the first nonzero higher cup product is illustrated
later on a larger block: in $TL_6(1)$, a block $B$ of type
$\Lambda(3)$ satisfies $[\Phi_B]^{\smile2}\ne0$.
The full cochain computations and the distinction between cochain
equalities and identities in cohomology are presented in
Section~\ref{vis:calculus}.
\FloatBarrier

\subsection{A minimal resolution in critical-path coordinates.}
For the doubled-line algebra $\Lambda(m)$, we state the projective terms
and the four-term differential first. The generators are the formal paths
$\mathrm C_{k;u,v}$ consisting of a descending run followed by an ascending
run. These include closed excursions when $u=v$, but also open paths:
the resolution is not made only of loops. A signed shuffle map proves
$\mathsf d^2=0$, and the augmented contraction
$\mathsf h=p\mathsf s\kappa$ proves exactness. Minimality follows from
exterior radical coefficients. The bound on the allowed indices makes
termination visible directly, without first truncating a longer complex
at a separately supplied global dimension. The Morse-theoretic origin
is recorded briefly, rather than reproduced as a long construction.
The same formulas provide the coefficient complex for any bimodule
(Appendix~\ref{app:coefficients}).

\subsection{Explicit cocycles and coboundaries.}
With regular coefficients, the parallel-path restriction reduces the
cochain equations to signed difference operators. We solve them and
exhibit the representatives $\eta_\ell$ and $\theta_\ell$, including
explicit primitives for coboundaries. This recovers the known dimensions
while specifying the actual values of the cochains on the resolution.
The formulas are valid over an arbitrary field in the quiver part.

\subsection{Products and brackets on chosen representatives.}
A signed shuffle diagonal computes the cup product. Its admissible cuts
and signs are part of the formula, not an informal geometric convention.
A weight derivation and a flat quiver deformation compute the
Gerstenhaber bracket, reproducing the known algebra in the chosen
coordinates. We distinguish strict cochain identities from identities
only in cohomology; Section~\ref{vis:calculus} illustrates this distinction
with an explicit coboundary in $TL_3(1)$.

\subsection{ Operations on Kauffman diagrams.}
Over $\mathbb C$ and at a nonzero nonsemisimple parameter, we identify
ordinary cochains $D_N$ and $\Phi_\beta$ without choosing every primitive
projector of every block. The first is a compatible derivation of the TL
tower. The second marks one circle created by a diagram product.
Their restrictions represent nonzero classes in degrees one and two on
every nonsemisimple block. Consequently their cup products give all
positive-degree classes, with central idempotents separating the blocks.
The formula for $D_N$ is made compact using relative full twists;
small examples are also drawn after complete planar expansion.

\subsection{A finite-lattice interpretation, with its scope explicit.}
The marked-circle rule is the first variation of the loop weight and,
under fixed local weights, inserts the loop-counting observable in a
finite partition sum. We derive the planar Potts normalization and
exhibit a two-spin TL representation whose local face operator has six
ice-rule matrix entries. The derivation $D$ gives an infinitesimal
automorphism of the gluing algebra, and the bracket describes the
infinitesimal transport of the marked-circle cochain. These are precise
algebraic statements relevant to loop and ice-type formulations, not a
calculation of critical exponents or a claim that every Hochschild class
is already a scalar physical observable. A deformable spin representation
can make $\rho\Phi$ a coboundary with $\operatorname{End}(V)$ coefficients;
we explain this explicitly in Section~\ref{phys:section}.

\subsection{The main algebraic statement}
The vertex number $m$ of a basic block is not the number $N$ of TL
strands. The coefficient bimodule is the regular one unless stated
otherwise.
\begin{theorem}[Overview]\label{tl:overview}
Let $\Lambda(m)$, $m\ge2$, be the doubled-line algebra with relations
\eqref{eq:ideal}.
\begin{enumerate}
\item Over a field $\Bbbk$, the complex of Theorem~\ref{thm:resolution}
is a minimal bimodule resolution of length $2m-2$, with generators
$\mathrm C_{k;u,v}$, $u+v=n$, $k+\max(u,v)\le m$, and differential
\eqref{eq:main-differential}. Its cohomology representatives, cup product,
and bracket are given in Theorems~\ref{thm:HH}, \ref{cup:thm:ring},
and \ref{gerst:complete-bracket}.
\item Let $A_N=TL_N(\beta)$ over $\mathbb C$, where
$\beta=q+q^{-1}\ne0$, $e=\operatorname{ord}(q^2)\ge3$, and
$\sigma=q^e\in\{1,-1\}$. Put
\[
 H_n=\frac{1+\sigma^{n+1}P_{e-2}(F_n)}e,
 \qquad D_N(u_i)=(H_i-H_{i-1})u_i.
\]
Then $D_N$ is a derivation, compatible with the tower inclusions.
If $d_1d_2=\beta^Ld_0$, define
\[
 \Phi_\beta(d_1,d_2)=
 \begin{cases}L\beta^{L-1}d_0,&L>0,\\0,&L=0.\end{cases}
\]
For each nonsemisimple block $B$,
$HH^1(B)=\mathbb C[D_B]$ and $HH^2(B)=\mathbb C[\Phi_B]$:
both displayed cochains represent nonzero classes.
\item If $B$ has $m$ simple modules, its positive cohomology has
representatives $\Phi_B^{\smile\ell}$ for $1\le\ell\le m-1$ and
$\Phi_B^{\smile\ell}\smile D_B$ for $0\le\ell\le m-2$.
Thus $HH^1(A_N)=Z(A_N)[D_N]$, $HH^2(A_N)=Z(A_N)[\Phi_\beta]$,
and these classes generate the positive cup algebra over the centre.
\end{enumerate}
\end{theorem}
The centre remains degree zero; it is not generated by positive cup
powers. Nor does cyclicity over the centre mean that $HH^1(A_N)$ or
$HH^2(A_N)$ is one-dimensional when several blocks are present.
The even-strand family at $\beta=0$ is expressly excluded from the
blockwise generator theorem and treated separately below.

\begin{figure}[htbp]
\centering
\begin{tikzpicture}[line width=.78pt,line cap=round,line join=round]
\begin{scope}[shift={(0.00,0)}]
\draw (0.000,4.800) -- (0.000,2.700);
\draw (0.000,2.100) -- (0.000,0);
\draw (0.660,4.800) -- (0.660,2.700);
\draw (0.660,2.100) -- (0.660,0);
\draw (0.660,2.700) .. controls (0.660,2.590) and (0.000,2.510) .. (0.000,2.400);
\draw[white,line width=3pt] (0.000,2.700) .. controls (0.000,2.590) and (0.660,2.510) .. (0.660,2.400);
\draw (0.000,2.700) .. controls (0.000,2.590) and (0.660,2.510) .. (0.660,2.400);
\draw (0.660,2.400) .. controls (0.660,2.290) and (0.000,2.210) .. (0.000,2.100);
\draw[white,line width=3pt] (0.000,2.400) .. controls (0.000,2.290) and (0.660,2.210) .. (0.660,2.100);
\draw (0.000,2.400) .. controls (0.000,2.290) and (0.660,2.210) .. (0.660,2.100);
\fill (0.000,4.800) circle (.85pt);
\fill (0.000,0) circle (.85pt);
\fill (0.660,4.800) circle (.85pt);
\fill (0.660,0) circle (.85pt);
\node[font=\small] at (0.330,-0.55) {$C_2=t_1^2$};
\end{scope}
\begin{scope}[shift={(4.20,0)}]
\draw (0.000,4.800) -- (0.000,3.300);
\draw (0.000,1.500) -- (0.000,0);
\draw (0.660,4.800) -- (0.660,3.300);
\draw (0.660,1.500) -- (0.660,0);
\draw (1.320,4.800) -- (1.320,3.300);
\draw (1.320,1.500) -- (1.320,0);
\draw (1.320,3.300) -- (1.320,3.000);
\draw (0.660,3.300) .. controls (0.660,3.190) and (0.000,3.110) .. (0.000,3.000);
\draw[white,line width=3pt] (0.000,3.300) .. controls (0.000,3.190) and (0.660,3.110) .. (0.660,3.000);
\draw (0.000,3.300) .. controls (0.000,3.190) and (0.660,3.110) .. (0.660,3.000);
\draw (0.000,3.000) -- (0.000,2.700);
\draw (1.320,3.000) .. controls (1.320,2.890) and (0.660,2.810) .. (0.660,2.700);
\draw[white,line width=3pt] (0.660,3.000) .. controls (0.660,2.890) and (1.320,2.810) .. (1.320,2.700);
\draw (0.660,3.000) .. controls (0.660,2.890) and (1.320,2.810) .. (1.320,2.700);
\draw (1.320,2.700) -- (1.320,2.400);
\draw (0.660,2.700) .. controls (0.660,2.590) and (0.000,2.510) .. (0.000,2.400);
\draw[white,line width=3pt] (0.000,2.700) .. controls (0.000,2.590) and (0.660,2.510) .. (0.660,2.400);
\draw (0.000,2.700) .. controls (0.000,2.590) and (0.660,2.510) .. (0.660,2.400);
\draw (0.000,2.400) -- (0.000,2.100);
\draw (1.320,2.400) .. controls (1.320,2.290) and (0.660,2.210) .. (0.660,2.100);
\draw[white,line width=3pt] (0.660,2.400) .. controls (0.660,2.290) and (1.320,2.210) .. (1.320,2.100);
\draw (0.660,2.400) .. controls (0.660,2.290) and (1.320,2.210) .. (1.320,2.100);
\draw (1.320,2.100) -- (1.320,1.800);
\draw (0.660,2.100) .. controls (0.660,1.990) and (0.000,1.910) .. (0.000,1.800);
\draw[white,line width=3pt] (0.000,2.100) .. controls (0.000,1.990) and (0.660,1.910) .. (0.660,1.800);
\draw (0.000,2.100) .. controls (0.000,1.990) and (0.660,1.910) .. (0.660,1.800);
\draw (0.000,1.800) -- (0.000,1.500);
\draw (1.320,1.800) .. controls (1.320,1.690) and (0.660,1.610) .. (0.660,1.500);
\draw[white,line width=3pt] (0.660,1.800) .. controls (0.660,1.690) and (1.320,1.610) .. (1.320,1.500);
\draw (0.660,1.800) .. controls (0.660,1.690) and (1.320,1.610) .. (1.320,1.500);
\fill (0.000,4.800) circle (.85pt);
\fill (0.000,0) circle (.85pt);
\fill (0.660,4.800) circle (.85pt);
\fill (0.660,0) circle (.85pt);
\fill (1.320,4.800) circle (.85pt);
\fill (1.320,0) circle (.85pt);
\node[font=\small] at (0.660,-0.55) {$C_3=(t_1t_2)^3$};
\end{scope}
\begin{scope}[shift={(8.90,0)}]
\draw (0.000,4.800) -- (0.000,4.200);
\draw (0.000,0.600) -- (0.000,0);
\draw (0.660,4.800) -- (0.660,4.200);
\draw (0.660,0.600) -- (0.660,0);
\draw (1.320,4.800) -- (1.320,4.200);
\draw (1.320,0.600) -- (1.320,0);
\draw (1.980,4.800) -- (1.980,4.200);
\draw (1.980,0.600) -- (1.980,0);
\draw (1.320,4.200) -- (1.320,3.900);
\draw (1.980,4.200) -- (1.980,3.900);
\draw (0.660,4.200) .. controls (0.660,4.090) and (0.000,4.010) .. (0.000,3.900);
\draw[white,line width=3pt] (0.000,4.200) .. controls (0.000,4.090) and (0.660,4.010) .. (0.660,3.900);
\draw (0.000,4.200) .. controls (0.000,4.090) and (0.660,4.010) .. (0.660,3.900);
\draw (0.000,3.900) -- (0.000,3.600);
\draw (1.980,3.900) -- (1.980,3.600);
\draw (1.320,3.900) .. controls (1.320,3.790) and (0.660,3.710) .. (0.660,3.600);
\draw[white,line width=3pt] (0.660,3.900) .. controls (0.660,3.790) and (1.320,3.710) .. (1.320,3.600);
\draw (0.660,3.900) .. controls (0.660,3.790) and (1.320,3.710) .. (1.320,3.600);
\draw (0.000,3.600) -- (0.000,3.300);
\draw (0.660,3.600) -- (0.660,3.300);
\draw (1.980,3.600) .. controls (1.980,3.490) and (1.320,3.410) .. (1.320,3.300);
\draw[white,line width=3pt] (1.320,3.600) .. controls (1.320,3.490) and (1.980,3.410) .. (1.980,3.300);
\draw (1.320,3.600) .. controls (1.320,3.490) and (1.980,3.410) .. (1.980,3.300);
\draw (1.320,3.300) -- (1.320,3.000);
\draw (1.980,3.300) -- (1.980,3.000);
\draw (0.660,3.300) .. controls (0.660,3.190) and (0.000,3.110) .. (0.000,3.000);
\draw[white,line width=3pt] (0.000,3.300) .. controls (0.000,3.190) and (0.660,3.110) .. (0.660,3.000);
\draw (0.000,3.300) .. controls (0.000,3.190) and (0.660,3.110) .. (0.660,3.000);
\draw (0.000,3.000) -- (0.000,2.700);
\draw (1.980,3.000) -- (1.980,2.700);
\draw (1.320,3.000) .. controls (1.320,2.890) and (0.660,2.810) .. (0.660,2.700);
\draw[white,line width=3pt] (0.660,3.000) .. controls (0.660,2.890) and (1.320,2.810) .. (1.320,2.700);
\draw (0.660,3.000) .. controls (0.660,2.890) and (1.320,2.810) .. (1.320,2.700);
\draw (0.000,2.700) -- (0.000,2.400);
\draw (0.660,2.700) -- (0.660,2.400);
\draw (1.980,2.700) .. controls (1.980,2.590) and (1.320,2.510) .. (1.320,2.400);
\draw[white,line width=3pt] (1.320,2.700) .. controls (1.320,2.590) and (1.980,2.510) .. (1.980,2.400);
\draw (1.320,2.700) .. controls (1.320,2.590) and (1.980,2.510) .. (1.980,2.400);
\draw (1.320,2.400) -- (1.320,2.100);
\draw (1.980,2.400) -- (1.980,2.100);
\draw (0.660,2.400) .. controls (0.660,2.290) and (0.000,2.210) .. (0.000,2.100);
\draw[white,line width=3pt] (0.000,2.400) .. controls (0.000,2.290) and (0.660,2.210) .. (0.660,2.100);
\draw (0.000,2.400) .. controls (0.000,2.290) and (0.660,2.210) .. (0.660,2.100);
\draw (0.000,2.100) -- (0.000,1.800);
\draw (1.980,2.100) -- (1.980,1.800);
\draw (1.320,2.100) .. controls (1.320,1.990) and (0.660,1.910) .. (0.660,1.800);
\draw[white,line width=3pt] (0.660,2.100) .. controls (0.660,1.990) and (1.320,1.910) .. (1.320,1.800);
\draw (0.660,2.100) .. controls (0.660,1.990) and (1.320,1.910) .. (1.320,1.800);
\draw (0.000,1.800) -- (0.000,1.500);
\draw (0.660,1.800) -- (0.660,1.500);
\draw (1.980,1.800) .. controls (1.980,1.690) and (1.320,1.610) .. (1.320,1.500);
\draw[white,line width=3pt] (1.320,1.800) .. controls (1.320,1.690) and (1.980,1.610) .. (1.980,1.500);
\draw (1.320,1.800) .. controls (1.320,1.690) and (1.980,1.610) .. (1.980,1.500);
\draw (1.320,1.500) -- (1.320,1.200);
\draw (1.980,1.500) -- (1.980,1.200);
\draw (0.660,1.500) .. controls (0.660,1.390) and (0.000,1.310) .. (0.000,1.200);
\draw[white,line width=3pt] (0.000,1.500) .. controls (0.000,1.390) and (0.660,1.310) .. (0.660,1.200);
\draw (0.000,1.500) .. controls (0.000,1.390) and (0.660,1.310) .. (0.660,1.200);
\draw (0.000,1.200) -- (0.000,0.900);
\draw (1.980,1.200) -- (1.980,0.900);
\draw (1.320,1.200) .. controls (1.320,1.090) and (0.660,1.010) .. (0.660,0.900);
\draw[white,line width=3pt] (0.660,1.200) .. controls (0.660,1.090) and (1.320,1.010) .. (1.320,0.900);
\draw (0.660,1.200) .. controls (0.660,1.090) and (1.320,1.010) .. (1.320,0.900);
\draw (0.000,0.900) -- (0.000,0.600);
\draw (0.660,0.900) -- (0.660,0.600);
\draw (1.980,0.900) .. controls (1.980,0.790) and (1.320,0.710) .. (1.320,0.600);
\draw[white,line width=3pt] (1.320,0.900) .. controls (1.320,0.790) and (1.980,0.710) .. (1.980,0.600);
\draw (1.320,0.900) .. controls (1.320,0.790) and (1.980,0.710) .. (1.980,0.600);
\fill (0.000,4.800) circle (.85pt);
\fill (0.000,0) circle (.85pt);
\fill (0.660,4.800) circle (.85pt);
\fill (0.660,0) circle (.85pt);
\fill (1.320,4.800) circle (.85pt);
\fill (1.320,0) circle (.85pt);
\fill (1.980,4.800) circle (.85pt);
\fill (1.980,0) circle (.85pt);
\node[font=\small] at (0.990,-0.55) {$C_4=(t_1t_2t_3)^4$};
\end{scope}
\end{tikzpicture}
\caption{Full twists on two, three, and four strands. In each panel,
$C_n=(t_1\cdots t_{n-1})^n$ is read from top to bottom; the braid
returns every strand to its original endpoint. Crossings abbreviate
$t_i=1-qu_i$ in the Temperley--Lieb algebra; the pictures are braid
representatives, not additional planar basis elements.}
\label{fig:intro-full-twists}
\end{figure}
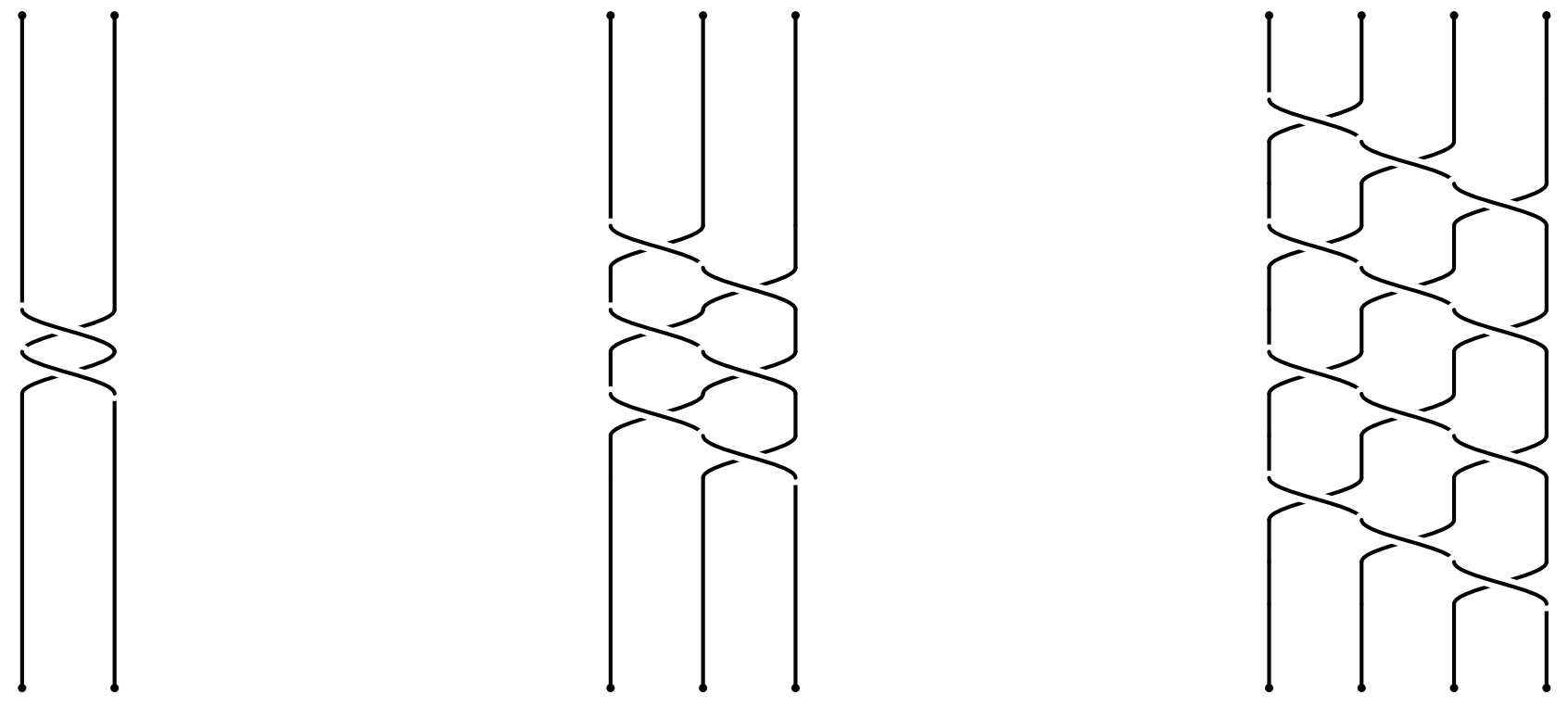

For the compact diagram rule set $t_i=1-qu_i$,
$C_n=(t_1\cdots t_{n-1})^n$, $J_n=C_nC_{n-1}^{-1}$, and
$R_n=q^{-n-1}J_n$ (with $R_0=q^{-1}$). Then
\begin{equation}\label{intro:twist-rule}
 F_n=R_n+R_n^{-1},\qquad
 P_{e-2}(F_n)=\sum_{j=0}^{e-2}R_n^{e-2-2j}.
\end{equation}
Figure~\ref{fig:intro-full-twists} displays the first three full twists,
which underlie the compact expression for the derivation.

The symbols $C_n$ for full twists and $\mathrm C_{k;u,v}$ for formal
resolution generators describe different objects. A cochain is a
multilinear rule on diagrams, not one unmarked diagram. Keeping its
inputs, coefficients, and signs is essential to every illustration.

\subsection{Comparison with previous work and organization}
The useful distinction from \cite{LiXuChen} is the choice of coordinates
and comparison data, not a denial of the explicitness or minimality of
that earlier resolution. Our terms are indexed directly by critical
paths; the differential is nonrecursive in those indices and has at
most four terms, and its length is read off from their inequalities.
The contracting map and the shuffle comparison make the passage to
ordinary bar operations explicit. Minimal resolutions are unique up
to chain isomorphism; a different presentation should not be confused
with a different derived invariant.

Relative bar constructions and algebraic Morse transfer are standard
\cite{Cibils,CLMS,Skoldberg,JollenbeckWelker,CLZ}; they are used here,
not claimed as new general methods. The moving-hole viewpoint is related to earlier Morse-theoretic
calculations \cite{Lopatkin}; the present proof verifies the
differential directly using comparison maps. We then compute its regular-coefficient
cohomology and operations, return to TL by Morita invariance, and
construct $D,\Phi$ directly on diagrams. The proof of blockwise
nonvanishing precedes the worked visual calculus. The final lattice-model
section combines established loop/ice dictionaries with consequences of
the cochain formulas, stating which parts are standard and which are
our deductions. No new phase diagram or solution of a lattice model
is claimed.

\section{Parameters, blocks, and cohomological conventions}\label{tl:prelim}
The defining diagrams, multiplication, and relations
\eqref{tl:relations} were given in the introduction. We use the
same loop-free basis for all parameters $\beta$, so only its multiplication
varies with $\beta$. Its dimension is always the Catalan number
$\mathrm{Cat}_N=\binom{2N}{N}/(N+1)$ \cite[Theorem~2.4]{RSA}.
For the tower write $A_n=TL_n(\beta)$, $A_0=A_1=\mathbb C$.

\subsection{Semisimplicity and reflection blocks}
Write $\beta=q+q^{-1}$ with $q\ne0$. Over $\mathbb C$, $TL_N(\beta)$
is semisimple for all $N$ when $q$ is not a root of unity, and also
at $\beta=\pm2$. If $e=\operatorname{ord}(q^2)\ge3$, it is
semisimple for $N<e$ and nonsemisimple for $N\ge e$; $\beta=0$
is exceptional, with semisimplicity for odd $N$
\cite[Corollaries~4.6, 4.8 and Proposition~8.2]{RSA}.
In the semisimple case $HH^{>0}(TL_N(\beta))=0$, while
$HH^0(TL_N(\beta))=Z(TL_N(\beta))$.
The nonsemisimple family includes $2\cos(r\pi/e)$ with
$\gcd(r,e)=1$, not only $r=1$; the sign $\sigma=q^e\in\{\pm1\}$
will be retained in formulas involving roots of unity.

The standard module $V_{N,p}$ has $p$ links and $N-2p$ defects;
set
\[
 \dim V_{N,p}=\binom Np-\binom N{p-1},\qquad
 h_p=N-2p+1\qquad(0\le p\le\lfloor N/2\rfloor).
\]
Its diagrammatic action and invariant Gram form are as in
\cite{RSA}. These representation-theoretic results are used,
not reproved below. A critical label, $e\mid h_p$, gives a
separate simple block. The noncritical labels in one block form
a reflection orbit, characterized within the fixed parity by
\begin{equation}\label{tl:blockrule}
 h_p\equiv\pm h_{p'}\pmod{2e}.
\end{equation}
Each nonsemisimple block with $m$ simple modules has basic algebra
$\Lambda(m)$, as presented in \eqref{intro:quiver-relations}
\cite[Section~10]{Westbury}; the quiver vertex count $m$ need not
equal the strand count $N$.

\subsection{Hochschild cochains and coefficient conventions}
For a unital algebra $A$, let $A^e=A\otimes A^{\mathrm{op}}$;
with regular coefficients,
$HH^n(A)=\operatorname{Ext}_{A^e}^n(A,A)$.
An ordinary cochain $f\colon A^{\otimes r}\to A$ has coboundary
\begin{align*}
 (\delta f)(a_1,\ldots,a_{r+1})={}&
 a_1f(a_2,\ldots,a_{r+1})\\
 &+\sum_{j=1}^r(-1)^j f(a_1,\ldots,a_ja_{j+1},\ldots,a_{r+1})\\
 &+(-1)^{r+1}f(a_1,\ldots,a_r)a_{r+1}.
\end{align*}
Thus $HH^0(A)=Z(A)$ and
$HH^1(A)=\operatorname{Der}(A)/\operatorname{Inn}(A)$.
The cup product uses multiplication of cochain outputs; on our
small resolution it is implemented by the explicit shuffle diagonal.

The quiver computation is over a field $\Bbbk$, whereas the block
and lattice-model interpretation uses $\mathbb C$. For a general
bimodule $M$, the same resolution computes $HH^*(A,M)$, but the
three families of nonzero regular-coefficient cochain coordinates
need not suffice; see Appendix~\ref{app:coefficients}.
In the quiver sections, paths multiply from left to right.
The integer base index $q$ of $\mathrm C_{q;u,v}$ is distinct from
the complex root parameter $q$ in the TL sections. A formal
chain with $n$ tensor slots has homological degree $n$; it is not
the algebra product of its arrow labels.

\ResearchPart{An explicit minimal Anick resolution}{This part works entirely with the bound quiver algebra $\Lambda(m)$.
It presents a finite four-term bimodule resolution, verifies its
exactness by a relative-bar contraction, and computes explicit Hochschild
cochains, cup products, and the Gerstenhaber bracket. The derivation of
the resolution by Morse paths is not required for the proof presented here.}

	\section{A quadratic Gr\"obner--Shirshov basis}
	\label{sec:GS}
	
	Let $\Bbbk$ be a field and fix $m\geq2$.  Let $Q_m$ be the quiver
	with vertices $1,\ldots,m$ and arrows
	\[
	\mathbf{a}_i:i+1\longrightarrow i,
	\qquad \mathbf{b}_i:i\longrightarrow i+1
	\qquad(1\leq i\leq m-1).
	\]
	Let $I_m$ be the following ideal of the path algebra $\Bbbk Q_m$:
	\begin{equation}\label{eq:ideal}
		I_m=\left\langle
		\mathbf{a}_{i+1}\mathbf{a}_i,\ \mathbf{b}_i\mathbf{b}_{i+1},\ \mathbf{a}_i\mathbf{b}_i-\mathbf{b}_{i+1}\mathbf{a}_{i+1}
		\ (1\leq i\leq m-2),\quad \mathbf{a}_{m-1}\mathbf{b}_{m-1}
		\right\rangle.
	\end{equation}
	For each $m$ we write $\Lambda(m)=\Bbbk Q_m/I_m$.  Whenever $m$ is
	fixed, we abbreviate $\Lambda=\Lambda(m)$.  Since the defining
	relations are homogeneous of path length two, $\Lambda$ has its natural
	grading
	\[
	\Lambda=\bigoplus_{d\geq0}\Lambda_d,
	\qquad \Lambda_0=\bigoplus_{i=1}^m\Bbbk e_i.
	\]
	The unit is $1=e_1+\cdots+e_m$.  Path relations are built into
	$\Bbbk Q_m$; we do not introduce an additional independent unit
	alongside the vertices.
	
	Order paths first by length and then lexicographically, using
	\begin{equation}\label{eq:order}
		\mathbf{a}_1>\cdots>\mathbf{a}_{m-1}>\mathbf{b}_1>\cdots>\mathbf{b}_{m-1}.
	\end{equation}
	The order between length-zero paths may be chosen arbitrarily.
	This is a well-founded order compatible with multiplication whenever the
	products are defined.  The reduction rules are
	\begin{equation}\label{eq:rules}
		\begin{aligned}
			\mathbf{a}_{i+1}\mathbf{a}_i&\longrightarrow0,&
			\mathbf{b}_i\mathbf{b}_{i+1}&\longrightarrow0
			&&(1\leq i\leq m-2),\\
			\mathbf{a}_i\mathbf{b}_i&\longrightarrow \mathbf{b}_{i+1}\mathbf{a}_{i+1}
			&&&&(1\leq i\leq m-2),\\
			\mathbf{a}_{m-1}\mathbf{b}_{m-1}&\longrightarrow0.
		\end{aligned}
	\end{equation}
	
	\begin{theorem}\label{thm:GS}
		The displayed generators of $I_m$, taken in monic form with respect to
		\eqref{eq:order}, are a reduced Gr\"obner--Shirshov basis.  Their leading
		paths are
		\begin{equation}\label{eq:tips}
			\mathcal T=
			\{\mathbf{a}_{i+1}\mathbf{a}_i,\mathbf{b}_i\mathbf{b}_{i+1}:1\leq i\leq m-2\}
			\cup\{\mathbf{a}_i\mathbf{b}_i:1\leq i\leq m-1\}.
		\end{equation}
		In particular, the completion introduces no cubic leading paths.
	\end{theorem}

    \begin{remark}[Path-algebra versus free-algebra leading words]
        The set $\mathcal T$ is the complete set of leading
        \emph{paths} for the Gr\"obner--Shirshov basis of the ideal
        $I_m\subseteq\Bbbk Q_m$.  It should not be confused with the
        leading words of a presentation in a free associative algebra
        with independent symbols for all vertices and arrows.
        Such a presentation must also encode vertex orthogonality and
        endpoint compatibility; for example, $\mathbf a_i\mathbf a_k=0$
        for $i\ne k+1$, and $e_s\mathbf a_k=0$ whenever $s\ne k+1$.
        These products are already zero in the path algebra, so they are
        not additional leading paths in $\mathcal T$.  Moreover,
        noncomposable tensors disappear upon balancing over $\Lambda_0$
        in the relative bar complex, before the Morse matching is
        applied (Section~\ref{sec:matching}).
    \end{remark}

	\begin{proof}
		Each nonzero right-hand side of \eqref{eq:rules} is smaller than its
		left-hand side, so reductions terminate.  Distinct leading paths have
		length two, and hence there are no nontrivial inclusion ambiguities.
		The complete list of overlap paths is
		\begin{equation}\label{eq:overlap-list}
			\begin{array}{c|c}
				\text{overlap path}&\text{range}\\\hline
				\mathbf{a}_{i+2}\mathbf{a}_{i+1}\mathbf{a}_i&1\leq i\leq m-3\\
				\mathbf{b}_i\mathbf{b}_{i+1}\mathbf{b}_{i+2}&1\leq i\leq m-3\\
				\mathbf{a}_{i+1}\mathbf{a}_i\mathbf{b}_i&1\leq i\leq m-2\\
				\mathbf{a}_i\mathbf{b}_i\mathbf{b}_{i+1}&1\leq i\leq m-2.
			\end{array}
		\end{equation}
		To see completeness, a leading pair of type $aa$ may be followed in an
		overlap only by $aa$ or $ab$, a pair of type $ab$ only by $bb$, and a pair
		of type $bb$ only by $bb$.  The path indices then force exactly the words
		in \eqref{eq:overlap-list}.
		
		Both reductions of an $aaa$ or a $bbb$ overlap are immediately zero.
		For the third type, one branch is
		\[
		(\mathbf{a}_{i+1}\mathbf{a}_i)\mathbf{b}_i\longrightarrow0.
		\]
		The other is
		\[
		\mathbf{a}_{i+1}(\mathbf{a}_i\mathbf{b}_i)
		\longrightarrow \mathbf{a}_{i+1}\mathbf{b}_{i+1}\mathbf{a}_{i+1}.
		\]
		If $i\leq m-3$, this reduces further as
		\[
		\mathbf{a}_{i+1}\mathbf{b}_{i+1}\mathbf{a}_{i+1}
		\longrightarrow \mathbf{b}_{i+2}\mathbf{a}_{i+2}\mathbf{a}_{i+1}
		\longrightarrow0.
		\]
		If $i=m-2$, the boundary relation $\mathbf{a}_{m-1}\mathbf{b}_{m-1}=0$ gives zero instead.
		For the fourth type, the branch using $\mathbf{b}_i\mathbf{b}_{i+1}=0$ is zero, while
		\[
		(\mathbf{a}_i\mathbf{b}_i)\mathbf{b}_{i+1}
		\longrightarrow \mathbf{b}_{i+1}\mathbf{a}_{i+1}\mathbf{b}_{i+1}
		\longrightarrow \mathbf{b}_{i+1}\mathbf{b}_{i+2}\mathbf{a}_{i+2}
		\longrightarrow0
		\]
		when $i\leq m-3$.  For $i=m-2$, the middle word contains the boundary
		relation and is again zero.  Thus every intersection composition is trivial modulo the defining
		relations.  Since all leading paths have length two, there are no
		inclusion compositions.  By the Composition--Diamond Lemma
		\cite{Bokut,Bergman}, the irreducible paths (those containing none
		of the leading paths in \eqref{eq:tips} as a subpath) form a
		$\Bbbk$-basis of $\Lambda$.
		
		Finally, no leading path contains another one, and each nonzero
		right-hand side $\mathbf{b}_{i+1}\mathbf{a}_{i+1}$ contains no leading path.  Thus the basis
		is reduced.  The argument includes $m=2$, where only $\mathbf{a}_1\mathbf{b}_1=0$ remains.
	\end{proof}
	
	\begin{corollary}\label{cor:basis}
		Write $\boldsymbol{\omega}_i=\mathbf{b}_i\mathbf{a}_i$ for $1\leq i\leq m-1$.  A normal-path basis is
		\begin{equation}\label{eq:basis}
			\mathcal B=
			\{e_i:1\leq i\leq m\}
			\cup\{\mathbf{a}_i,\mathbf{b}_i:1\leq i\leq m-1\}
			\cup\{\boldsymbol{\omega}_i:1\leq i\leq m-1\}.
		\end{equation}
		Consequently $\dim_{\Bbbk}\Lambda=4m-3$, and the grading is
		concentrated in degrees $0,1,2$:
		\[
		\begin{aligned}
			\Lambda_0&=\bigoplus_{i=1}^m\Bbbk e_i,\\
			\Lambda_1&=\operatorname{span}_{\Bbbk}\{\mathbf{a}_i,\mathbf{b}_i:1\leq i<m\},\\
			\Lambda_2&=\operatorname{span}_{\Bbbk}\{\boldsymbol{\omega}_i:1\leq i<m\},\\
			\Lambda_d&=0\qquad(d\geq3).
		\end{aligned}
		\]
		Moreover, put $J=\Lambda_{>0}=\bigoplus_{d\ge1}\Lambda_d\,=\Lambda_1\oplus\Lambda_2$; it
		equals $\operatorname{rad}\Lambda$, and $J^3=0$.
	\end{corollary}
	
	\begin{proof}
		The only irreducible paths of length two are $\mathbf{b}_i\mathbf{a}_i$.  An irreducible
		path of length three would have two consecutive irreducible pairs, both
		of type $ba$, which is impossible.  Longer irreducible paths therefore do
		not exist either.  Theorem~\ref{thm:GS} gives \eqref{eq:basis}.
		The span of the positive paths is a nilpotent ideal and the quotient by it
		is $\Lambda_0\cong\Bbbk^m$; it is therefore the radical.
	\end{proof}
	
	For subsequent computations, the complete list of nonzero products of
	two arrows is
	\begin{equation}\label{eq:arrow-products}
		\mathbf{b}_i\mathbf{a}_i=\boldsymbol{\omega}_i\quad(1\leq i\leq m-1),\qquad
		\mathbf{a}_i\mathbf{b}_i=\boldsymbol{\omega}_{i+1}\quad(1\leq i\leq m-2).
	\end{equation}
	Every product of a $\boldsymbol{\omega}_i$ with an arrow, on either side, and every product
	of two $\boldsymbol{\omega}_i$'s is zero.  Idempotents act according to endpoints.

\section{The explicit Anick resolution and its contraction}\label{sec:differential}
The algebra and its normal basis are fixed by \eqref{eq:ideal} and
\eqref{eq:basis}. We now state the resolution, then verify that it is
one. Its original derivation by directed movement of a marked normal
pair is not reproduced.

\subsection{Modules and the four-term map}
For
\begin{equation}\label{eq:indices}
 k\ge1,\quad u,v\ge0,\quad u+v=n,\quad k+\max(u,v)\le m,
\end{equation}
let $\mathrm C_{k;u,v}$ have endpoints $k+u,k+v$. Equivalently, using
$q$ for the base vertex within the quiver sections, write
\begin{equation}\label{eq:endpoints}
 s(\mathrm C_{q;u,v})=q+u,\qquad r(\mathrm C_{q;u,v})=q+v.
\end{equation}
The letter $q$ here is an integer vertex index, not the root parameter
of the later TL chapters. Put
\begin{equation}\label{eq:Pn}
 \mathsf A_n=\bigoplus_{\substack{q\ge1,\ u,v\ge0\\
 u+v=n,\ q+\max(u,v)\le m}}
 \Lambda e_{q+u}\otimes_{\Bbbk}\Bbbk\mathrm C_{q;u,v}
 \otimes_{\Bbbk}e_{q+v}\Lambda.
\end{equation}
The distinguished tensor attached to a generator is
\begin{equation}\label{eq:chains}
 \mathsf{[\mathbf a_{q+u-1}|\cdots|\mathbf a_q|
 \mathbf b_q|\cdots|\mathbf b_{q+v-1}]}.
\end{equation}
Bars are tensor separators, not multiplications. In particular, a
nonzero formal cell can have zero product in $\Lambda$.
In degree zero the tensor is empty at vertex $q$ and
$\mu(\mathrm C_{q;0,0})=e_q$. Define, bimodule-linearly,
\begin{equation}\label{eq:main-differential}
 \boxed{\begin{aligned}
 \mathsf d_n(\mathrm C_{q;u,v})={}&
 \mathbf a_{q+u-1}\mathrm C_{q;u-1,v}
 +(-1)^{u+v}\mathrm C_{q;u,v-1}\mathbf b_{q+v-1}\\
 &+(-1)^u\mathbf b_{q+u}\mathrm C_{q+1;u,v-1}
 +(-1)^u\mathrm C_{q+1;u-1,v}\mathbf a_{q+v}.
 \end{aligned}}
\end{equation}
Every inadmissible generator or nonexistent arrow means zero.
The generator count is
\begin{equation}\label{eq:count}
 g_n=\sum_{u=0}^n\max\{0,m-\max(u,n-u)\}.
\end{equation}
At this point square-zero and exactness have not been assumed.

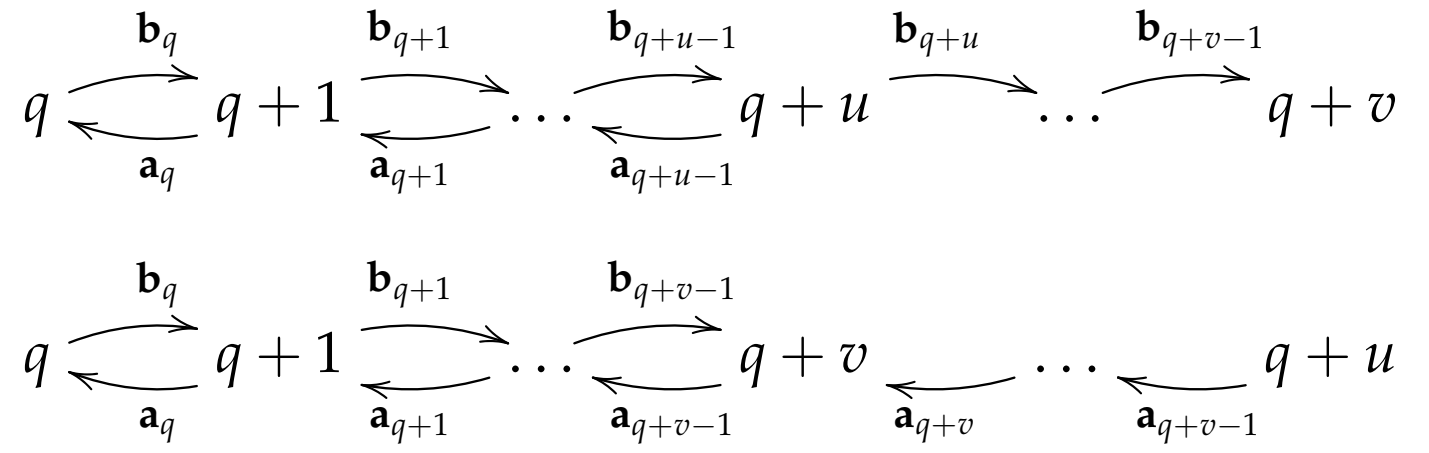
\begin{figure}
	\[
	\xymatrix{
		q  \ar@/^/[r]^{\mathbf{b}_q}   &\ar@/^/[l]^{\mathbf{a}_q} q+1 \ar@/^/[r]^{\mathbf{b}_{q+1}}  & \ar@/^/[l]^{\mathbf{a}_{q+1}}  \ldots \ar@/^/[r]^{\mathbf{b}_{q+u-1}}  & q+u \ar@/^/[l]^{\mathbf{a}_{q+u-1}} \ar@/^/[r]^{\mathbf{b}_{q+u}}  & \ldots \ar@/^/[r]^{\mathbf{b}_{q+v-1}} & q+v & \mbox{$u\le v$}
	}
	\]
	\[
	\xymatrix{
		q  \ar@/^/[r]^{\mathbf{b}_q}   &\ar@/^/[l]^{\mathbf{a}_q} q+1 \ar@/^/[r]^{\mathbf{b}_{q+1}}  & \ar@/^/[l]^{\mathbf{a}_{q+1}}  \ldots \ar@/^/[r]^{\mathbf{b}_{q+v-1}}  & q+v \ar@/^/[l]^{\mathbf{a}_{q+v-1}}  & \ar@/^/[l]^{\mathbf{a}_{q+v}}  \ldots &  \ar@/^/[l]^{\mathbf{a}_{q+v-1}} q+u & \mbox{$u\ge v$}
	}
	\]
	\caption{In both cases, the chain \(\mathrm C_{q;u,v}\) starts at \(q+u\), follows the \(a\)-arrows to \(q\), and then follows the \(b\)-arrows to \(q+v\). Thus its initial and terminal vertices are \(q+u\) and \(q+v\), respectively.
	}
\end{figure}

\subsection{Relative bar data}\label{sec:normalization}
Write $\Lambda=\Lambda_0\oplus J$ and $x=x^{(0)}+x^+$. Set
\begin{equation}\label{eq:relative-bar}
 \mathsf B_n^{\Lambda_0}=
 \Lambda\otimes_{\Lambda_0}J^{\otimes_{\Lambda_0}n}
 \otimes_{\Lambda_0}\Lambda,
 \qquad \mathsf B_0^{\Lambda_0}=\Lambda\otimes_{\Lambda_0}\Lambda.
\end{equation}
The internal normal-path slots are composable. With degree-zero
augmentation given by multiplication, the differential is
\begin{equation}\label{eq:bar-d}
 \begin{aligned}
 \mathsf d_{\mathsf B}(x_0\mathsf{[x_1|\cdots|x_n]}x_{n+1})={}&
 x_0x_1\mathsf{[x_2|\cdots|x_n]}x_{n+1}\\
 &+\sum_{j=1}^{n-1}(-1)^j
 x_0\mathsf{[x_1|\cdots|x_jx_{j+1}|\cdots|x_n]}x_{n+1}\\
 &+(-1)^nx_0\mathsf{[x_1|\cdots|x_{n-1}]}x_nx_{n+1}.
 \end{aligned}
\end{equation}
Products are expanded in the normal basis; zero products contribute zero.

\begin{lemma}\label{lem:relative-exact}
The augmented relative bar complex is a projective bimodule resolution.
It has the right-$\Lambda$-linear, left-$\Lambda_0$-linear contraction
\begin{equation}\label{eq:bar-contraction}
 \mathsf s_{-1}(x)=1\otimes_{\Lambda_0}x,\qquad
 \mathsf s_n(x\mathsf{[w_1|\cdots|w_n]}y)
 =1\mathsf{[x^+|w_1|\cdots|w_n]}y.
\end{equation}
\end{lemma}
\begin{proof}
Balancing is respected because $x\mapsto x^+$ is $\Lambda_0$-bilinear.
Adjacent terms cancel in $\mathsf d_{\mathsf B}\mathsf s+
\mathsf s\mathsf d_{\mathsf B}$. In positive degree the remainder is
\[
 x^+\mathsf{[w_1|\cdots|w_n]}y+
 1\mathsf{[x^{(0)}w_1|w_2|\cdots|w_n]}y
 =x\mathsf{[w_1|\cdots|w_n]}y.
\]
In degree zero the same identity is
$x^+\otimes y-1\otimes x^+y+1\otimes xy=x\otimes y$, and
$\mu_{\mathsf B}\mathsf s_{-1}=1$.
Each normal cell with endpoints $s,t$ gives the summand
$\Lambda e_s\otimes_\Bbbk e_t\Lambda\cong
\Lambda^e(e_s\otimes e_t)$, proving projectivity.
\end{proof}
For comparison with ordinary Hochschild cochains, the absolute bar map
replaces each internal factor $x_j$ by $x_j^+$ and balances over
$\Lambda_0$. The terms with vertex-algebra entries cancel in adjacent
positions, so this is a chain map lifting the identity. A section on
composable positive normal cells inserts their endpoint idempotents;
in degree zero it sends $x\otimes_{\Lambda_0}y$ to
$\sum_i xe_i\otimes_\Bbbk e_i y$. This is the usual relative comparison
for the split separable vertex algebra \cite{Cibils,CLMS}.

\subsection{The shuffle comparison map}

Recall that $\Lambda_1\subseteq J$ is the $\Lambda_0$-bimodule
spanned by the arrows. No additional arrow space is needed. For an
admissible triple $(q,u,v)$ satisfying \eqref{eq:indices}, put
$n=u+v$ and let $\mathcal W_{q;u,v}$ be the set of composable arrow
sequences with $u$ descending steps and $v$ ascending steps, starting
at vertex $q+u$ and remaining in $\{1,\ldots,m\}$. A descending step
from $h$ is $\mathbf a_{h-1}$, an ascending step from $h$ is
$\mathbf b_h$, and every sequence ends at $q+v$. Such a sequence is
regarded as an element of $\Lambda_1^{\otimes_{\Lambda_0}n}$;
its product in $\Lambda$ need not be nonzero.

For $w=(w_1,\ldots,w_n)\in\mathcal W_{q;u,v}$, let
\[
 \operatorname{inv}(w)
 =\#\{(i,j):i<j,\ \text{step $i$ is ascending and
 step $j$ is descending}\}.
\]
Write
\[
 \mathsf{[w]}_{\Lambda_0}
 =e_{q+u}\otimes_{\Lambda_0}w_1\otimes_{\Lambda_0}\cdots
 \otimes_{\Lambda_0}w_n\otimes_{\Lambda_0}e_{q+v}
 \ \in\mathsf B_n^{\Lambda_0}.
\]
The distinguished word, with all descending steps preceding all
ascending steps, is the critical word defining $\mathrm C_{q;u,v}$.
It has inversion number zero. For $u=v=0$, the unique word is the
empty sequence at $q$.

\begin{proposition}[Shuffle comparison]\label{cup:prop:comparison}
For every $n\geq0$, the formula
\begin{equation}\label{cup:eq:kappa}
 \boxed{\displaystyle
 \kappa_n(\mathrm C_{q;u,v})
 =\sum_{w\in\mathcal W_{q;u,v}}
 (-1)^{\operatorname{inv}(w)}\mathsf{[w]}_{\Lambda_0}}
\end{equation}
defines an injective, degreewise split $\Lambda$-bimodule map
\[
 \kappa_n:\mathsf A_n\longrightarrow\mathsf B_n^{\Lambda_0}.
\]
These maps satisfy the direct identity
\begin{equation}\label{cup:eq:kappa-chain}
 \mathsf d_{\mathsf B}\kappa=\kappa\mathsf d,
 \qquad \kappa_0=\operatorname{id}_{\mathsf A_0}.
\end{equation}
At this stage we have not used exactness of $\mathsf A$; the identity is
a verification of the displayed four-term map against the bar differential.
\end{proposition}

\begin{proof}
All summands in \eqref{cup:eq:kappa} have initial vertex $q+u$
and terminal vertex $q+v$, so the formula extends bimodule-linearly.
There is an explicit retraction on graded bimodules: for $n\geq1$,
let $\mathsf r_n$ send a basis bar cell to $\mathrm C_{q;u,v}$ when
its internal blocks are exactly the critical arrow sequence
$\mathbf a_{q+u-1}|\cdots|\mathbf a_q|
\mathbf b_q|\cdots|\mathbf b_{q+v-1}$, and to zero otherwise.
Extend this map bimodule-linearly, and put
$\mathsf r_0=\operatorname{id}$. Every shuffle sum contains its
distinguished critical tensor with coefficient $1$, while no other
shuffle sum contains that tensor. Hence
\begin{equation}\label{cup:eq:kappa-retraction}
 \boxed{\mathsf r_n\kappa_n=\operatorname{id}_{\mathsf A_n}.}
\end{equation}
This proves injectivity and splitting in each degree. The retraction
$\mathsf r$ is not asserted to be a chain map and is not the Morse
projection $p$.

For the chain identity, separate internal and exterior terms of the
bar differential. Internal multiplications of two $a$-arrows or two
$b$-arrows are zero. At an interior vertex $2\leq h\leq m-1$,
the two remaining products are equal:
\[
 \mathbf a_{h-1}\mathbf b_{h-1}
 =\boldsymbol\omega_h
 =\mathbf b_h\mathbf a_h.
\]
Their contributions are represented by the local diagram
\[
\begin{array}{ccc}
 \mathsf{[\cdots|\mathbf a_{h-1}|\mathbf b_{h-1}|\cdots]}
 &&\mathsf{[\cdots|\mathbf b_h|\mathbf a_h|\cdots]}\\[-1mm]
 &\searrow\qquad\swarrow&\\[-1mm]
 &\mathsf{[\cdots|\boldsymbol\omega_h|\cdots]}&
\end{array}
\]
The two upper words differ by one adjacent interchange and have
opposite shuffle signs. Since the merged pair occupies the same bar
position $j$, their combined contribution is
\begin{equation}\label{cup:eq:internal-cancellation}
 (-1)^{j+\operatorname{inv}(w)}
 \mathsf{[\cdots|\boldsymbol\omega_h|\cdots]}
 +(-1)^{j+\operatorname{inv}(w)+1}
 \mathsf{[\cdots|\boldsymbol\omega_h|\cdots]}=0.
\end{equation}
Both words remain admissible because the interchange has the same
initial and terminal vertex. At $h=m$, the unmatched product
$\mathbf a_{m-1}\mathbf b_{m-1}$ is zero. The other possible
unmatched boundary, a $ba$ pair at vertex $1$, never occurs: if $d$
downward and $a$ upward steps precede such a pair, then $d\leq u-1$
and the current vertex $h=q+u-d+a\geq q+1\geq2$. Thus all nonzero
internal terms cancel.

It remains to collect the four kinds of exterior terms. With
$n=u+v$, conditioning on the first or last step gives
\begin{equation}\label{cup:eq:four-exterior}
\begin{array}{c|l}
\text{first step $a$}
 &\mathbf a_{q+u-1}\cdot\kappa(\mathrm C_{q;u-1,v})\\
\text{first step $b$}
 &(-1)^u\mathbf b_{q+u}\cdot\kappa(\mathrm C_{q+1;u,v-1})\\
\text{last step $b$}
 &(-1)^n\kappa(\mathrm C_{q;u,v-1})\cdot\mathbf b_{q+v-1}\\
\text{last step $a$}
 &(-1)^{n+v}\kappa(\mathrm C_{q+1;u-1,v})\cdot\mathbf a_{q+v}.
\end{array}
\end{equation}
Here the extra $u$ inversions in the second row arise from moving its
initial $b$ past the remaining $u$ descending steps. In the last row,
the final $a$ follows all $v$ ascending steps and hence contributes
$v$ additional inversions. Because $n+v\equiv u\pmod2$, the four
rows agree term for term with \eqref{eq:main-differential}, including
its boundary and inadmissible-index conventions. Therefore
$\mathsf d_{\mathsf B}\kappa=\kappa\mathsf d$.
In degree zero the empty word gives $\kappa_0=\operatorname{id}$,
compatible with the augmentation $\mu$.
\end{proof}

For example,
\[
 \kappa_2(\mathrm C_{i;1,1})
 =\mathsf{[\mathbf a_i|\mathbf b_i]}_{\Lambda_0}
 -\mathsf{[\mathbf b_{i+1}|\mathbf a_{i+1}]}_{\Lambda_0}
 \qquad(1\leq i\leq m-2),
\]
whereas at $i=m-1$ only the first term occurs. In particular, a
formal obstruction generator need not be represented by a single
bar tensor.

Since $\kappa$ is injective, the chain identity just checked gives
\[
 \kappa\mathsf d^2=\mathsf d_{\mathsf B}^{\,2}\kappa=0,
 \qquad\text{hence}\qquad\mathsf d^2=0.
\]
This establishes square-zero independently of exactness.

\subsection{The finite projection and the augmented contraction}\label{sec:matching}
We specify the comparison projection used below, rather than leaving the
contracting homotopy as an existence statement. Scan a composable normal
bar cell from the left. Its first defect is a loop $\boldsymbol\omega_i$
or a normal pair $\mathbf b_i|\mathbf a_i$. If there is no defect,
the cell is a distinguished tensor \eqref{eq:chains}. If the first defect
is a loop in position $j$, denote the cell by $L$ and let $U$ be the
cell obtained by splitting this loop into $\mathbf b_i|\mathbf a_i$.
Define, bimodule-linearly,
\begin{equation}\label{eq:p-rec}
 \begin{aligned}
 p(C)&=\mathrm C_{q;u,v}&&\text{for a distinguished tensor }C,\\
 p(U)&=0&&\text{if the first defect is a normal pair},\\
 p(L)&=-(-1)^j p(\mathsf d_{\mathsf B}U-(-1)^jL)
 &&\text{if the first defect is a loop at }j.
 \end{aligned}
\end{equation}
Put $p_0=\operatorname{id}$ under the vertex identification
$\mathsf B_0^{\Lambda_0}=\mathsf A_0$.

\begin{lemma}\label{lem:matching}\label{lem:pkappa}
The recursion is finite, defines an augmented bimodule chain map
$p:\mathsf B^{\Lambda_0}\to\mathsf A$ for \eqref{eq:main-differential},
and satisfies $p\kappa=\operatorname{id}_{\mathsf A}$.
\end{lemma}
\begin{proof}
Pair $U$ to $L$ using its internal multiplication coefficient $(-1)^j$.
This is an isomorphism of equal endpoint projectives. The first-defect
operations are inverse, so no cell is paired twice. To check finite
acyclicity expand loops formally as $b_i a_i$, and order cells by
$(W,\varepsilon,\ell)$: $W$ is the expanded word with the path order
\eqref{eq:order}; $\varepsilon=1$ for a first-loop defect and $0$
otherwise; $\ell$ is the number of blocks. A reversed matching edge
lowers $\varepsilon$. An exterior multiplication shortens $W$;
$a_i b_i\mapsto b_{i+1}a_{i+1}$ lowers $W$; any remaining nonzero
unmatched merger preserves the first defect and lowers $\ell$.
Only finitely many words and decompositions are reachable. Thus the
path-sum comparison theorem applies to the endpoint projectives
\cite[Appendix~B]{JollenbeckWelker}, \cite[Theorem~3.3]{CLZ}.
It gives a differential $\partial$ on the distinguished generators and
a projection with \eqref{eq:p-rec}, satisfying
$\partial p=p\mathsf d_{\mathsf B}$.
Every nondistinguished arrow shuffle has a first normal pair $b_i|a_i$
and is killed by $p$; the distinguished word occurs with coefficient one.
Hence $p\kappa=1$. Using the independently proved shuffle identity gives
\[
 \partial=\partial p\kappa=p\mathsf d_{\mathsf B}\kappa
 =p\kappa\mathsf d=\mathsf d.
\]
This identifies the map needed for the contraction without using
exactness of the proposed complex or classifying the transferred paths.
\end{proof}

\begin{theorem}[Explicit resolution and contraction]\label{thm:resolution}
The modules \eqref{eq:Pn}, differential \eqref{eq:main-differential},
and augmentation $\mu$ form a projective bimodule resolution.
Its augmented underlying right-module complex has contraction
\begin{equation}\label{eq:explicit-contraction}
 \boxed{\mathsf h_{-1}=p_0\mathsf s_{-1},\qquad
 \mathsf h_n=p_{n+1}\mathsf s_n\kappa_n\quad(n\ge0).}
\end{equation}
Thus $\mu\mathsf h_{-1}=1$ and, with $\mathsf d_0=\mu$,
$\mathsf d_{n+1}\mathsf h_n+\mathsf h_{n-1}\mathsf d_n=1$ for $n\ge0$.
\end{theorem}
\begin{proof}
Projectivity follows from the endpoint summands. Square-zero was
proved with the injective shuffle map. Since $p,\kappa$ are chain maps,
\[
 \mathsf d\mathsf h+\mathsf h\mathsf d
 =p(\mathsf d_{\mathsf B}\mathsf s+
       \mathsf s\mathsf d_{\mathsf B})\kappa=p\kappa=1.
\]
In degree zero this includes $\mathsf h_{-1}\mu$; on the augmented
target, $\mu p_0\mathsf s_{-1}=1$. This proves exactness.
\end{proof}
For endpoint-compatible exterior coefficients, the full formula is
\begin{equation}\label{eq:h-expanded}
 \mathsf h_n(x\mathrm C_{q;u,v}y)=
 \sum_{w\in\mathcal W_{q;u,v}}(-1)^{\operatorname{inv}(w)}
 p_{n+1}(1\mathsf{[x^+|w_1|\cdots|w_n]}y).
\end{equation}
It is right-$\Lambda$-linear and left-$\Lambda_0$-linear, not generally
bimodule-linear. In particular it vanishes on bare critical generators
but not on all their left multiples. For $m=2$, write
$\alpha=\mathrm C_{1;1,0}$, $\gamma=\mathrm C_{1;0,1}$,
$\xi=\mathrm C_{1;1,1}$, and $a=\mathbf a_1$, $b=\mathbf b_1$.
Then $\mathsf d_2\xi=a\gamma+\alpha b$, whereas
$\mathsf h_1(a\gamma)=\xi$ and $\mathsf h_1(\alpha b)=0$.
Thus $\mathsf h_1\mathsf d_2\xi=\xi$, illustrating why exterior
coefficients cannot be pulled through the bar contraction.

\begin{remark}[Morse origin]
The generators are Anick chains for the quadratic leading paths.
The formula was found by Morse reduction and directed movement of a
marked normal pair, in the spirit of \cite{Lopatkin}. Here this history
is replaced by verification of the displayed complex. The finite
projection above is retained solely because it defines the contraction
and justifies its identities.
\end{remark}

\subsection{Minimality and termination}
	\label{sec:minimality}
	
	\begin{theorem}\label{thm:minimal}
		The resolution in Theorem~\ref{thm:resolution} is minimal and has length
		$2m-2$:
		\begin{equation}\label{eq:finite-resolution}
			0\longrightarrow \mathsf A_{2m-2}\longrightarrow\cdots\longrightarrow
			\mathsf A_1\longrightarrow \mathsf A_0\xrightarrow{\mu}\Lambda\longrightarrow0.
		\end{equation}
		Moreover,
		\[
		\mathsf A_{2m-2}\cong\Lambda e_m\otimes_{\Bbbk} e_m\Lambda
		\]
		with its generator labelled by
		$\mathsf{[\mathbf{a}_{m-1}|\cdots|\mathbf{a}_1|\mathbf{b}_1|\cdots|\mathbf{b}_{m-1}]}$.
	\end{theorem}
	
	\begin{proof}
		Every nonzero exterior coefficient in \eqref{eq:main-differential} is an
		arrow.  Therefore
		\[
		\mathsf d_n(\mathsf A_n)\subseteq J\mathsf A_{n-1}+\mathsf A_{n-1}J.
		\]
		The ideal
		$J\otimes_{\Bbbk}\Lambda^{\mathrm{op}}+
		\Lambda\otimes_{\Bbbk}J^{\mathrm{op}}$ is the radical of $\Lambda^e$:
		it is nilpotent and its quotient is $(\Lambda_0)^e\cong \Bbbk^{m^2}$.
		Thus the displayed inclusion is exactly the radical criterion for a
		minimal projective $\Lambda^e$-resolution.  Equivalently, it also follows
		from homogeneity: each critical generator in degree $n$ has internal
		arrow degree $n$, and a differential to degree $n-1$ has exterior degree
		one.
		
		The index inequalities \eqref{eq:indices} imply
		\[
		n=u+v\leq2\max(u,v)\leq2(m-q)\leq2m-2.
		\]
		Equality forces $q=1$ and $u=v=m-1$, giving the unique top generator and
		the stated endpoint type.  Its bimodule summand is nonzero.
		No assumption on global dimension was used to discard higher terms:
		there simply are no critical generators there.
	\end{proof}
	
	\begin{corollary}\label{cor:gldim}
		The left and right global dimensions of $\Lambda$ are both $2m-2$.
	\end{corollary}
	
	\begin{proof}
		Put $d=2m-2$.  Viewed as a complex of right $\Lambda$-modules, the
		augmented resolution \eqref{eq:finite-resolution} is split exact:
		$\Lambda$ is projective and its successive kernels split inductively.
		Tensoring it over $\Lambda$ with any left module $M$ is therefore exact.
		Each resulting term is a direct sum of modules of the form
		\[
		\Lambda e_s\otimes_{\Bbbk} e_tM,
		\]
		which are projective as left $\Lambda$-modules.  Hence every left module
		has projective dimension at most $d$.
		
		For the reverse inequality, tensor the resolution first on the right
		with $\Lambda_0$ and then on the left with $\Lambda_0$.  The first operation gives a
		projective left-module resolution of $\Lambda_0$.  The differentials become
		zero after the second operation, since every term of
		\eqref{eq:main-differential} has an exterior arrow.  Consequently
		\[
		\operatorname{Tor}^{\Lambda}_d(\Lambda_0,\Lambda_0)
		\cong \Lambda_0\otimes_\Lambda \mathsf A_d\otimes_\Lambda \Lambda_0
		\cong \Bbbk e_m\otimes_{\Bbbk} \Bbbk e_m\ne0.
		\]
		Thus the left global dimension is at least $d$.  The argument for right
		modules is symmetric.  This proves the claimed value as a consequence
		of the constructed resolution, rather than as an input to it.
	\end{proof}

\begin{example}[$m=3$]\label{ex:m3}
		The positive-degree critical generators are
		\[
		\begin{array}{c|l}
			\mathsf A_1&\mathsf{[\mathbf{a}_1]},\ \mathsf{[\mathbf{a}_2]},\ \mathsf{[\mathbf{b}_1]},\ \mathsf{[\mathbf{b}_2]}\\
			\mathsf A_2&\mathsf{[\mathbf{a}_2|\mathbf{a}_1]},\ \mathsf{[\mathbf{b}_1|\mathbf{b}_2]},\ \mathsf{[\mathbf{a}_1|\mathbf{b}_1]},\ \mathsf{[\mathbf{a}_2|\mathbf{b}_2]}\\
			\mathsf A_3&U=\mathsf{[\mathbf{a}_2|\mathbf{a}_1|\mathbf{b}_1]},\quad V=\mathsf{[\mathbf{a}_1|\mathbf{b}_1|\mathbf{b}_2]}\\
			\mathsf A_4&W=\mathsf{[\mathbf{a}_2|\mathbf{a}_1|\mathbf{b}_1|\mathbf{b}_2]}.
		\end{array}
		\]
		There are no generators in degree five.  Put
		\[
		R_1=\mathsf{[\mathbf{a}_2|\mathbf{a}_1]},\quad R_2=\mathsf{[\mathbf{b}_1|\mathbf{b}_2]},\quad
		R_3=\mathsf{[\mathbf{a}_1|\mathbf{b}_1]},\quad R_4=\mathsf{[\mathbf{a}_2|\mathbf{b}_2]}.
		\]
		Then
		\begin{align*}
			\mathsf d_2R_1&=\mathbf{a}_2\mathsf{[\mathbf{a}_1]}+\mathsf{[\mathbf{a}_2]}\mathbf{a}_1,&
			\mathsf d_2R_2&=\mathbf{b}_1\mathsf{[\mathbf{b}_2]}+\mathsf{[\mathbf{b}_1]}\mathbf{b}_2,\\
			\mathsf d_2R_3&=\mathbf{a}_1\mathsf{[\mathbf{b}_1]}+\mathsf{[\mathbf{a}_1]}\mathbf{b}_1-\mathbf{b}_2\mathsf{[\mathbf{a}_2]}-\mathsf{[\mathbf{b}_2]}\mathbf{a}_2,&
			\mathsf d_2R_4&=\mathbf{a}_2\mathsf{[\mathbf{b}_2]}+\mathsf{[\mathbf{a}_2]}\mathbf{b}_2,\\
			\mathsf d_3U&=\mathbf{a}_2R_3-R_1\mathbf{b}_1+R_4\mathbf{a}_2,&
			\mathsf d_3V&=\mathbf{a}_1R_2-\mathbf{b}_2R_4-R_3\mathbf{b}_2,\\
			\mathsf d_4W&=\mathbf{a}_2V+U\mathbf{b}_2.
		\end{align*}
		For example,
		\[
		\kappa(\mathrm U)=\mathsf{[\mathbf{a}_2|\mathbf{a}_1|\mathbf{b}_1]}-\mathsf{[\mathbf{a}_2|\mathbf{b}_2|\mathbf{a}_2]},
		\]
		and applying the bar differential makes the correction $R_4\mathbf{a}_2$ explicit.
		The displayed complex is obtained by Anick matching; no cohomology
		calculation is needed to identify its differentials.
	\end{example}
	
	\section{Parallel paths and the Hochschild cochain spaces}
	\label{sec:cochains}
	
	We now compute Hochschild cohomology on the relative Anick resolution
	constructed above.  No further change of resolution is made.
	We write $\mathrm{C}_{q;u,v}$ for its formal generators and
	$\mathsf C^n$ for the vector space of cochains. The coboundary
	$\mathsf d^n$ below belongs to the Anick cochain complex; later,
	on ordinary bar cochains, we use the notation $\delta$. Put
	$\Lambda^e=\Lambda\otimes_{\Bbbk}\Lambda^{\mathrm{op}}$ and
	\begin{equation}\label{eq:cochain-complex}
		\mathsf C^n=\operatorname{Hom}_{\Lambda^e}(\mathsf A_n,\Lambda),\qquad
		(\mathsf d^n\phi)(\gamma)=\phi(\mathsf d_{n+1}\gamma).
	\end{equation}
	Thus $\mathsf d^{n+1}\mathsf d^n=0$, and
	\begin{equation}\label{eq:HH-definition}
		Z^n=\ker\mathsf d^n,\qquad B^n=\operatorname{im}\mathsf d^{n-1},\qquad
		HH^n(\Lambda)=Z^n/B^n,
	\end{equation}
	where $B^0=0$.  These are the standard cochain constructions on a
	projective bimodule resolution; see, for example, \cite[Section~2]{Volkov}.
	The notation $B^n$ for coboundaries is distinguished from the bar
	resolutions $\mathsf B_n^{\Bbbk},\mathsf B_n^{\Lambda_0}$ by its superscript and the absence of a base
	algebra label.  The spaces $Z^n,B^n$ in what follows refer to this
	particular relative cochain complex, not to all cochains on the original
	unnormalized bar resolution.
	In particular, $\mathsf C^0\cong\bigoplus_{q=1}^m e_q\Lambda e_q$;
	it is not identified with the whole vector space $\Lambda$,
	as happens for the unnormalized bar cochain model.
	
	\subsection{Coefficient functionals and parallel paths}
	
	For a fixed cochain $\phi$ and a critical generator $\gamma$, write
	\begin{equation}\label{eq:functional-expansion}
		\phi(\gamma)=\sum_{w\in\mathcal B}F_w^\phi(\gamma)w.
	\end{equation}
	Here $F_w^\phi(\gamma)\in \Bbbk$ is the coefficient of the normal basis path $w$
	in $\phi(\gamma)$.  We omit the superscript $\phi$ when the cochain is fixed.
	The functions $F_w^\phi$ are scalar coefficient functionals; they are
	not separately required to be bimodule maps to $\Bbbk$.
	
	\begin{definition}\label{def:parallel}
		A normal basis path $w$ and a critical generator $\gamma$ are \emph{parallel}
		if they have the same source and the same target.  This condition concerns
		the endpoints of the formal chain $\gamma$; its product in $\Lambda$ need not
		be nonzero.
	\end{definition}
	
	\begin{lemma}[Parallel-path constraint]\label{lem:parallel}
		Let $\gamma=\mathrm{C}_{q;u,v}$.  Then
		\begin{equation}\label{eq:parallel-value}
			\phi(\gamma)\in e_{q+u}\Lambda e_{q+v},
			\qquad
			F_w^\phi(\gamma)=0\quad\text{unless }w\text{ is parallel to }\gamma.
		\end{equation}
		Conversely, an arbitrary choice of an element of
		$e_{q+u}\Lambda e_{q+v}$ for every critical generator of a fixed degree
		defines a unique cochain of that degree.
	\end{lemma}
	
	\begin{proof}
		The distinguished generator satisfies
		$\gamma=e_{q+u}\gamma e_{q+v}$.  Bimodule-linearity gives
		\[
		\phi(\gamma)=e_{q+u}\phi(\gamma)e_{q+v}
		=\sum_{w\in\mathcal B}F_w^\phi(\gamma)e_{q+u}we_{q+v}.
		\]
		The last multiplication fixes $w$ exactly when its endpoints agree with
		those of $\gamma$, and is zero otherwise.  Compare coefficients in the normal
		basis.  For the converse, on the summand indexed by $\gamma$ extend the
		chosen value by $\phi(x\gamma y)=x\phi(\gamma)y$.  Its endpoint condition makes this
		well-defined, and \eqref{eq:Pn} is a direct sum of these summands.
	\end{proof}
	
	From \eqref{eq:basis} and \eqref{eq:arrow-products}, the endpoint spaces
	are explicitly
	\begin{equation}\label{eq:corner-spaces}
		e_s\Lambda e_t=
		\begin{cases}
			\Bbbk e_s\oplus \Bbbk\boldsymbol{\omega}_s,&s=t<m,\\
			\Bbbk e_m,&s=t=m,\\
			\Bbbk\mathbf{a}_t,&s=t+1,\\
			\Bbbk\mathbf{b}_s,&t=s+1,\\
			0,&|s-t|\geq2.
		\end{cases}
	\end{equation}
	In particular, $\phi(\mathrm{C}_{q;u,v})=0$ whenever $|u-v|\geq2$.
	This is a restriction on \emph{every cochain}, before imposing the
	cocycle condition.  It does not permit us to delete the corresponding
	projective summands from the resolution: their boundaries can contain
	other critical generators.
	
	For reference, if
	$\mathsf d(\gamma)=\sum_t\lambda_t x_t\gamma_ty_t$ is expanded with basis-path exterior
	coefficients, then \eqref{eq:functional-expansion} yields
	\begin{equation}\label{eq:factorization-cochain}
		F_M^{\mathsf d\phi}(\gamma)
		=\sum_t\lambda_t
		\sum_{\substack{w\in\mathcal B\\x_twy_t=M}}
		F_w^\phi(\gamma_t).
	\end{equation}
	Products in this formula are taken in $\Lambda$ and put into normal
	form.  It applies here because a product of normal basis paths is either
	zero or a single normal basis path.  Thus the coefficient equations are
	precisely controlled by factorization equations with parallel outputs.
	
	\subsection{The three families of nonzero cochain coordinates}
	
	Define
	\begin{equation}\label{eq:DUV-definitions}
		\mathrm{D}_q^{(\ell)}=\mathrm{C}_{q;\ell,\ell},\qquad
		\mathrm{U}_q^{(\ell)}=\mathrm{C}_{q;\ell+1,\ell},\qquad
		\mathrm{V}_q^{(\ell)}=\mathrm{C}_{q;\ell,\ell+1}.
	\end{equation}
	In expanded form,
	\begin{align}
		\mathrm{D}_q^{(\ell)}
		&=\mathsf{[\mathbf{a}_{q+\ell-1}|\cdots|\mathbf{a}_q\mid
			\mathbf{b}_q|\cdots|\mathbf{b}_{q+\ell-1}]},\label{eq:D-expanded}\\
		\mathrm{U}_q^{(\ell)}
		&=\mathsf{[\mathbf{a}_{q+\ell}|\cdots|\mathbf{a}_q\mid
			\mathbf{b}_q|\cdots|\mathbf{b}_{q+\ell-1}]},\label{eq:U-expanded}\\
		\mathrm{V}_q^{(\ell)}
		&=\mathsf{[\mathbf{a}_{q+\ell-1}|\cdots|\mathbf{a}_q\mid
			\mathbf{b}_q|\cdots|\mathbf{b}_{q+\ell}]}.\label{eq:V-expanded}
	\end{align}
	Empty arrow blocks are omitted.  Thus
	$\mathrm{D}_q^{(0)}=\mathrm{C}_{q;0,0}$ is a degree-zero vertex generator, while
	$\mathrm{U}_q^{(0)}=\mathsf{[\mathbf{a}_q]}$ and $\mathrm{V}_q^{(0)}=\mathsf{[\mathbf{b}_q]}$.  The degrees and ranges are
	\begin{equation}\label{eq:DUV-table}
		\begin{array}{c|c|c|c|c}
			\text{generator}&\text{degree}&\text{source}&\text{target}&q\text{-range}\\\hline
			\mathrm{D}_q^{(\ell)}&2\ell&q+\ell&q+\ell&1\leq q\leq m-\ell\\
			\mathrm{U}_q^{(\ell)}&2\ell+1&q+\ell+1&q+\ell&1\leq q\leq m-\ell-1\\
			\mathrm{V}_q^{(\ell)}&2\ell+1&q+\ell&q+\ell+1&1\leq q\leq m-\ell-1.
		\end{array}
	\end{equation}
	The family $D$ is used for $0\leq\ell\leq m-1$; the families $\mathrm{U},\mathrm{V}$
	are used for $0\leq\ell\leq m-2$.
	
	\begin{proposition}\label{prop:cochain-coordinates}
		Put $s=m-\ell$.
		An arbitrary cochain $\phi\in \mathsf C^{2\ell}$ is uniquely determined by
		\begin{equation}\label{eq:even-values}
			\phi(\mathrm{D}_q^{(\ell)})=
			\begin{cases}
				E_qe_{q+\ell}+\Omega_q\boldsymbol{\omega}_{q+\ell},&1\leq q<s,\\
				E_se_m,&q=s,
			\end{cases}
		\end{equation}
		and is zero on all other critical generators in degree $2\ell$.
		Thus
		\begin{equation}\label{eq:even-coordinate-space}
			\mathsf C^{2\ell}\cong \Bbbk^s\oplus \Bbbk^{s-1},\qquad
			\phi\longleftrightarrow(\mathbf E,\boldsymbol\Omega),
		\end{equation}
		where $\mathbf E=(E_1,\ldots,E_s)$ and
		$\boldsymbol\Omega=(\Omega_1,\ldots,\Omega_{s-1})$.
		
		An arbitrary cochain $\psi\in \mathsf C^{2\ell+1}$ is uniquely determined by
		\begin{equation}\label{eq:odd-values}
			\psi(\mathrm{U}_q^{(\ell)})=\alpha_q \mathbf{a}_{q+\ell},\qquad
			\psi(\mathrm{V}_q^{(\ell)})=\beta_q \mathbf{b}_{q+\ell}
			\quad(1\leq q\leq s-1),
		\end{equation}
		and is zero on all other critical generators in that degree.  Hence
		\begin{equation}\label{eq:odd-coordinate-space}
			\mathsf C^{2\ell+1}\cong \Bbbk^{s-1}\oplus \Bbbk^{s-1},\qquad
			\psi\longleftrightarrow(\boldsymbol\alpha,\boldsymbol\beta).
		\end{equation}
		All displayed coordinates are independent before imposing $\mathsf d\phi=0$
		or $\mathsf d\psi=0$.
	\end{proposition}
	
	\begin{proof}
		If $u+v=2\ell$, the integer $u-v$ is even.  By
		\eqref{eq:corner-spaces}, a nonzero value is possible only for $u=v=\ell$.
		The corresponding corner has dimension two except at vertex $m$,
		where it has dimension one.  If $u+v=2\ell+1$, the difference is odd,
		and the only possibilities are $(u,v)=(\ell+1,\ell)$ and
		$(\ell,\ell+1)$.  Their corners are the one-dimensional arrow spaces
		in \eqref{eq:corner-spaces}.  Lemma~\ref{lem:parallel} proves that these
		choices extend uniquely and impose no further conditions on a cochain.
	\end{proof}
	
	The coordinates in this proposition have the precise functional meaning
	\begin{equation}\label{eq:coordinate-functionals}
		\begin{aligned}
			E_q&=F_{e_{q+\ell}}^\phi(\mathrm{D}_q^{(\ell)}),&
			\Omega_q&=F_{\boldsymbol{\omega}_{q+\ell}}^\phi(\mathrm{D}_q^{(\ell)}),\\
			\alpha_q&=F_{\mathbf{a}_{q+\ell}}^\psi(\mathrm{U}_q^{(\ell)}),&
			\beta_q&=F_{\mathbf{b}_{q+\ell}}^\psi(\mathrm{V}_q^{(\ell)}).
		\end{aligned}
	\end{equation}
	Superscripts indicating the cohomological degree or the cochain
	will usually be suppressed inside a fixed calculation.
	
	Counting these coordinates gives
	\begin{equation}\label{eq:cochain-dimensions}
		\dim_{\Bbbk}\mathsf C^n=
		\begin{cases}
			2m-1-n,&0\leq n\leq2m-2,\\
			0,&n>2m-2.
		\end{cases}
	\end{equation}
	Thus, as vector spaces, the cochain complex has successive dimensions
	$2m-1,2m-2,\ldots,1$.  The minimality of $\mathsf A_\bullet$ does not imply
	that its Hochschild coboundaries are zero: the target is $\Lambda$,
	not its semisimple quotient $\Lambda_0$.
	
	\section{From chain boundaries to the Hochschild coboundary}
	\label{sec:coboundary}
	
	We first record the complete boundaries of $D,U,V$.  We then apply a
	cochain and use parallel paths to remove precisely those terms whose
	values must vanish.  This order of operations is important: the three
	families in \eqref{eq:DUV-definitions} describe all nonzero cochain
	coordinates but do not, by themselves, describe a subcomplex of
	$\mathsf A_\bullet$.
	
	\subsection{The complete boundaries of \texorpdfstring{$\mathrm{D}$, $\mathrm{U}$, and $\mathrm{V}$}{D, U, and V}}
	
	For $\ell\geq1$ introduce two auxiliary critical generators
	\begin{equation}\label{eq:RS-definitions}
		\mathrm{R}_q^{(\ell)}=\mathrm{C}_{q;\ell+1,\ell-1},\qquad
		\mathrm{S}_q^{(\ell)}=\mathrm{C}_{q;\ell-1,\ell+1}.
	\end{equation}
	They belong to $\mathsf A_{2\ell}$ and have endpoints
	\[
	\mathrm{R}_q^{(\ell)}:q+\ell+1\longrightarrow q+\ell-1,
	\qquad
	\mathrm{S}_q^{(\ell)}:q+\ell-1\longrightarrow q+\ell+1.
	\]
	These symbols denote formal chains. In formulas containing
	$\ell=0$, terms involving $\mathrm{R}^{(0)}$ or $\mathrm{S}^{(0)}$ are omitted.
	As throughout the paper, every term with an inadmissible chain or a
	nonexistent exterior arrow is omitted as well.
	
	\begin{lemma}\label{lem:DUV-boundaries}
		For $\ell\geq1$,
		\begin{equation}\label{eq:dD}
			\begin{aligned}
				\mathsf d_{2\ell}\mathrm{D}_q^{(\ell)}={}&
				\mathbf{a}_{q+\ell-1}\mathrm{V}_q^{(\ell-1)}
				+\mathrm{U}_q^{(\ell-1)}\mathbf{b}_{q+\ell-1}\\
				&+(-1)^\ell \mathbf{b}_{q+\ell}\mathrm{U}_{q+1}^{(\ell-1)}
				+(-1)^\ell \mathrm{V}_{q+1}^{(\ell-1)}\mathbf{a}_{q+\ell}.
			\end{aligned}
		\end{equation}
		For $\ell\geq0$, put $r=q+\ell$.  Then
		\begin{equation}\label{eq:dU}
			\begin{aligned}
				\mathsf d_{2\ell+1}\mathrm{U}_q^{(\ell)}={}&
				\mathbf{a}_r\mathrm{D}_q^{(\ell)}
				+(-1)^{\ell+1}\mathrm{D}_{q+1}^{(\ell)}\mathbf{a}_r\\
				&-\mathrm{R}_q^{(\ell)}\mathbf{b}_{r-1}
				+(-1)^{\ell+1}\mathbf{b}_{r+1}\mathrm{R}_{q+1}^{(\ell)},
			\end{aligned}
		\end{equation}
		and
		\begin{equation}\label{eq:dV}
			\begin{aligned}
				\mathsf d_{2\ell+1}\mathrm{V}_q^{(\ell)}={}&
				-\mathrm{D}_q^{(\ell)}\mathbf{b}_r+(-1)^\ell \mathbf{b}_r\mathrm{D}_{q+1}^{(\ell)}\\
				&+\mathbf{a}_{r-1}\mathrm{S}_q^{(\ell)}
				+(-1)^\ell \mathrm{S}_{q+1}^{(\ell)}\mathbf{a}_{r+1}.
			\end{aligned}
		\end{equation}
	\end{lemma}
	
	\begin{proof}
		Substitute respectively $(u,v)=(\ell,\ell)$, $(\ell+1,\ell)$ and
		$(\ell,\ell+1)$ into \eqref{eq:main-differential}.  For $\mathrm{D}$ the ordinary
		right-boundary sign is $(-1)^{2\ell}=1$.  The two shifted coefficients
		have sign $(-1)^\ell$.  For $\mathrm{U}$ and $\mathrm{V}$ the ordinary right-boundary
		sign is $-1$; their shifted signs are $(-1)^{\ell+1}$ and
		$(-1)^\ell$, respectively.  Identify the resulting critical generators
		using \eqref{eq:DUV-definitions} and \eqref{eq:RS-definitions}.
		For $\ell=0$ the same substitution gives
		\[
		\mathsf d_1\mathrm{U}_q^{(0)}=\mathbf{a}_q\mathrm{D}_q^{(0)}-\mathrm{D}_{q+1}^{(0)}\mathbf{a}_q,
		\qquad
		\mathsf d_1\mathrm{V}_q^{(0)}=\mathbf{b}_q\mathrm{D}_{q+1}^{(0)}-\mathrm{D}_q^{(0)}\mathbf{b}_q.
		\]
	\end{proof}
	
	\subsection{Evaluating on \texorpdfstring{$\mathrm{D}$}{D}: an odd cochain gives a loop-valued coboundary}
	\label{subsec:odd-to-even}
	
	Fix $1\leq\ell\leq m-1$ and let $\phi\in \mathsf C^{2\ell-1}$ be arbitrary.
	The input cochain has degree $2\ell-1$, while its argument
	$\mathrm{D}_q^{(\ell)}$ in the following calculation belongs to $\mathsf A_{2\ell}$.
	By definition and then by bimodule-linearity,
	\begin{equation}\label{eq:deltaD-before-values}
		\begin{aligned}
			(\mathsf d^{2\ell-1}\phi)(\mathrm{D}_q^{(\ell)})
			={}&\phi(\mathsf d\,\mathrm{D}_q^{(\ell)})\\
			={}&\mathbf{a}_{q+\ell-1}\phi(\mathrm{V}_q^{(\ell-1)})
			+\phi(\mathrm{U}_q^{(\ell-1)})\mathbf{b}_{q+\ell-1}\\
			&+(-1)^\ell \mathbf{b}_{q+\ell}\phi(\mathrm{U}_{q+1}^{(\ell-1)})\\
			&+(-1)^\ell\phi(\mathrm{V}_{q+1}^{(\ell-1)})\mathbf{a}_{q+\ell}.
		\end{aligned}
	\end{equation}
	Bimodule-linearity has so far been used only to take the exterior
	arrows outside $\phi$.  We now use it a second time, to determine
	$\phi$ on each remaining argument.
	
	Put $r=q+\ell-1$ and first suppose $q<m-\ell$, so that all four terms
	exist.  Since
	$\mathrm{V}_q^{(\ell-1)}=e_r\mathrm{V}_q^{(\ell-1)}e_{r+1}$,
	\begin{align}
		\phi(\mathrm{V}_q^{(\ell-1)})
		&=e_r\phi(\mathrm{V}_q^{(\ell-1)})e_{r+1}\notag\\
		&=\sum_{w\in\mathcal B}F_w^\phi(\mathrm{V}_q^{(\ell-1)})e_rwe_{r+1}\notag\\
		&=F_{\mathbf{b}_r}^\phi(\mathrm{V}_q^{(\ell-1)})\mathbf{b}_r.
		\label{eq:parallel-V-calculation}
	\end{align}
	Only the parallel arrow $\mathbf{b}_r$ survives.  Similarly, write
	\begin{equation}\label{eq:four-odd-values}
		\begin{aligned}
			\phi(\mathrm{U}_q^{(\ell-1)})&=\alpha_q \mathbf{a}_r,&
			\phi(\mathrm{V}_q^{(\ell-1)})&=\beta_q \mathbf{b}_r,\\
			\phi(\mathrm{U}_{q+1}^{(\ell-1)})&=\alpha_{q+1}\mathbf{a}_{r+1},&
			\phi(\mathrm{V}_{q+1}^{(\ell-1)})&=\beta_{q+1}\mathbf{b}_{r+1}.
		\end{aligned}
	\end{equation}
	For example, $\alpha_q=F_{\mathbf{a}_r}^\phi(\mathrm{U}_q^{(\ell-1)})$.
	The coefficients in \eqref{eq:four-odd-values} are those of the fixed
	input cochain of degree $2\ell-1$.
	
	Substitution into \eqref{eq:deltaD-before-values} gives
	\begin{align}
		(\mathsf d^{2\ell-1}\phi)(\mathrm{D}_q^{(\ell)})
		={}&\beta_q \mathbf{a}_r\mathbf{b}_r+\alpha_q \mathbf{a}_r\mathbf{b}_r\notag\\
		&+(-1)^\ell\alpha_{q+1}\mathbf{b}_{r+1}\mathbf{a}_{r+1}
		+(-1)^\ell\beta_{q+1}\mathbf{b}_{r+1}\mathbf{a}_{r+1}\notag\\
		={}&\bigl(\alpha_q+\beta_q
		+(-1)^\ell(\alpha_{q+1}+\beta_{q+1})\bigr)\boldsymbol{\omega}_{q+\ell}.
		\label{eq:deltaD-collected}
	\end{align}
	The final equality uses
	$\mathbf{a}_r\mathbf{b}_r=\mathbf{b}_{r+1}\mathbf{a}_{r+1}=\boldsymbol{\omega}_{r+1}$.
	All four terms contribute to the same loop; they must be collected
	before any coefficient is set equal to zero.  No vertex coefficient
	occurs in the result.
	
	At $q=m-\ell$ one has $r=m-1$.  The two shifted terms of
	\eqref{eq:dD} are absent, and
	\begin{equation}\label{eq:deltaD-boundary}
		(\mathsf d^{2\ell-1}\phi)(\mathrm{D}_{m-\ell}^{(\ell)})
		=(\alpha_{m-\ell}+\beta_{m-\ell})\mathbf{a}_{m-1}\mathbf{b}_{m-1}=0.
	\end{equation}
	This is zero for every input cochain.  It imposes no equation
	$\alpha_{m-\ell}+\beta_{m-\ell}=0$.
	
	Since $\mathsf d\phi$ is itself a bimodule map, it vanishes on all the
	other generators of $\mathsf A_{2\ell}$ by Lemma~\ref{lem:parallel}.
	Consequently \eqref{eq:deltaD-collected} and
	\eqref{eq:deltaD-boundary} describe its complete set of values.
	
	\subsection{Evaluating on \texorpdfstring{$\mathrm{U},\mathrm{V}$}{U,V}: an even cochain gives opposite arrow coefficients}
	\label{subsec:even-to-odd}
	
	Fix $0\leq\ell\leq m-2$ and let $\phi\in \mathsf C^{2\ell}$ be arbitrary.
	Set $r=q+\ell$.  By \eqref{eq:dU} and bimodule-linearity,
	\begin{equation}\label{eq:deltaU-full}
		\begin{aligned}
			(\mathsf d^{2\ell}\phi)(\mathrm{U}_q^{(\ell)})
			={}&\mathbf{a}_r\phi(\mathrm{D}_q^{(\ell)})
			+(-1)^{\ell+1}\phi(\mathrm{D}_{q+1}^{(\ell)})\mathbf{a}_r\\
			&-\phi(\mathrm{R}_q^{(\ell)})\mathbf{b}_{r-1}
			+(-1)^{\ell+1}\mathbf{b}_{r+1}\phi(\mathrm{R}_{q+1}^{(\ell)}).
		\end{aligned}
	\end{equation}
	The last two values are zero: whenever the arguments exist,
	\[
	\phi(\mathrm{R}_q^{(\ell)})
	\in e_{r+1}\Lambda e_{r-1}=0,
	\qquad
	\phi(\mathrm{R}_{q+1}^{(\ell)})
	\in e_{r+2}\Lambda e_r=0.
	\]
	This is why these two terms disappear \emph{after} applying a cochain.
	They were retained in the chain boundary \eqref{eq:dU}.
	
	On a closed argument, bimodule-linearity instead gives
	\begin{equation}\label{eq:parallel-D-calculation}
		\begin{aligned}
			\phi(\mathrm{D}_q^{(\ell)})
			&=e_r\phi(\mathrm{D}_q^{(\ell)})e_r\\
			&=\sum_{w\in\mathcal B}F_w^\phi(\mathrm{D}_q^{(\ell)})e_rwe_r
			=E_qe_r+\Omega_q\boldsymbol{\omega}_r.
		\end{aligned}
	\end{equation}
	Similarly,
	$\phi(\mathrm{D}_{q+1}^{(\ell)})=E_{q+1}e_{r+1}+\Omega_{q+1}\boldsymbol{\omega}_{r+1}$
	when $r+1<m$; if $r+1=m$ only $E_{q+1}e_m$ remains.
	Thus, away from that endpoint, \eqref{eq:deltaU-full} reduces to
	\begin{align}
		(\mathsf d^{2\ell}\phi)(\mathrm{U}_q^{(\ell)})
		={}&\mathbf{a}_r(E_qe_r+\Omega_q\boldsymbol{\omega}_r)\notag\\
		&+(-1)^{\ell+1}(E_{q+1}e_{r+1}+\Omega_{q+1}\boldsymbol{\omega}_{r+1})\mathbf{a}_r\notag\\
		={}&\bigl(E_q+(-1)^{\ell+1}E_{q+1}\bigr)\mathbf{a}_r.
		\label{eq:deltaU-collected}
	\end{align}
	We used $\mathbf{a}_r e_r=e_{r+1}\mathbf{a}_r=\mathbf{a}_r$ and
	$\mathbf{a}_r\boldsymbol{\omega}_r=\boldsymbol{\omega}_{r+1}\mathbf{a}_r=0$.  The latter products have arrow length three.
	
	Likewise, \eqref{eq:dV} gives
	\begin{equation}\label{eq:deltaV-full}
		\begin{aligned}
			(\mathsf d^{2\ell}\phi)(\mathrm{V}_q^{(\ell)})
			={}&-\phi(\mathrm{D}_q^{(\ell)})\mathbf{b}_r
			+(-1)^\ell \mathbf{b}_r\phi(\mathrm{D}_{q+1}^{(\ell)})\\
			&+\mathbf{a}_{r-1}\phi(\mathrm{S}_q^{(\ell)})
			+(-1)^\ell\phi(\mathrm{S}_{q+1}^{(\ell)})\mathbf{a}_{r+1}.
		\end{aligned}
	\end{equation}
	Both $S$-values vanish because their endpoints differ by two.  Therefore
	\begin{align}
		(\mathsf d^{2\ell}\phi)(\mathrm{V}_q^{(\ell)})
		={}&-(E_qe_r+\Omega_q\boldsymbol{\omega}_r)\mathbf{b}_r\notag\\
		&+(-1)^\ell \mathbf{b}_r(E_{q+1}e_{r+1}+\Omega_{q+1}\boldsymbol{\omega}_{r+1})\notag\\
		={}&-\bigl(E_q+(-1)^{\ell+1}E_{q+1}\bigr)\mathbf{b}_r.
		\label{eq:deltaV-collected}
	\end{align}
	Here $e_r\mathbf{b}_r=\mathbf{b}_r e_{r+1}=\mathbf{b}_r$ and $\boldsymbol{\omega}_r\mathbf{b}_r=\mathbf{b}_r\boldsymbol{\omega}_{r+1}=0$.
	
	At $r=m-1$, the missing loop $\boldsymbol{\omega}_m$ is simply omitted from the input
	value.  The same final formulas hold because $\mathbf{a}_{m-1}$ and $\mathbf{b}_{m-1}$
	are nonzero, whereas their products with $\boldsymbol{\omega}_{m-1}$ are zero.
	Thus the last neighboring relation in
	\eqref{eq:deltaU-collected} does \emph{not} disappear, in contrast to
	\eqref{eq:deltaD-boundary}.
	For $\ell=0$ the formulas follow directly from the two boundaries at the
	end of Lemma~\ref{lem:DUV-boundaries}; they read
	\[
	(\mathsf d^0\phi)(\mathsf{[\mathbf{a}_q]})=(E_q-E_{q+1})\mathbf{a}_q,
	\qquad
	(\mathsf d^0\phi)(\mathsf{[\mathbf{b}_q]})=-(E_q-E_{q+1})\mathbf{b}_q.
	\]
	Finally, all other generators of $\mathsf A_{2\ell+1}$ have endpoints differing
	by at least three, so $\mathsf d\phi$ vanishes on them by parallelism.
	We have therefore computed the whole coboundary, not only a subsystem
	of its equations.
	
	\subsection{A signed difference operator}
	\label{subsec:nabla}
	
	The preceding computations can now be recorded as linear maps between
	coordinate spaces.  The difference notation is introduced only at this
	point: it abbreviates the coefficients already derived from the chain
	boundaries and multiplication in $\Lambda$.
	
	For $s\geq1$ and $\ell\geq0$, define
	\begin{equation}\label{eq:nabla-definition}
		\nabla_{\ell,s}:\Bbbk^s\longrightarrow \Bbbk^{s-1},\qquad
		(\nabla_{\ell,s}x)_q=x_q+(-1)^{\ell+1}x_{q+1}
		\quad(1\leq q<s).
	\end{equation}
	When $s=1$ the codomain is zero and this is the zero map.
	For $s\geq2$ its matrix has $1$ on the main diagonal,
	$(-1)^{\ell+1}$ on the next diagonal, and zeros elsewhere:
	\begin{equation}\label{eq:nabla-matrix}
		[\nabla_{\ell,s}]=
		\begin{pmatrix}
			1&(-1)^{\ell+1}&0&\cdots&0\\
			0&1&(-1)^{\ell+1}&\ddots&\vdots\\
			\vdots&\ddots&\ddots&\ddots&0\\
			0&\cdots&0&1&(-1)^{\ell+1}
		\end{pmatrix}.
	\end{equation}
	Thus $\nabla_{0,3}(x_1,x_2,x_3)=(x_1-x_2,x_2-x_3)$, whereas
	$\nabla_{1,3}(x_1,x_2,x_3)=(x_1+x_2,x_2+x_3)$.
	The second index records the length of the input vector; it is not a
	cohomological degree.
	
	\begin{theorem}[The complete cochain differential]\label{thm:delta-nabla}
		Use the coordinates in Proposition~\ref{prop:cochain-coordinates} and
		put $s=m-\ell$.
		For $0\leq\ell\leq m-1$,
		\begin{equation}\label{eq:delta-even-nabla}
			\boxed{\mathsf d^{2\ell}(\mathbf E,\boldsymbol\Omega)
				=\bigl(\nabla_{\ell,s}\mathbf E,
				-\nabla_{\ell,s}\mathbf E\bigr).}
		\end{equation}
		The two output vectors are respectively the coefficients at the
		$U$- and $V$-generators.  At $\ell=m-1$ both are empty, since
		$\mathsf C^{2m-1}=0$.
		For $0\leq\ell\leq m-2$,
		\begin{equation}\label{eq:delta-odd-nabla}
			\boxed{\mathsf d^{2\ell+1}(\boldsymbol\alpha,\boldsymbol\beta)
				=\bigl(0,\nabla_{\ell,s-1}
				(\boldsymbol\alpha+\boldsymbol\beta)\bigr).}
		\end{equation}
		Here the zero vector belongs to $\Bbbk^{s-1}$ and consists of the vertex
		coefficients, while the second vector belongs to $\Bbbk^{s-2}$ and consists
		of the loop coefficients.  In particular, the maps have types
		\begin{equation}\label{eq:delta-types}
			\begin{aligned}
				\mathsf d^{2\ell}:&\ \Bbbk^s\oplus \Bbbk^{s-1}
				\longrightarrow \Bbbk^{s-1}\oplus \Bbbk^{s-1},\\
				\mathsf d^{2\ell+1}:&\ \Bbbk^{s-1}\oplus \Bbbk^{s-1}
				\longrightarrow \Bbbk^{s-1}\oplus \Bbbk^{s-2}.
			\end{aligned}
		\end{equation}
		The second line is used only for $s\geq2$.
	\end{theorem}
	
	\begin{proof}
		The two components of \eqref{eq:delta-even-nabla} are exactly
		\eqref{eq:deltaU-collected} and \eqref{eq:deltaV-collected}, with the
		exterior arrows factored off.  For \eqref{eq:delta-odd-nabla}, replace
		$\ell$ by $\ell+1$ in \eqref{eq:deltaD-collected}.  The resulting loop
		coefficient is
		\[
		\alpha_q+\beta_q
		+(-1)^{\ell+1}(\alpha_{q+1}+\beta_{q+1}),
		\qquad 1\leq q\leq s-2.
		\]
		There is no vertex coefficient, and the value on the last closed chain
		is zero by \eqref{eq:deltaD-boundary}.  Parallelism accounts for all
		other generators.  This proves both formulas, including their boundary
		ranges.
	\end{proof}
	
	\begin{remark}\label{rem:delta-square-coordinates}
		These formulas display $\mathsf d^{n+1}\mathsf d^n=0$ although successive difference
		maps need not have zero composition: for $s\geq3$, the map
		$\nabla_{\ell,s-1}\nabla_{\ell,s}$ is generally nonzero.  Indeed,
		\[
		(\mathbf E,\boldsymbol\Omega)\longmapsto(\mathbf t,-\mathbf t)
		\longmapsto(0,0),
		\qquad
		(\boldsymbol\alpha,\boldsymbol\beta)\longmapsto(0,\boldsymbol\Xi)
		\longmapsto(0,0).
		\]
		The first composition vanishes because the two arrow vectors sum to
		zero.  The second vanishes because the next map uses only vertex
		coordinates.  A nonzero pair $(\mathbf t,-\mathbf t)$ is not itself
		the zero cochain: its components are evaluated on different generators.
	\end{remark}
	
	\begin{lemma}\label{lem:nabla-kernel-image}
		The map $\nabla_{\ell,s}$ is surjective, and its kernel is the line
		spanned by the vector
		\begin{equation}\label{eq:nabla-kernel-vector}
			\mathbf v_{\ell,s}
			=\bigl((-1)^{\ell(q-1)}\bigr)_{q=1}^{s}.
		\end{equation}
		These assertions hold over every field, including characteristic two.
	\end{lemma}
	
	\begin{proof}
		The homogeneous equations give
		$x_{q+1}=(-1)^\ell x_q$, and hence
		$x_q=(-1)^{\ell(q-1)}x_1$.
		For any prescribed $y=(y_1,\ldots,y_{s-1})$, choose $x_1=0$ and define
		\begin{equation}\label{eq:nabla-right-inverse}
			x_{q+1}=(-1)^{\ell+1}(y_q-x_q),\qquad 1\leq q<s.
		\end{equation}
		Then $\nabla_{\ell,s}x=y$.  The case $s=1$ has zero codomain and the
		same kernel description.  No division is used.
	\end{proof}
	
	\begin{remark}[Index shifts]\label{rem:index-shifts}
		For any critical generator $\gamma$, let $\gamma[t]$ denote the
		generator obtained by shifting every arrow index by $t$, provided
		the resulting chain remains admissible. Then
		$\mathrm{D}_{q+1}^{(\ell)}=\mathrm{D}_q^{(\ell)}[1]$, and likewise for $\mathrm{U},\mathrm{V}$.
		The relation $\nabla_{\ell,s}\mathbf E=0$ therefore propagates a vertex
		coefficient along this finite sequence of shifts, with the output
		idempotent shifted at the same time.  For odd cochains it propagates
		$\alpha_q+\beta_q$, not $\alpha_q$ and $\beta_q$ separately.
		Writing $x_q=(-1)^{\ell(q-1)}y_q$ gives
		\[
		(\nabla_{\ell,s}x)_q
		=(-1)^{\ell(q-1)}(y_q-y_{q+1}).
		\]
		After this sign change the equations simply require constancy along
		adjacent shifts.  The propagation stops at $q=1$ on the left; it is not
		an action of a translation automorphism on the entire bound quiver.
	\end{remark}
	
	\section{Cocycles, coboundaries, and explicit cohomology classes}
	\label{sec:cohomology}
	
	We now solve the complete system.  Throughout this section, coordinate
	vectors are evaluated on the generators and in the order specified in
	Proposition~\ref{prop:cochain-coordinates}.  In particular, a displayed
	formula for a cochain is completed by setting it to zero on all other
	critical generators and extending bimodule-linearly.
	
	\subsection{All cocycles}
	
	\begin{proposition}\label{prop:all-cocycles}
		For $0\leq\ell\leq m-1$, put $s=m-\ell$.
		Every even cocycle is uniquely of the form
		\begin{equation}\label{eq:all-even-cocycles}
			\phi(\mathrm{D}_q^{(\ell)})=
			\begin{cases}
				(-1)^{\ell(q-1)}\epsilon_\ell e_{q+\ell}
				+\Omega_q \boldsymbol{\omega}_{q+\ell},&1\leq q<s,\\
				(-1)^{\ell(s-1)}\epsilon_\ell e_m,&q=s,
			\end{cases}
		\end{equation}
		where $\epsilon_\ell,\Omega_1,\ldots,\Omega_{s-1}\in \Bbbk$ are arbitrary.
		Hence $\dim_{\Bbbk}Z^{2\ell}=s$.
		
		For $0\leq\ell\leq m-2$, put $h=m-\ell-1$.
		Every odd cocycle is uniquely of the form
		\begin{equation}\label{eq:all-odd-cocycles}
			\begin{aligned}
				\psi(\mathrm{U}_q^{(\ell)})&=\lambda_q \mathbf{a}_{q+\ell},\\
				\psi(\mathrm{V}_q^{(\ell)})&=
				\bigl((-1)^{\ell(q-1)}\sigma_\ell-\lambda_q\bigr)\mathbf{b}_{q+\ell},
				\qquad 1\leq q\leq h,
			\end{aligned}
		\end{equation}
		where $\sigma_\ell,\lambda_1,\ldots,\lambda_h\in \Bbbk$ are arbitrary.
		Hence $\dim_{\Bbbk}Z^{2\ell+1}=h+1=m-\ell$.
	\end{proposition}
	
	\begin{proof}
		By \eqref{eq:delta-even-nabla}, an even cochain is a cocycle exactly
		when $\mathbf E\in\ker\nabla_{\ell,s}$; its loop vector is unrestricted.
		Lemma~\ref{lem:nabla-kernel-image} gives \eqref{eq:all-even-cocycles}.
		The initial scalar is
		\[
		\epsilon_\ell=F_{e_{\ell+1}}^\phi(\mathrm{D}_1^{(\ell)}).
		\]
		For an odd cochain the condition is
		$\boldsymbol\alpha+\boldsymbol\beta\in\ker\nabla_{\ell,h}$.
		Consequently
		\[
		\alpha_q+\beta_q=(-1)^{\ell(q-1)}\sigma_\ell,
		\quad
		\sigma_\ell=\alpha_1+\beta_1.
		\]
		Choose the $\alpha_q=\lambda_q$ freely, and solve for the $\beta_q$.
		Theorem~\ref{thm:delta-nabla} was obtained on all generators, so there
		are no additional cocycle conditions.  Counting the independent
		parameters proves the dimensions.
	\end{proof}
	
	For example, even cocycles in degree two have vertex coefficients
	$\epsilon_1,-\epsilon_1,\epsilon_1,\ldots$, whereas in degree four the
	vertex coefficients are all equal.  In degree one the sums of the two
	opposite arrow coefficients are constant; in degree three these sums
	alternate in sign.  In characteristic two these sign distinctions
	disappear, but the number of parameters is unchanged.
	
	\subsection{All coboundaries, with explicit primitives}
	
	\begin{proposition}\label{prop:all-coboundaries}
		For $1\leq\ell\leq m-1$ and $s=m-\ell$,
		\begin{equation}\label{eq:all-even-boundaries}
			B^{2\ell}
			=\{(\mathbf E,\boldsymbol\Omega):
			\mathbf E=0,\ \boldsymbol\Omega\in \Bbbk^{s-1}\}.
		\end{equation}
		In particular, $\dim_{\Bbbk}B^{2\ell}=s-1$.
		For $0\leq\ell\leq m-2$ and $h=m-\ell-1$,
		\begin{equation}\label{eq:all-odd-boundaries}
			B^{2\ell+1}
			=\{(\boldsymbol\alpha,\boldsymbol\beta):
			\boldsymbol\alpha=\mathbf t,\
			\boldsymbol\beta=-\mathbf t,\ \mathbf t\in \Bbbk^h\}.
		\end{equation}
		Thus $\dim_{\Bbbk}B^{2\ell+1}=h$.
	\end{proposition}
	
	\begin{proof}
		An even coboundary is the image of a cochain of degree $2\ell-1$.
		Using \eqref{eq:delta-odd-nabla} with $\ell-1$ in place of $\ell$, it
		has zero vertex coefficients and loop coefficients
		\[
		\boldsymbol\Omega
		=\nabla_{\ell-1,s}(\boldsymbol\alpha+\boldsymbol\beta).
		\]
		Surjectivity in Lemma~\ref{lem:nabla-kernel-image} proves that every loop
		vector occurs.  More explicitly, for a prescribed $\Omega$ choose
		$t_1=0$ and recursively set
		\begin{equation}\label{eq:even-primitive-recursion}
			t_{q+1}=(-1)^\ell(\Omega_q-t_q),\qquad1\leq q<s.
		\end{equation}
		Define $\chi\in \mathsf C^{2\ell-1}$ by
		\begin{equation}\label{eq:even-primitive}
			\chi(\mathrm{U}_q^{(\ell-1)})=t_q \mathbf{a}_{q+\ell-1},\qquad
			\chi(\mathrm{V}_q^{(\ell-1)})=0\quad(1\leq q\leq s).
		\end{equation}
		Then $\mathsf d\chi=(0,\boldsymbol\Omega)$.
		At $s=1$ there are no loop coordinates and this image is zero.
		
		An odd coboundary has the form
		$(\nabla_{\ell,h+1}\mathbf E,-\nabla_{\ell,h+1}\mathbf E)$ by
		\eqref{eq:delta-even-nabla}.  Every $\mathbf t\in \Bbbk^h$ occurs.
		An explicit primitive is obtained by setting $x_1=0$ and
		\begin{equation}\label{eq:odd-primitive-recursion}
			x_{q+1}=(-1)^{\ell+1}(t_q-x_q),\qquad1\leq q\leq h,
		\end{equation}
		then defining $\chi\in \mathsf C^{2\ell}$ by
		\begin{equation}\label{eq:odd-primitive}
			\chi(\mathrm{D}_q^{(\ell)})=x_q e_{q+\ell}
			\quad(1\leq q\leq h+1),
		\end{equation}
		with all loop coefficients zero.  Its coboundary has arrow coordinates
		$(\mathbf t,-\mathbf t)$.  These constructions prove both the claimed
		images and the stated dimensions over every field.
	\end{proof}
	
	The primitives above are arbitrary cochains, not required to be
	cocycles.  Requiring them to satisfy their own cocycle equations would
	incorrectly shrink the space of coboundaries.
	
	\subsection{Representatives \texorpdfstring{$\eta_\ell$ and $\theta_\ell$}{eta and theta}}
	
	For $0\leq\ell\leq m-1$, define an even cochain
	$\eta_\ell\in \mathsf C^{2\ell}$ by
	\begin{equation}\label{eq:eta-representatives}
		\boxed{\eta_\ell(\mathrm{D}_q^{(\ell)})
			=(-1)^{\ell(q-1)}e_{q+\ell}
			\quad(1\leq q\leq m-\ell).}
	\end{equation}
	For $0\leq\ell\leq m-2$, define
	$\theta_\ell\in \mathsf C^{2\ell+1}$ by
	\begin{equation}\label{eq:theta-representatives}
		\boxed{
			\theta_\ell(\mathrm{U}_q^{(\ell)})=0,\qquad
			\theta_\ell(\mathrm{V}_q^{(\ell)})
			=(-1)^{\ell(q-1)}\mathbf{b}_{q+\ell}
			\quad(1\leq q\leq m-\ell-1).}
	\end{equation}
	Values on other critical generators are zero.  These are well-defined
	bimodule maps by parallelism, and
	Proposition~\ref{prop:all-cocycles} shows that they are cocycles.
	The indices on these symbols label the chosen representatives; no
	identity expressing them as cup powers or cup products is being
	asserted.
	
	\begin{theorem}[Hochschild cohomology and representatives]\label{thm:HH}
		For the algebra \eqref{eq:ideal} over any field $\Bbbk$,
		\begin{equation}\label{eq:HH-result}
			HH^n(\Lambda)=
			\begin{cases}
				\Bbbk\cdot1\ \oplus\displaystyle\bigoplus_{i=1}^{m-1}\Bbbk\boldsymbol{\omega}_i,&n=0,\\[1mm]
				\Bbbk\cdot[\eta_\ell],&n=2\ell,\quad1\leq\ell\leq m-1,\\
				\Bbbk\cdot[\theta_\ell],&n=2\ell+1,\quad0\leq\ell\leq m-2,\\
				0,&n\geq2m-1.
			\end{cases}
		\end{equation}
		The first line identifies $HH^0(\Lambda)$ with the center of $\Lambda$;
		the other lines specify bases of the corresponding vector spaces.
		In particular,
		\begin{equation}\label{eq:HH-dimensions}
			\boxed{\dim_{\Bbbk}HH^n(\Lambda)=
				\begin{cases}
					m,&n=0,\\
					1,&1\leq n\leq2m-2,\\
					0,&n\geq2m-1.
			\end{cases}}
		\end{equation}
		Every positive even cocycle is cohomologous to a unique multiple of
		$\eta_\ell$, and every odd cocycle to a unique multiple of $\theta_\ell$.
	\end{theorem}
	
	\begin{proof}
		In degree zero, \eqref{eq:all-even-cocycles} gives
		$E_1=\cdots=E_m$, with arbitrary coefficients of $\boldsymbol{\omega}_1,\ldots,\boldsymbol{\omega}_{m-1}$.
		There are no degree-zero coboundaries.  The cochain $\eta_0$ corresponds
		to $1=e_1+\cdots+e_m$.  For each $i<m$, the cochain with value $\boldsymbol{\omega}_i$ on
		$\mathrm{D}_i^{(0)}$ and zero on the other vertex generators corresponds to the
		central element $\boldsymbol{\omega}_i$.  These loops commute with all vertices and
		annihilate all positive paths on either side.  Conversely, commuting
		with the vertices forces a central element to be diagonal in endpoint
		type, and commuting with every $\mathbf{a}_i,\mathbf{b}_i$ makes all its vertex
		coefficients equal.  This proves the first line of \eqref{eq:HH-result}.
		
		For $\ell\geq1$, subtract $\epsilon_\ell\eta_\ell$ from an even
		cocycle \eqref{eq:all-even-cocycles}.  The difference has zero vertex
		coefficients and hence is a coboundary by
		\eqref{eq:all-even-boundaries}, with the explicit primitive
		\eqref{eq:even-primitive}.  The cocycle $\eta_\ell$ is not a
		coboundary: its value on $\mathrm{D}_1^{(\ell)}$ has vertex coefficient one,
		whereas every even coboundary has zero vertex coefficients.  This
		also proves uniqueness of $\epsilon_\ell$.
		
		For an odd cocycle \eqref{eq:all-odd-cocycles}, subtract
		$\sigma_\ell\theta_\ell$.  The resulting arrow coordinates are
		$(\boldsymbol\lambda,-\boldsymbol\lambda)$, which form a coboundary
		by \eqref{eq:all-odd-boundaries}, with primitive
		\eqref{eq:odd-primitive}.  On the other hand, any odd coboundary has
		$\alpha_1+\beta_1=0$, while this sum is one for $\theta_\ell$.
		Thus $[\theta_\ell]\ne0$ and $\sigma_\ell$ is unique.
		
		Finally, $\mathsf A_n=0$ for $n>2m-2$, so the cochain spaces and their
		cohomology vanish there.  All arguments use only addition and signs;
		no restriction on the characteristic is needed.
	\end{proof}
	
	The quotient coordinates can be recorded directly in terms of the
	original functionals.  In positive even degree the isomorphism
	$HH^{2\ell}(\Lambda)\to \Bbbk$ is
	\begin{equation}\label{eq:even-quotient-coordinate}
		[\phi]\longmapsto
		F_{e_{\ell+1}}^\phi(\mathrm{D}_1^{(\ell)})=\epsilon_\ell,
	\end{equation}
	and in odd degree the isomorphism $HH^{2\ell+1}(\Lambda)\to \Bbbk$ is
	\begin{equation}\label{eq:odd-quotient-coordinate}
		[\psi]\longmapsto
		F_{\mathbf{a}_{\ell+1}}^\psi(\mathrm{U}_1^{(\ell)})
		+F_{\mathbf{b}_{\ell+1}}^\psi(\mathrm{V}_1^{(\ell)})=\sigma_\ell.
	\end{equation}
	The proof shows that these maps are well-defined on cohomology and that
	their kernels on cocycles are exactly the coboundaries.
	
	For clarity, all dimensions are collected below.  The spaces $Z^n$
	and $B^n$ are those of \eqref{eq:cochain-complex}.
	\begin{equation}\label{eq:dimension-table}
		\begin{array}{c|c|c|c|c}
			n&\dim \mathsf C^n&\dim Z^n&\dim B^n&\dim HH^n\\\hline
			0&2m-1&m&0&m\\
			2\ell\ (1\leq\ell\leq m-1)
			&2m-2\ell-1&m-\ell&m-\ell-1&1\\
			2\ell+1\ (0\leq\ell\leq m-2)
			&2m-2\ell-2&m-\ell&m-\ell-1&1\\
			n\geq2m-1&0&0&0&0.
		\end{array}
	\end{equation}
	In particular, the Hilbert series of the graded cohomology vector space is
	\begin{equation}\label{eq:HH-series}
		\sum_{n\geq0}\dim_{\Bbbk}HH^n(\Lambda)t^n
		=m+t+t^2+\cdots+t^{2m-2}.
	\end{equation}
	This is an additive computation with explicit representatives.  A
	multiplicative description and a Gerstenhaber bracket require additional
	comparison or diagonal data; these are supplied in the next two sections.
	
	\section{Some examples}
	
	\begin{example}[Three vertices: the entire cochain complex]
		\label{ex:HH-m3}
		Let $m=3$.  Use the generators of Example~\ref{ex:m3}, now with
		$\mathrm{D}_1^{(1)}=\mathsf{[\mathbf{a}_1|\mathbf{b}_1]}$, $\mathrm{D}_2^{(1)}=\mathsf{[\mathbf{a}_2|\mathbf{b}_2]}$,
		$\mathrm{U}_1^{(1)}=\mathrm{U}$, $\mathrm{V}_1^{(1)}=\mathrm{V}$, and $\mathrm{D}_1^{(2)}=\mathrm{W}$.
		The successive coordinate orders are
		\[
		\begin{array}{c|l}
			\mathsf C^0&(E_1,E_2,E_3,\Omega_1,\Omega_2)\\
			\mathsf C^1&(\alpha_1,\alpha_2,\beta_1,\beta_2)\\
			\mathsf C^2&(E_1,E_2,\Omega_1)\\
			\mathsf C^3&(\alpha_1,\beta_1)\\
			\mathsf C^4&(E_1).
		\end{array}
		\]
		Repeated coordinate names in different rows denote different cochains
		and different arguments.  For instance, the coordinates in $\mathsf C^2$ mean
		\[
		\phi(\mathsf{[\mathbf{a}_1|\mathbf{b}_1]})=E_1e_2+\Omega_1\boldsymbol{\omega}_2,
		\qquad \phi(\mathsf{[\mathbf{a}_2|\mathbf{b}_2]})=E_2e_3,
		\]
		while $\phi(\mathsf{[\mathbf{a}_2|\mathbf{a}_1]})=\phi(\mathsf{[\mathbf{b}_1|\mathbf{b}_2]})=0$ by parallelism.
		
		With these orders the four coboundary matrices are
		\begin{equation}\label{eq:m3-cochain-matrices}
			\begin{aligned}
				[\mathsf d^0]&=
				\begin{pmatrix}
					1&-1&0&0&0\\0&1&-1&0&0\\
					-1&1&0&0&0\\0&-1&1&0&0
				\end{pmatrix},&
				[\mathsf d^1]&=
				\begin{pmatrix}
					0&0&0&0\\0&0&0&0\\1&-1&1&-1
				\end{pmatrix},\\[2mm]
				[\mathsf d^2]&=
				\begin{pmatrix}1&1&0\\-1&-1&0\end{pmatrix},&
				[\mathsf d^3]&=\begin{pmatrix}0&0\end{pmatrix}.
			\end{aligned}
		\end{equation}
		For example, evaluating a degree-one cochain on $\mathsf d\mathsf{[\mathbf{a}_1|\mathbf{b}_1]}$ gives
		\[
		\begin{aligned}
			(\mathsf d^1\phi)(\mathsf{[\mathbf{a}_1|\mathbf{b}_1]})
			&=\mathbf{a}_1(\beta_1\mathbf{b}_1)+(\alpha_1\mathbf{a}_1)\mathbf{b}_1
			-\mathbf{b}_2(\alpha_2\mathbf{a}_2)-(\beta_2\mathbf{b}_2)\mathbf{a}_2\\
			&=(\alpha_1+\beta_1-\alpha_2-\beta_2)\boldsymbol{\omega}_2.
		\end{aligned}
		\]
		The coefficient on $\mathsf{[\mathbf{a}_2|\mathbf{b}_2]}$ is zero because $\mathbf{a}_2\mathbf{b}_2=0$.
		This explains both the nonzero row and the zero vertex rows of
		$[\mathsf d^1]$.
		
		In the next degree, $\mathsf d\,\mathrm{U}=\mathbf{a}_2\mathrm{D}_1^{(1)}-\mathrm{R}_1\mathbf{b}_1+\mathrm{D}_2^{(1)}\mathbf{a}_2$.
		The value of a degree-two cochain on $\mathrm{R}_1=\mathsf{[\mathbf{a}_2|\mathbf{a}_1]}$ is zero.  Hence
		\[
		(\mathsf d^2\phi)(\mathrm{U})=(E_1+E_2)\mathbf{a}_2,
		\qquad
		(\mathsf d^2\phi)(\mathrm{V})=-(E_1+E_2)\mathbf{b}_2,
		\]
		and $\Omega_1$ disappears because the exterior products have arrow
		length three.  Finally $\mathsf d\,\mathrm{W}=\mathbf{a}_2\mathrm{V}+\mathrm{U}\mathbf{b}_2$ gives
		$(\mathsf d^3\psi)(\mathrm{W})=(\alpha_1+\beta_1)\mathbf{a}_2\mathbf{b}_2=0$.
		
		The full spaces of cocycles and coboundaries, in the displayed orders,
		are
		\begin{equation}\label{eq:m3-ZB}
			\begin{array}{c|c|c}
				n&Z^n&B^n\\\hline
				0&\{(c,c,c,u,v)\}&0\\
				1&\{\alpha_1+\beta_1=\alpha_2+\beta_2\}
				&\{(t_1,t_2,-t_1,-t_2)\}\\
				2&\{(c,-c,u)\}&\{(0,0,u)\}\\
				3&\Bbbk^2&\{(t,-t)\}\\
				4&\Bbbk&0.
			\end{array}
		\end{equation}
		The scalars in the parametrized rows are arbitrary; in degree one
		the displayed linear relation is the only constraint.
		Thus the cohomology dimensions are $(3,1,1,1,1)$.
		The representatives from Theorem~\ref{thm:HH} specialize to
		\begin{align*}
			\theta_0(\mathsf{[\mathbf{a}_1]})&=\theta_0(\mathsf{[\mathbf{a}_2]})=0,&
			\theta_0(\mathsf{[\mathbf{b}_1]})&=\mathbf{b}_1,&\theta_0(\mathsf{[\mathbf{b}_2]})&=\mathbf{b}_2,\\
			\eta_1(\mathrm{D}_1^{(1)})&=e_2,&
			\eta_1(\mathrm{D}_2^{(1)})&=-e_3,\\
			\theta_1(\mathrm{U})&=0,&\theta_1(\mathrm{V})&=\mathbf{b}_2,\\
			\eta_2(\mathrm{W})&=e_3.
		\end{align*}
		All unspecified values are zero.  In particular, the two values of
		$\eta_1$ must be taken together: a map supported only on $\mathrm{D}_1^{(1)}$
		with value $e_2$ would have nonzero coboundary.
	\end{example}
	
	\begin{example}[Four vertices: neighboring shifts and the last boundary]
		\label{ex:HH-m4}
		Let $m=4$ and begin with a degree-two cochain:
		\[
		\begin{aligned}
			\phi(\mathrm{D}_1^{(1)})&=E_1e_2+\Omega_1\boldsymbol{\omega}_2,\\
			\phi(\mathrm{D}_2^{(1)})&=E_2e_3+\Omega_2\boldsymbol{\omega}_3,\\
			\phi(\mathrm{D}_3^{(1)})&=E_3e_4.
		\end{aligned}
		\]
		For example, the full chain boundary is
		\[
		\mathsf d_3\mathrm{U}_1^{(1)}=\mathbf{a}_2\mathrm{D}_1^{(1)}+\mathrm{D}_2^{(1)}\mathbf{a}_2
		-\mathrm{R}_1^{(1)}\mathbf{b}_1+\mathbf{b}_3\mathrm{R}_2^{(1)}.
		\]
		After applying $\phi$, the $R$-values vanish by parallelism, and both
		loop products vanish.  Thus $(\mathsf d^2\phi)(\mathrm{U}_1^{(1)})=(E_1+E_2)\mathbf{a}_2$.
		Repeating this calculation gives
		\[
		\mathsf d^2(\mathbf E,\boldsymbol\Omega)
		=\bigl((E_1+E_2,E_2+E_3),(-E_1-E_2,-E_2-E_3)\bigr).
		\]
		The cocycle condition is $\mathbf E=(t,-t,t)$; both loop coefficients
		are arbitrary and form coboundaries.
		
		For a degree-three cochain with coordinates
		$\boldsymbol\alpha=(\alpha_1,\alpha_2)$ and
		$\boldsymbol\beta=(\beta_1,\beta_2)$, one obtains
		\[
		\begin{aligned}
			(\mathsf d^3\psi)(\mathrm{D}_1^{(2)})
			&=(\alpha_1+\beta_1+\alpha_2+\beta_2)\boldsymbol{\omega}_3,\\
			(\mathsf d^3\psi)(\mathrm{D}_2^{(2)})&=0.
		\end{aligned}
		\]
		The second equality follows from the boundary product $\mathbf{a}_3\mathbf{b}_3=0$,
		not from an additional equation on $\alpha_2+\beta_2$.
		Every degree-three cocycle is consequently of the form
		\[
		\boldsymbol\alpha=(\lambda_1,\lambda_2),\qquad
		\boldsymbol\beta=(\sigma-\lambda_1,-\sigma-\lambda_2).
		\]
		Subtracting a coboundary removes both $\lambda_q$ and leaves
		$\sigma\theta_1$, where
		$\theta_1(\mathrm{V}_1^{(1)})=\mathbf{b}_2$ and $\theta_1(\mathrm{V}_2^{(1)})=-\mathbf{b}_3$.
	\end{example}
	
	\begin{remark}[The smallest quiver and characteristic two]
		For $m=2$ the complex has dimensions $3,2,1$.
		Its only nonzero possible coboundary is
		\[
		\mathsf d^0(E_1,E_2,\Omega_1)=(E_1-E_2,-E_1+E_2),
		\qquad\mathsf d^1=0.
		\]
		The latter equality is the boundary relation $\mathbf{a}_1\mathbf{b}_1=0$.
		Thus $(\dim HH^0,\dim HH^1,\dim HH^2)=(2,1,1)$, in agreement with
		Theorem~\ref{thm:HH}.  In characteristic two the sign $-1$ becomes $1$,
		but the kernel and surjectivity proof for $\nabla$ and the primitives
		\eqref{eq:even-primitive-recursion}--\eqref{eq:odd-primitive} remain
		valid.  No extra classes arise from this change of signs.
	\end{remark}

\section{A shuffle coproduct and the cup product}\label{cup:sec}
The shuffle map was constructed in Proposition~\ref{cup:prop:comparison}
and used in the contraction of Theorem~\ref{thm:resolution}.
Its image is stable under bar deconcatenation, giving a strictly
coassociative counital coproduct. The resulting cup algebra agrees
with \cite[Theorem~4.2]{LiXuChen}; the formulas below retain the signs
of our chosen representatives over an arbitrary field $\Bbbk$.

\subsection{A counital coassociative coproduct}

Let $\Delta_{\mathrm{bar}}$ denote the usual deconcatenation diagonal
of the relative two-sided bar resolution. On a composable normal
bar cell it cuts the list of internal blocks in every position and
inserts the idempotent of the intermediate vertex:
\begin{equation}\label{cup:eq:bar-deconcatenation}
 \Delta_{\mathrm{bar}}\mathsf{[w_1|\cdots|w_n]}_{\Lambda_0}
 =\sum_{j=0}^n
 \mathsf{[w_1|\cdots|w_j]}_{\Lambda_0}
 \otimes_\Lambda
 \mathsf{[w_{j+1}|\cdots|w_n]}_{\Lambda_0}.
\end{equation}
Here the outside coefficients are the compatible endpoint idempotents
and the empty prefixes or suffixes represent degree-zero bar cells.
The diagonal is a chain map, strictly coassociative and counital
with respect to the bar augmentation.

\begin{theorem}[Explicit coproduct]\label{cup:thm:diagonal}
The shuffle image is stable under deconcatenation. It therefore
inherits a unique coproduct
\[
 \Delta:\mathsf A_\bullet\longrightarrow
 \mathsf A_\bullet\otimes_\Lambda\mathsf A_\bullet
\]
characterized by the commutative square
\begin{equation}\label{cup:eq:diagonal-comparison}
 \boxed{\Delta_{\mathrm{bar}}\kappa
 =(\kappa\otimes_\Lambda\kappa)\Delta.}
\end{equation}
On critical generators,
\begin{equation}\label{cup:eq:diagonal}
 \boxed{\displaystyle
 \Delta(\mathrm C_{q;u,v})
 =\sum_{\substack{0\leq r\leq u,\ 0\leq s\leq v\\q+u-r+s\leq m}}
 (-1)^{s(u-r)}
 \mathrm C_{q+u-r;r,s}\otimes_\Lambda
 \mathrm C_{q+s;u-r,v-s}.}
\end{equation}
The intermediate vertex $q+u-r+s$ must not exceed $m$; with this
restriction both factors are admissible. Their degrees are
$r+s$ and $u+v-r-s$. With
\[
 \mathsf d_\otimes(x\otimes y)
 =\mathsf d_{\mathsf A}x\otimes y
 +(-1)^{|x|}x\otimes\mathsf d_{\mathsf A}y,
\]
the map $\Delta$ is a chain map. Moreover,
\begin{equation}\label{cup:eq:coassociativity}
\begin{gathered}
 (\Delta\otimes_\Lambda\operatorname{id})\Delta
 =(\operatorname{id}\otimes_\Lambda\Delta)\Delta,\\
 (\mu\otimes_\Lambda\operatorname{id})\Delta
 =\operatorname{id}_{\mathsf A}
 =(\operatorname{id}\otimes_\Lambda\mu)\Delta,
\end{gathered}
\end{equation}
under the canonical identifications with $\Lambda\otimes_\Lambda
\mathsf A$ and $\mathsf A\otimes_\Lambda\Lambda$, where
$\mu:\mathsf A_0\to\Lambda$ is the augmentation and is zero
in positive degrees. Thus the Anick resolution is a counital
differential graded $\Lambda$-coring.
\end{theorem}

\begin{proof}
Consider a shuffle word with $u$ downward and $v$ upward steps.
Cut it after a prefix containing $r$ downward and $s$ upward steps.
The intermediate vertex is
\[
 h=q+u-r+s.
\]
The prefix is an admissible word for
$\mathrm C_{q+u-r;r,s}$, with endpoints $q+u,h$, and the
suffix is an admissible word for
$\mathrm C_{q+s;u-r,v-s}$, with endpoints $h,q+v$.
Indeed,
\[
\begin{aligned}
 (q+u-r)+\max(r,s)&=\max(q+u,h)\leq m,\\
 (q+s)+\max(u-r,v-s)&=\max(h,q+v)\leq m.
\end{aligned}
\]
Conversely, two such admissible words concatenate to an admissible
shuffle word. The inversion number satisfies
\begin{equation}\label{cup:eq:inversion-split}
 \boxed{\operatorname{inv}(w'w'')
 =\operatorname{inv}(w')+\operatorname{inv}(w'')+s(u-r),}
\end{equation}
since the last term counts ascending steps of the prefix that precede
descending steps of the suffix. Grouping all cuts in
$\Delta_{\mathrm{bar}}\kappa(\mathrm C_{q;u,v})$ by $r,s$
proves both \eqref{cup:eq:diagonal-comparison} and
\eqref{cup:eq:diagonal}. The compatibility may be displayed as
\[
\begin{array}{ccc}
 \mathsf A_\bullet &\xrightarrow{\ \kappa\ }&\mathsf B^{\Lambda_0}_\bullet\\[2mm]
 {\scriptstyle\Delta}\downarrow &&\downarrow{\scriptstyle\Delta_{\mathrm{bar}}}\\[2mm]
 \mathsf A_\bullet\otimes_\Lambda\mathsf A_\bullet
 &\xrightarrow{\ \kappa\otimes_\Lambda\kappa\ }&
 \mathsf B^{\Lambda_0}_\bullet\otimes_\Lambda
 \mathsf B^{\Lambda_0}_\bullet.
\end{array}
\]

The retraction \eqref{cup:eq:kappa-retraction} implies
\[
 (\mathsf r\otimes_\Lambda\mathsf r)
 (\kappa\otimes_\Lambda\kappa)=\operatorname{id},
\]
so $\kappa\otimes_\Lambda\kappa$ is injective. Since
$\Delta_{\mathrm{bar}}$ and $\kappa$ are chain maps,
\eqref{cup:eq:diagonal-comparison} implies
$\mathsf d_\otimes\Delta=\Delta\mathsf d_{\mathsf A}$.

For coassociativity, split a word into three consecutive blocks
with descending counts $(r_1,r_2,r_3)$ and ascending counts
$(s_1,s_2,s_3)$. Both ways of parenthesizing the cuts have the same
cross-inversion exponent:
\[
 s_1(r_2+r_3)+s_2r_3
 =s_1r_2+(s_1+s_2)r_3.
\]
They therefore give identical terms and signs, proving strict
coassociativity. The left and right counits select, respectively,
$(r,s)=(0,0)$ and $(r,s)=(u,v)$ in \eqref{cup:eq:diagonal}.
Their coefficients are $+1$, and the corresponding endpoint
idempotents act as the identity on $\mathrm C_{q;u,v}$.
\end{proof}

For bimodule cochains $f\in\mathsf C^p$ and
$g\in\mathsf C^r$, define
\begin{equation}\label{cup:eq:product-on-cochains}
 \boxed{f\smile g
 =\mu_\Lambda(f\otimes_\Lambda g)\Delta_{p,r},}
\end{equation}
where $\Delta_{p,r}$ denotes the component with factor degrees $p,r$,
and $\mu_\Lambda:\Lambda\otimes_\Lambda\Lambda\to\Lambda$
is multiplication. We evaluate
$(f\otimes_\Lambda g)(x\otimes_\Lambda y)=
f(x)\otimes_\Lambda g(y)$, with no additional sign in this
evaluation convention. The chain, coassociativity and counit
identities give
\begin{equation}\label{cup:eq:dga}
\begin{aligned}
 \mathsf d^{p+r}(f\smile g)
 &= (\mathsf d^p f)\smile g
 +(-1)^p f\smile(\mathsf d^r g),\\
 (f\smile g)\smile h&=f\smile(g\smile h),
 \qquad \eta_0\smile f=f=f\smile\eta_0.
\end{aligned}
\end{equation}
Thus $\mathsf C^\bullet$ becomes an associative differential graded
algebra whose induced product on $HH^*(\Lambda)$ is the usual
Hochschild cup product.

\subsection{Products of the chosen representatives}

\begin{theorem}[Multiplication constants]\label{cup:thm:products}
For the representatives \eqref{eq:eta-representatives} and
\eqref{eq:theta-representatives},
\begin{align}
 [\eta_a]\smile[\eta_b]
 &=(-1)^{ab}[\eta_{a+b}]
 &&(a+b\leq m-1),\label{cup:eq:etaeta}\\
 [\eta_a]\smile[\theta_b]
 &=[\theta_b]\smile[\eta_a]
 =(-1)^{ab}[\theta_{a+b}]
 &&(a+b\leq m-2),\label{cup:eq:etatheta}\\
 [\theta_a]\smile[\theta_b]&=0.
 \label{cup:eq:thetatheta}
\end{align}
Here $a,b\geq0$ and the representatives are taken only in their
admissible ranges. Products outside the degree bounds vanish,
and $[\eta_0]=1$. For the displayed representatives and our
coproduct, these identities hold already at the cochain level
whenever both sides are defined.
\end{theorem}

\begin{proof}
Only cuts on which both cochains are nonzero contribute.
The possible nonzero pairs in \eqref{cup:eq:diagonal} are
\begin{equation}\label{cup:eq:surviving-cuts}
\begin{array}{c|c|c}
 \text{product}&\text{relevant tensor factor}&\text{coefficient}\\\hline
 \eta_a\smile\eta_b&
 \mathrm D_{q+b}^{(a)}\otimes_\Lambda\mathrm D_{q+a}^{(b)}
 &(-1)^{ab}\\
 \eta_a\smile\theta_b&
 \mathrm D_{q+b}^{(a)}\otimes_\Lambda\mathrm V_{q+a}^{(b)}
 &(-1)^{ab}\\
 \theta_a\smile\eta_b&
 \mathrm V_{q+b}^{(a)}\otimes_\Lambda\mathrm D_{q+a+1}^{(b)}
 &(-1)^{b(a+1)}.
\end{array}
\end{equation}
The total input generator is $\mathrm D_q^{(a+b)}$ in the
first row and $\mathrm V_q^{(a+b)}$ in the second and third.

For $\eta_a\smile\eta_b$ the product of the two values is an
idempotent at vertex $q+a+b$. Its sign is
\[
 ab+a(q+b-1)+b(q+a-1)
 \equiv ab+(a+b)(q-1)\pmod2.
\]
Thus $\eta_a\smile\eta_b=(-1)^{ab}\eta_{a+b}$.
The same calculation applies to $\eta_a\smile\theta_b$,
with the final output an ascending arrow. For the reverse product,
its total exponent satisfies
\[
 b(a+1)+a(q+b-1)+b(q+a)
 \equiv ab+(a+b)(q-1)\pmod2.
\]
Hence both mixed products equal $(-1)^{ab}\theta_{a+b}$.

If both factors are $\theta$-cochains, their nonzero values are
ascending arrows. The only endpoint-compatible product is of the
form $\mathbf b_i\mathbf b_{i+1}=0$. Therefore
$\theta_a\smile\theta_b=0$ already as a cochain. Finally, the
counit identifies $\eta_0$ as the cup unit, and the finite range
of the resolution makes products in higher degrees zero.
\end{proof}

Put $X=[\eta_1]$ and $Y=[\theta_0]$. Iterating the preceding
formulas gives
\begin{equation}\label{cup:eq:powers}
 \boxed{
 X^\ell=(-1)^{\binom\ell2}[\eta_\ell],
 \qquad
 X^\ell Y=(-1)^{\binom\ell2}[\theta_\ell].}
\end{equation}
The signs are units in every field, and the indicated powers are
nonzero precisely within the ranges of the representatives.

\begin{example}[The degree-four product for $m=3$]\label{cup:ex:m3}
Let $\mathrm W=\mathrm C_{1;2,2}=\mathrm D_1^{(2)}$. The only two
admissible shuffle words yield
\[
 \kappa_4(\mathrm W)
 =\mathsf{[\mathbf a_2|\mathbf a_1|\mathbf b_1|\mathbf b_2]}_{\Lambda_0}
 -\mathsf{[\mathbf a_2|\mathbf b_2|\mathbf a_2|\mathbf b_2]}_{\Lambda_0}.
\]
The degree-$(2,2)$ component of the coproduct is
\[
 \Delta_{2,2}(\mathrm W)
 =\mathrm C_{1;2,0}\otimes_\Lambda\mathrm C_{1;0,2}
 -\mathrm D_2^{(1)}\otimes_\Lambda\mathrm D_2^{(1)}.
\]
Since $\eta_1$ vanishes on the two pure generators and
$\eta_1(\mathrm D_2^{(1)})=-e_3$, we obtain
\[
 (\eta_1\smile\eta_1)(\mathrm W)
 =-(-e_3)(-e_3)=-e_3=-\eta_2(\mathrm W).
\]
Consequently, $[\eta_1]^2=-[\eta_2]\ne0$.
\end{example}

\section{The ring presentation}

\begin{theorem}[Recovery of the known ring]\label{cup:thm:ring}
Let $Z_1,\ldots,Z_{m-1},X,Y$ be commuting homogeneous
indeterminates of degrees $0,\ldots,0,2,1$, respectively.
There is an isomorphism of graded $\Bbbk$-algebras
\begin{equation}\label{cup:eq:ring}
 \boxed{
 HH^*(\Lambda)\cong
 \frac{\Bbbk[Z_1,\ldots,Z_{m-1},X,Y]}
 {(Z_iZ_j,\ Z_iX,\ Z_iY,\ X^m,\ Y^2,\ X^{m-1}Y)}.}
\end{equation}
It maps $Z_i$ to the class of the central loop
$\boldsymbol\omega_i$, $X$ to $[\eta_1]$, and $Y$ to
$[\theta_0]$. This is the presentation of
\cite[Theorem~4.2]{LiXuChen}, expressed in our normalization.
\end{theorem}

\begin{proof}
The loops satisfy $\boldsymbol\omega_i\boldsymbol\omega_j=0$.
Moreover, multiplication of an even positive-degree representative
by a central loop produces a cocycle with only loop-valued
coordinates. Such a cocycle is a coboundary by
\eqref{eq:all-even-boundaries}. Multiplication of an odd representative
by a loop is zero already on cochains, since every loop annihilates
positive paths. Thus
\[
 Z_iZ_j=Z_iX=Z_iY=0\quad\text{in }HH^*(\Lambda).
\]
The remaining relations follow from Theorem~\ref{cup:thm:products}
and the degree bound. The surviving products commute in cohomology:
mixed even--odd products agree by \eqref{cup:eq:etatheta}, and all
odd--odd products vanish. Consequently, the assignments in the
statement define an algebra homomorphism
\[
 \Psi:\mathcal R\longrightarrow HH^*(\Lambda),
 \quad
 \mathcal R=
 \frac{\Bbbk[Z_1,\ldots,Z_{m-1},X,Y]}
 {(Z_iZ_j,\ Z_iX,\ Z_iY,\ X^m,\ Y^2,\ X^{m-1}Y)}.
\]
The source has basis
\[
 1,\ Z_1,\ldots,Z_{m-1},\
 X,\ldots,X^{m-1},\
 Y,XY,\ldots,X^{m-2}Y.
\]
Its images are, respectively,
\begin{equation}\label{cup:eq:ring-basis-map}
\begin{array}{c|c}
 \text{monomial}&\text{image under }\Psi\\\hline
 1&1\\
 Z_i&[\boldsymbol\omega_i]\\
 X^\ell&(-1)^{\binom\ell2}[\eta_\ell]
 \quad(1\leq\ell\leq m-1)\\
 X^\ell Y&(-1)^{\binom\ell2}[\theta_\ell]
 \quad(0\leq\ell\leq m-2).
\end{array}
\end{equation}
These images are a basis by Theorem~\ref{thm:HH}, so $\Psi$ is an
isomorphism. No restriction on the characteristic is needed; in
particular, $Y^2=0$ was proved by multiplication, not by dividing
by $2$.
\end{proof}

\begin{remark}[Comparison with the 2015 normalization]
\label{cup:rem:source-normalization}
Under the arrow identification with \cite{LiXuChen}, their
degree-zero generators $x_i^{\mathrm{paper}}$ represent the same
central loops as $Z_i$, whereas their generators of degrees one
and two correspond to
\[
 y_{\mathrm{paper}}\longmapsto -Y,
 \qquad z_{\mathrm{paper}}\longmapsto -X.
\]
The signs reflect their degree-one generators in the $a$-direction,
$g^1=-a$, and the sign of their degree-two comparison generators,
$g^2=-\kappa(\mathrm C_{q;1,1})$. The degree-one representative
supported on the negative $a$-direction and $-\theta_0$ differ by the
inner derivation $-\operatorname{ad}h$, where
\[
 h=\sum_{i=1}^m i e_i,
 \qquad \operatorname{ad}h(\lambda)=h\lambda-\lambda h,
\]
and $\operatorname{ad}h(\mathbf a_i)=\mathbf a_i$,
$\operatorname{ad}h(\mathbf b_i)=-\mathbf b_i$.
These changes of generators account for the different representative
normalizations without changing the abstract ring.
\end{remark}

\section{The Gerstenhaber bracket and a deformation representative}
\label{gerst:section}

The Gerstenhaber structure of the Hochschild cohomology of
Temperley--Lieb algebras was determined in \cite{LyuWu}.  We record it
in our normalization, emphasizing two explicit representatives: a
weight derivation for the degree-one class and a flat deformation for
the degree-two class.  The aim is not a new abstract bracket
presentation, but a comparison with the representatives used in this
paper.  All calculations in this section are valid over the field
$\Bbbk$ without a characteristic restriction.

Put
\[
 X=[\eta_1],\qquad Y=[\theta_0],\qquad |X|=2,\quad |Y|=1,
\]
and denote the classes of the central loops by
$\boldsymbol\omega_i\in HH^0(\Lambda)$.  The cup-product theorem gives
\begin{equation}\label{gerst:ring-relations}
 X^m=Y^2=X^{m-1}Y=0,\qquad
 \boldsymbol\omega_iX=\boldsymbol\omega_iY
 =\boldsymbol\omega_i\boldsymbol\omega_j=0.
\end{equation}
The chosen representatives are normalized by
\begin{equation}\label{gerst:basis-conversion}
 [\eta_r]=(-1)^{\binom r2}X^r,\qquad
 [\theta_r]=(-1)^{\binom r2}X^rY,
\end{equation}
for $0\le r\le m-1$ and $0\le r\le m-2$, respectively;
$\eta_0$ represents the unit.

\subsection{The grading derivation and its action}

We use the standard Gerstenhaber convention on ordinary Hochschild
cochains.  For cochains $F$ and $G$ of degrees $r\ge1$ and $s\ge0$,
write
\begin{equation}\label{gerst:insertion}
 (F\circ G)(x_1,\ldots,x_{r+s-1})
 =\sum_{i=0}^{r-1}(-1)^{i(s-1)}
 F(x_1,\ldots,x_i,G(x_{i+1},\ldots,x_{i+s}),
 x_{i+s+1},\ldots,x_{r+s-1}).
\end{equation}
For $s=0$ the inserted value is an algebra element, and for $r=0$
we set $F\circ G=0$.  The bracket is
\begin{equation}\label{gerst:bar-bracket}
 [F,G]=F\circ G-(-1)^{(r-1)(s-1)}G\circ F.
\end{equation}
In cohomology it satisfies graded antisymmetry and the Leibniz rule:
\begin{align}
 [F,G]&=-(-1)^{(|F|-1)(|G|-1)}[G,F],
 \label{gerst:antisymmetry}\\
 [F,G\smile H]&=[F,G]\smile H
 +(-1)^{(|F|-1)|G|}G\smile[F,H].
 \label{gerst:leibniz}
\end{align}
In particular, for a derivation $D$ and a cochain $F$ of degree $r$,
\begin{equation}\label{gerst:derivation-bar-action}
 [D,F](x_1,\ldots,x_r)
 =D(F(x_1,\ldots,x_r))
 -\sum_{i=1}^r F(x_1,\ldots,D(x_i),\ldots,x_r).
\end{equation}

Define the \emph{$b$-weight derivation} of $\Lambda$ by
\begin{equation}\label{gerst:derivation}
 D(e_i)=D(\mathbf a_i)=0,\qquad
 D(\mathbf b_i)=\mathbf b_i,\qquad
 D(\boldsymbol\omega_i)=\boldsymbol\omega_i.
\end{equation}
The defining relations are homogeneous for the weight
$\mathrm{wt}(\mathbf a_i)=0$, $\mathrm{wt}(\mathbf b_i)=1$.
Consequently $D$ is a derivation.  Its degree-one class is
$Y=[\theta_0]$.

\begin{lemma}\label{gerst:derivation-lift}
For $n=u+v$, the formula
\begin{equation}\label{gerst:lift-formula}
 \widetilde D_n(x\cdot\mathrm C_{q;u,v}\cdot y)
 =D(x)\cdot\mathrm C_{q;u,v}\cdot y
 +v\,x\cdot\mathrm C_{q;u,v}\cdot y
 +x\cdot\mathrm C_{q;u,v}\cdot D(y)
\end{equation}
defines a $D^e$-operator on $\mathsf A_n$, where
$D^e=D\otimes1+1\otimes D$ acts on $\Lambda^e$.
It commutes with $\mathsf d_{\mathsf A}$ and lifts $D$ through the
augmentation.  For a bimodule cochain $\phi\in\mathsf C^n$, the map
\begin{equation}\label{gerst:lift-action}
 \mathcal L_D\phi=D\circ\phi-\phi\circ\widetilde D_n
\end{equation}
is again a bimodule cochain, commutes with the cochain differential,
and induces the Gerstenhaber action $[Y,-]$ on cohomology.
\end{lemma}

\begin{proof}
The generator $\mathrm C_{q;u,v}$ contains precisely $v$ arrows of
type $\mathbf b$; give it weight $v$.  Each of the four terms of
\eqref{eq:main-differential} has the same total $b$-weight:
\[
 0+v,\qquad (v-1)+1,\qquad 1+(v-1),\qquad v+0.
\]
Thus $\mathsf d_{\mathsf A}\widetilde D
=\widetilde D\mathsf d_{\mathsf A}$.
The endpoint idempotents have weight zero, so the lift is well-defined;
at degree zero it is compatible with multiplication.  The derivation
terms cancel in \eqref{gerst:lift-action} when it is evaluated on
$x\cdot\mathrm C\cdot y$, proving bimodule-linearity.  Commutation
with the cochain differential follows from commutation with
$\mathsf d_{\mathsf A}$.  The identification of this lifted action
with the Gerstenhaber action is the derivation-lifting construction of
\cite{SuarezAlvarez}; on the ordinary bar complex it is exactly
\eqref{gerst:derivation-bar-action}.
\end{proof}

\begin{proposition}\label{gerst:degree-one-action}
For the cohomology classes chosen in
\eqref{eq:eta-representatives}--\eqref{eq:theta-representatives},
\begin{align}
 [Y,[\eta_\ell]]&=-\ell[\eta_\ell]
 &&(0\le\ell\le m-1),\label{gerst:eta-weight}\\
 [Y,[\theta_\ell]]&=-\ell[\theta_\ell]
 &&(0\le\ell\le m-2),\label{gerst:theta-weight}\\
 [Y,\boldsymbol\omega_i]&=\boldsymbol\omega_i,
 \qquad [Y,1]=0.\label{gerst:center-weight}
\end{align}
\end{proposition}

\begin{proof}
On a distinguished generator, $\mathcal L_D$ subtracts the
$b$-weight of the generator from that of the cochain value.
The values of $\eta_\ell$ are vertex idempotents of weight zero,
while $\mathrm D_q^{(\ell)}$ has weight $\ell$.
The values of $\theta_\ell$ are $b$-arrows of weight one,
while $\mathrm V_q^{(\ell)}$ has weight $\ell+1$.
Thus the first two identities hold already at the cochain level.
The last identity is the action of $D$ on the centre.
\end{proof}

\subsection{A flat deformation representing the even generator}

Let $s_i=(-1)^{i-1}$, and let $\tau$ be central.  Define a
$\Bbbk[\tau]$-algebra $\Lambda_\tau$ by the same quiver, with relations
\begin{align}
 \mathbf a_{i+1}\mathbf a_i&=0,
 &\mathbf b_i\mathbf b_{i+1}&=0 &&(1\le i\le m-2),
 \label{gerst:deformation-pure}\\
 \mathbf a_i\mathbf b_i
 &=\mathbf b_{i+1}\mathbf a_{i+1}+\tau s_i e_{i+1}
 &&&&(1\le i\le m-2),
 \label{gerst:deformation-interior}\\
 \mathbf a_{m-1}\mathbf b_{m-1}&=\tau s_{m-1}e_m.
 \label{gerst:deformation-boundary}
\end{align}

\begin{lemma}\label{gerst:flat-family}
These oriented relations form a monic Gr\"obner--Shirshov basis
over $\Bbbk[\tau]$.  Thus the normal paths $\mathcal B$ from
\eqref{eq:basis} are a $\Bbbk[\tau]$-basis of $\Lambda_\tau$.
In particular, $\Lambda_\tau$ is a flat deformation of $\Lambda$.
\end{lemma}

\begin{proof}
The leading paths and possible overlaps remain those of
\eqref{eq:tips} and \eqref{eq:overlap-list}; all reductions decrease
the path order.  Pure $aaa$ and $bbb$ overlaps vanish on both sides.
For the two mixed overlap types, the nontrivial reduction branches
give
\[
\begin{array}{c|c}
 \text{overlap}&\text{normal remainder}\\\hline
 \mathbf a_{i+1}\mathbf a_i\mathbf b_i
 &\tau(s_i+s_{i+1})\mathbf a_{i+1}=0\\
 \mathbf a_i\mathbf b_i\mathbf b_{i+1}
 &\tau(s_i+s_{i+1})\mathbf b_{i+1}=0.
\end{array}
\]
These formulas include $i=m-2$ by the boundary relation;
$s_i+s_{i+1}=0$ in every characteristic.  There are no inclusion
ambiguities.  The Composition--Diamond Lemma over the commutative
coefficient ring $\Bbbk[\tau]$ proves the asserted free normal basis.
\end{proof}

Identify $\Lambda_\tau$ with
$\Lambda\otimes_{\Bbbk}\Bbbk[\tau]$ via this normal basis.
For convenience only, put $\boldsymbol\omega_m=0$ in formulas at the
right endpoint.  The relevant products of positive normal paths are
\begin{align}
 \mathbf b_i\mathbf a_i&=\boldsymbol\omega_i,
 &\mathbf a_i\mathbf b_i&=\boldsymbol\omega_{i+1}+\tau s_i e_{i+1},
 \label{gerst:deformed-arrow-products}\\
 \mathbf a_i\boldsymbol\omega_i&=\tau s_i\mathbf a_i,
 &\boldsymbol\omega_i\mathbf b_i&=\tau s_i\mathbf b_i,
 &\boldsymbol\omega_i^2&=\tau s_i\boldsymbol\omega_i.
 \label{gerst:deformed-loop-products}
\end{align}
All other products of positive normal paths are zero.  Therefore the
transported multiplication is exactly $\mu_\tau=\mu+\tau F$, where
the normalized Hochschild two-cochain $F$ has the nonzero values
\begin{equation}\label{gerst:even-bar-representative}
\begin{aligned}
 F(\mathbf a_i,\mathbf b_i)&=s_i e_{i+1},
 &F(\mathbf a_i,\boldsymbol\omega_i)&=s_i\mathbf a_i,\\
 F(\boldsymbol\omega_i,\mathbf b_i)&=s_i\mathbf b_i,
 &F(\boldsymbol\omega_i,\boldsymbol\omega_i)&=s_i\boldsymbol\omega_i,
 \qquad (1\le i<m).
\end{aligned}
\end{equation}
It vanishes on all other pairs of normal basis paths, including
pairs with a vertex idempotent.

\begin{proposition}\label{gerst:even-self-bracket}
The cochain $F$ is a Hochschild two-cocycle satisfying
$F\circ F=0$, and $[F]=X=[\eta_1]$.  Furthermore,
\[
 [D,F]=-F
\]
as an identity of ordinary Hochschild cochains.
Consequently $[X,X]=0$ and $[Y,X]=-X$ in cohomology.
\end{proposition}

\begin{proof}
Associativity of $\mu+\tau F$ and comparison of the coefficients
of $\tau$ and $\tau^2$ yield the cocycle equation and
$F(F(x,y),z)=F(x,F(y,z))$, i.e. $F\circ F=0$.
The shuffle comparison of \eqref{cup:eq:kappa} gives
\begin{equation}\label{gerst:comparison-degree-two}
 \kappa_2(\mathrm C_{i;1,1})=
 \begin{cases}
  \mathsf{[\mathbf a_i|\mathbf b_i]}_{\Lambda_0}
  -\mathsf{[\mathbf b_{i+1}|\mathbf a_{i+1}]}_{\Lambda_0},
  &1\le i\le m-2,\\
  \mathsf{[\mathbf a_{m-1}|\mathbf b_{m-1}]}_{\Lambda_0},
  &i=m-1.
 \end{cases}
\end{equation}
The cochain $F$ vanishes on pure $aa$, $bb$ and normal $ba$
pairs.  Thus $\kappa_2^*F$ is $s_i e_{i+1}$ on
$\mathrm D_i^{(1)}$ and zero on all other critical generators,
which is exactly $\eta_1$.
Finally, in each nonzero value of $F$, the total $b$-weight of the
inputs exceeds that of the output by one.  Formula
\eqref{gerst:derivation-bar-action} therefore gives $[D,F]=-F$.
\end{proof}

\begin{proposition}\label{gerst:even-center-bracket}
For every $1\le j<m$, the bracket $[X,\boldsymbol\omega_j]$ vanishes
in $HH^1(\Lambda)$.
\end{proposition}

\begin{proof}
By the degree-zero insertion convention,
\[
 [F,\boldsymbol\omega_j](x)
 =F(\boldsymbol\omega_j,x)-F(x,\boldsymbol\omega_j).
\]
Its only nonzero basis values are
$[F,\boldsymbol\omega_j](\mathbf a_j)=-s_j\mathbf a_j$ and
$[F,\boldsymbol\omega_j](\mathbf b_j)=s_j\mathbf b_j$.
These are precisely the values of the inner derivation
$[h_j,-]$, where $h_j=s_j(e_1+\cdots+e_j)$; on all other
basis paths both maps vanish.  Hence the bracket is a coboundary.
\end{proof}

\subsection{The resulting bracket}

\begin{theorem}[Gerstenhaber bracket in our normalization]
\label{gerst:complete-bracket}
Together with graded antisymmetry and the Leibniz identity, the
Gerstenhaber bracket is determined by
\begin{equation}\label{gerst:generator-brackets}
 [Y,X]=-X,\qquad [Y,\boldsymbol\omega_i]=\boldsymbol\omega_i,
 \qquad
 [X,X]=[X,\boldsymbol\omega_i]=[Y,Y]
 = [\boldsymbol\omega_i,\boldsymbol\omega_j]=0.
\end{equation}
The unit has zero bracket with every class.  Equivalently,
\begin{align}
 [X^r,X^s]&=0,\label{gerst:power-even-even}\\
 [X^rY,X^s]&=-sX^{r+s},\label{gerst:power-odd-even}\\
 [X^rY,X^sY]&=(r-s)X^{r+s}Y,\label{gerst:power-odd-odd}\\
 [X^r,\boldsymbol\omega_i]&=0,
 & [X^rY,\boldsymbol\omega_i]&=
 \begin{cases}\boldsymbol\omega_i,&r=0,\\0,&r>0.
 \end{cases}\label{gerst:power-center}
\end{align}
where inadmissible powers vanish as in \eqref{gerst:ring-relations}.
In terms of the representatives of Theorem~\ref{thm:HH},
\begin{align}
 [[\eta_r],[\eta_s]]&=0,\label{gerst:original-even-even}\\
 [[\theta_r],[\eta_s]]&=-s(-1)^{rs}[\eta_{r+s}],
 \label{gerst:original-odd-even}\\
 [[\theta_r],[\theta_s]]&=(r-s)(-1)^{rs}[\theta_{r+s}],
 \label{gerst:original-odd-odd}
\end{align}
and
\begin{equation}\label{gerst:original-center}
 [[\eta_r],\boldsymbol\omega_i]=0,\qquad
 [[\theta_r],\boldsymbol\omega_i]
 =\begin{cases}\boldsymbol\omega_i,&r=0,\\0,&r>0.
 \end{cases}
\end{equation}
The classes on the right are understood to vanish outside their
ranges.  They give this Gerstenhaber structure in the normalization of our
representatives; for its prior determination see \cite{LyuWu}.
\end{theorem}

\begin{proof}
The brackets involving $Y$ follow from
Proposition~\ref{gerst:degree-one-action}; the brackets $[X,X]$ and
$[X,\boldsymbol\omega_i]$ follow from
Propositions~\ref{gerst:even-self-bracket} and
\ref{gerst:even-center-bracket}.  Also $[Y,Y]=0$ is the commutator of
a derivation with itself, while two degree-zero classes have a
bracket of degree $-1$ and hence bracket to zero.
The cup algebra is generated by $X,Y$ and the centre, so the Leibniz
rule yields the power formulas.  For example,
\[
 [X^rY,X^s]=-sX^{r+s},\qquad
 [X^rY,X^sY]=(-s+r)X^{r+s}Y.
\]
The centre formulas follow from the vanishing of positive-degree
products with $\boldsymbol\omega_i$.  Finally, using
\eqref{gerst:basis-conversion} and
$\binom{r+s}{2}=\binom r2+\binom s2+rs$ gives the displayed
formulas for $[\eta_r]$ and $[\theta_r]$.  No division by $2$ is
used, so the argument holds in every characteristic.
\end{proof}

\ResearchPart{From quiver cohomology to Kauffman diagrams}{This part returns to the original algebras $TL_N(\beta)$.
Morita equivalence identifies the basic blocks; central elements,
Chebyshev polynomials, and relative full twists turn their positive
Hochschild generators into operations on Kauffman diagrams. We then
prove blockwise nonvanishing, work through planar examples, and discuss
the finite-lattice interpretation and its limitations.}
	
	\section{Morita transfer and the distinction between cells and operations}\label{tl:transfer}
	\subsection{The full corner of a single block}
	Let $A=TL_N(\beta)$ over $\mathbb C$ with $\beta\ne0$, let
	$B=c_BA$ be a nonsemisimple block, and choose a full basic idempotent
	\[
	f=\varepsilon_1+\cdots+\varepsilon_m\in B,
	\qquad BfB=B,
	\qquad \rho:\Lambda(m)\xrightarrow{\sim}fBf.
	\]
	Here one primitive summand is selected from each isomorphism class of
	indecomposable projectives.  The equality $f=c_B$ is not required and
	usually fails: the regular module of $B$ contains those projectives
	with multiplicities given by the dimensions of their simple heads.
	The corner $fBf$ has dimension $4m-3$; the block itself need not have
	that dimension.
	
	The bimodule resolution transfers to
	\begin{equation}\label{tl:transferred-resolution}
		\widetilde{\mathsf A}_n=Bf\otimes_{\Lambda(m)}\mathsf A_n\otimes_{\Lambda(m)}fB.
	\end{equation}
	The equivalence of bimodule categories sends the augmentation to
	$Bf\otimes_{fBf}fB\cong B$.  It preserves projectives, exactness and,
	for these finite-dimensional algebras, minimal projective resolutions.
	A cochain $\phi:\mathsf A_n\to\Lambda(m)$ transfers by
	\begin{equation}\label{tl:transfer-cochain}
		\widetilde\phi(x\otimes\xi\otimes y)=x\rho(\phi(\xi))y,
		\qquad x\in Bf,\quad y\in fB.
	\end{equation}
	Thus $HH^*(B)\cong HH^*(\Lambda(m))$, compatibly with cup products
	and the Gerstenhaber structure; see \cite{Keller,Transfer}.
	This statement does not by itself specify a convenient formula for
	$\widetilde\phi$ on arbitrary tuples of original TL generators.
	
	In degree zero, even the inclusion of the corner must be interpreted
	correctly.  If $\sum_s x_sy_s=c_B$ with $x_s\in Bf$, $y_s\in fB$,
	then the inverse of restriction on centres is
	\begin{equation}\label{tl:center-transfer}
		Z(fBf)\longrightarrow Z(B),\qquad z\longmapsto\sum_s x_szy_s.
	\end{equation}
	Simply viewing $z\in fBf$ as an element of $B$ need not make it
	central in $B$.
	
	\subsection{Why a formal Anick cell is not a diagram product}
	For $m=3$, the generator $\mathrm C_{1;2,2}=\mathsf{[\mathbf{a}_2|\mathbf{a}_1|\mathbf{b}_1|\mathbf{b}_2]}$ has
	homological degree four and is nonzero in $\mathsf A_4$, although its product
	$\mathbf{a}_2\mathbf{a}_1\mathbf{b}_1\mathbf{b}_2$ in $\Lambda(3)$ is zero.  Its cochain value
	$\eta_2(\mathrm C_{1;2,2})=e_3$ is nonzero.  Erasing the bars would lose the
	relation and syzygy data which the cochain detects.
	
	The return elements $\boldsymbol{\omega}_i=\mathbf{b}_i\mathbf{a}_i$ are likewise not detached TL
	circles.  They satisfy $\boldsymbol{\omega}_i\ne0$ and $\boldsymbol{\omega}_i^2=0$, whereas a detached
	circle is removed with the scalar $\beta$.  A realization of $\boldsymbol{\omega}_i$
	in a corner is normally a linear combination of diagrams whose
	nilpotence results from cancellation.  The quiver index labels
	simple modules inside one block, not physical strands.
	
	\subsection{A concrete failure of the naive cochain extension}
	In $\Lambda(2)$, write $\mathbf{a}=\mathbf{a}_1$, $\mathbf{b}=\mathbf{b}_1$, $\boldsymbol{\omega}=\mathbf{b}\mathbf{a}$, so $\mathbf{a}\mathbf{b}=0$.
	The Anick class $\eta_1$ takes the value $e_2$ on $\mathsf{[\mathbf{a}|\mathbf{b}]}$.
	The ordinary bilinear function $F_0(\mathbf{a},\mathbf{b})=e_2$, zero on all other
	basis pairs, is not a cocycle:
	\[
	(\delta F_0)(\mathbf{a},\mathbf{b},\mathbf{a})=F_0(\mathbf{a},\boldsymbol{\omega})-e_2\mathbf{a}=-\mathbf{a}.
	\]
	The corrected ordinary cochain, already constructed in
	\eqref{gerst:even-bar-representative}, has the values
	\[
	F(\mathbf{a},\mathbf{b})=e_2,\quad F(\mathbf{a},\boldsymbol{\omega})=\mathbf{a},\quad F(\boldsymbol{\omega},\mathbf{b})=\mathbf{b},\quad F(\boldsymbol{\omega},\boldsymbol{\omega})=\boldsymbol{\omega}.
	\]
	These extra values are forced by the extension from formal cells to
	arbitrary inputs.  For a general resolution cochain $\phi$, an ordinary
	representative is obtained by precomposition with a comparison
	$\mathsf B^{\Bbbk}_\bullet\to\mathsf A_\bullet$, for instance the composite of relative
	normalization and the Morse projection.  Direct insertion into the
	slots of an Anick cochain is not defined before such a comparison.
	
	\subsection{Why the low-degree search is sufficient}
	By Theorem~\ref{thm:HH}, a basic nonsemisimple block has one-dimensional
	$HH^1$ and $HH^2$.  Therefore an ordinary derivation and an ordinary
	bilinear cocycle with nonzero classes suffice to identify its positive
	generators.  Their higher cup products then follow from
	Theorem~\ref{cup:thm:ring}.  This motivates looking for $D$ and $\Phi$
	directly in the original diagram algebra, rather than first realizing
	all the projective arrows separately.
	
	For several blocks, global nonvanishing is not enough.  We must prove
	nonvanishing on each block, because a vector $(Y_1,0)$ is nonzero in
	$HH^1(B_1)\oplus HH^1(B_2)$ but misses the second summand.  The
	blockwise arguments in Section~\ref{tl:blocknonzero} are consequently
	an essential part of the construction, not an optional normalization.
	
	\section{A direct low-degree test from the TL presentation}\label{tl:direct}
	Before deriving the closed formula, it is useful to isolate a finite
	linear test which motivated and checked the search.  This test does
	not use primitive idempotents or the quiver resolution as input.
	Throughout this section $N\ge2$. Let $v_i\in A=TL_N(\beta)$ be unknown elements and put
	\begin{align}
		S_i(v)&=v_i u_i+u_i v_i-\beta v_i,\label{tl:S}\\
		T_{ij}(v)&=v_i u_j u_i+u_i v_j u_i+u_i u_j v_i-v_i
		\quad(|i-j|=1),\label{tl:T}\\
		K_{ij}(v)&=v_i u_j+u_i v_j-v_j u_i-u_j v_i
		\quad(j>i+1).\label{tl:K}
	\end{align}
	Collect these into a linear map $\mathcal L_\beta:A^{N-1}\to A^{r_N}$,
	where
	\[
	r_N=(N-1)+2(N-2)+\frac{(N-2)(N-3)}2
	=\frac{N^2+N-4}2.
	\]
	Its matrix in the diagram basis has
	$r_N\mathrm{Cat}_N$ rows and $(N-1)\mathrm{Cat}_N$ columns.
	
	\begin{proposition}\label{tl:direct-test}
		The assignments $D(u_i)=v_i$ define a derivation precisely when
		$\mathcal L_\beta(v)=0$.  Inner derivations are the tuples
		$(hu_i-u_ih)_i$.  Moreover, let $s$ be the vector which is $u_i$
		in the $S_i$ slot and zero in every $T_{ij}$ and $K_{ij}$ slot.
		For the loop-parameter cocycle $\Phi$ defined in
		\eqref{tl:Phi-definition},
		\begin{equation}\label{tl:inhomogeneous-test}
			[\Phi]=0\quad\Longleftrightarrow\quad
			\mathcal L_\beta(v)=s\text{ has a solution}.
		\end{equation}
	\end{proposition}
	\begin{proof}
		The homogeneous equations are precisely the derivatives of all
		relations \eqref{tl:relations}.  They are necessary and sufficient
		for the derivation of the free algebra on the $u_i$ to preserve its
		relation ideal.
		
		For the second assertion work over $R=\mathbb C[t]/(t^2)$, comparing
		the constant algebra $A\otimes R$ and the diagram algebra
		$TL_N(\beta+t)$.  The assignment $u_i\mapsto u_i-tv_i$ preserves
		the constant square relation exactly when
		$S_i(v)=u_i$, and preserves all the other defining relations exactly
		when their corresponding right-hand sides are zero.  It then defines
		an $R$-algebra homomorphism which reduces to the identity and is
		therefore an isomorphism on the finite free diagram modules.  Write it
		as $\operatorname{id}-th$.  Multiplicativity to first order is exactly
		$\Phi=\delta h$.  Conversely, such a trivialization supplies
		$v_i=h(u_i)$ and the displayed equations.
	\end{proof}
	
	Thus the same presentation matrix finds derivations through its kernel
	and tests the particular known two-cocycle through an inhomogeneous
	system.  This is not a computation of every element of $HH^2$; an
	arbitrary degree-two cochain still has to satisfy the next differential.
	
	\begin{example}[A direct nonvanishing test at $\beta=1$]
		For $N\ge3$, the character $\chi(u_i)=1$ is defined on $TL_N(1)$.
		If the inhomogeneous system had a solution, $\chi(S_i(v))=1$
		would give $\chi(v_i)=1$.  Applying $\chi$ to $T_{i,i+1}(v)=0$
		would give $\chi(v_i)+\chi(v_{i+1})=0$, a contradiction over
		$\mathbb C$.  Hence $[\Phi_1]\ne0$ globally.  The analogous test
		at $\beta=-1$ uses $\chi(u_i)=-1$.  These character arguments by
		themselves do not separate several nonsemisimple blocks.
	\end{example}

	\section{Chebyshev polynomials and the central tower}\label{tl:central}
	\subsection{The two polynomial families}
	We use the normalizations
	\begin{equation}\label{tl:cheb-def}
		\begin{gathered}
			P_0(x)=1,\quad P_1(x)=x,\quad P_{r+1}(x)=xP_r(x)-P_{r-1}(x),\\
			\Theta_0(x)=2,\quad\Theta_1(x)=x,\quad
			\Theta_{r+1}(x)=x\Theta_r(x)-\Theta_{r-1}(x).
		\end{gathered}
	\end{equation}
	Thus $P_r(x)=U_r(x/2)$ and $\Theta_r(x)=2T_r(x/2)$ in the
	standard Chebyshev notation.  The symbol $\Theta_r$, rather than
	$C_r$, avoids confusion with the full twists introduced later.
	\[
	\begin{array}{c|c|c}
		r&P_r(x)&\Theta_r(x)\\\hline
		0&1&2\\1&x&x\\2&x^2-1&x^2-2\\
		3&x^3-2x&x^3-3x\\4&x^4-3x^2+1&x^4-4x^2+2
	\end{array}
	\]
	The recurrences and initial values directly imply
	\begin{equation}\label{tl:cheb-laurent}
		\begin{aligned}
			P_r(s+s^{-1})&=s^r+s^{r-2}+\cdots+s^{-r},\\
			\Theta_r(s+s^{-1})&=s^r+s^{-r}.
		\end{aligned}
	\end{equation}
	Multiplication of the first sum by $s+s^{-1}$ followed by subtraction
	of the sum for $r-1$ proves the next case, including its endpoints.
	In particular,
	\begin{equation}\label{tl:cheb-trig}
		P_r(2\cos\theta)=\frac{\sin((r+1)\theta)}{\sin\theta},\qquad
		\Theta_r(2\cos\theta)=2\cos(r\theta).
	\end{equation}
	At the endpoints the polynomials are still defined:
	$P_r(2)=r+1$ and $P_r(-2)=(-1)^r(r+1)$.
	For formulas at $r=0$ set $P_{-1}(x)=0$. We shall also use
	\begin{equation}\label{tl:cheb-extra}
		\Theta_r'(x)=rP_{r-1}(x),\qquad
		2P_r(x)-xP_{r-1}(x)=\Theta_r(x),\qquad
		\sum_{r\ge0}P_r(x)t^r=(1-xt+t^2)^{-1}.
	\end{equation}
	The final equality is an identity of formal power series; no analytic
	convergence is involved.  The first two follow either from
	\eqref{tl:cheb-trig} or by the recurrences.
	
	Write
	\[
	[r]_q=P_{r-1}(\beta)=\frac{q^r-q^{-r}}{q-q^{-1}}
	\quad(\beta=q+q^{-1},\ q\ne\pm1).
	\]
	If $e=\operatorname{ord}(q^2)\ge3$ and $\sigma=q^e$, then
	\begin{equation}\label{tl:root-values}
		[e]_q=0,\quad [e-1]_q=-\sigma,\quad
		[e-2]_q=-\sigma\beta,\quad
		\Theta_{e-1}(\beta)=\sigma\beta.
	\end{equation}
	For the original real family $q=\exp(\sqrt{-1}\pi/e)$ this says
	$[e-1]_q=1$, $[e-2]_q=\beta$ and
	$\Theta_{e-1}(\beta)=-\beta$.
	
	A polynomial in a central element need not be determined by its
	scalar eigenvalues on simple modules.  For example,
	$Q(\lambda1+z)=Q(\lambda)1+Q'(\lambda)z$ when $z^2=0$.
	We therefore establish the following identities in the algebra,
	not only on its simple modules.
	
	\subsection{The recursive central elements}
	Initially regard $\beta$ as an indeterminate and work in the
	universal diagram algebra over $\mathbb Z[\beta]$.  Set
	\begin{equation}\label{tl:Frec}
		\begin{aligned}
			F_0&=\beta1,\qquad F_1=(\beta^2-2)1,\\
			F_n&=\beta F_{n-1}-F_{n-2}-u_{n-1}F_{n-1}-F_{n-1}u_{n-1}
			+\beta F_{n-2}u_{n-1}\quad(n\ge2).
		\end{aligned}
	\end{equation}
	All smaller-level elements are embedded by adding through strands.
	For example,
	\begin{align}
		F_2&=(\beta^3-3\beta)1+(4-\beta^2)u_1,\label{tl:F2}\\
		F_3&=(\beta^4-4\beta^2+2)1
		-\beta(\beta^2-4)(u_1+u_2)
		+(\beta^2-4)(u_1u_2+u_2u_1).\label{tl:F3}
	\end{align}
	These agree with the normalization of Appendix~A of \cite{RSA}.
	We shall identify the whole sequence with the element used there
	after proving its recursive identities.
	
	\begin{lemma}\label{tl:F-identities}
		Put $\kappa=\beta^2-2$ and $\alpha=\beta^2-4$.  Then
		\begin{align}
			F_n&\in Z(A_n),\label{tl:Fcentral}\\
			u_nF_nu_n&=\kappa F_{n-1}u_n\quad(n\ge1),\label{tl:Fsandwich}\\
			F_nu_{n-1}&=F_{n-2}u_{n-1}\quad(n\ge2),\label{tl:Fcap}\\
			F_n^2-\beta F_{n-1}F_n+F_{n-1}^2+\alpha1&=0\quad(n\ge1).
			\label{tl:Fquad}
		\end{align}
		These are polynomial identities and survive every specialization of
		$\beta$.
	\end{lemma}
	\begin{proof}
		We give the details since these identities support the later uniform
		formula.  For the sandwich relation put $u=u_n$, $v=u_{n-1}$,
		$A=F_{n-1}$ and $B=F_{n-2}$.  Both $A$ and $B$ commute with $u$
		by distant commutation.  Hence
		\[
		\begin{aligned}
			uF_nu&=u(\beta A-B-vA-Av+\beta Bv)u\\
			&=\beta^2Au-\beta Bu-2Au+\beta Bu=\kappa Au.
		\end{aligned}
		\]
		The case $n=1$ is immediate.  This calculation does not presuppose
		centrality of $F_n$.
		
		For centrality use induction.  Commutation with $u_j$ for
		$j\le n-3$ follows term by term.  With $u=u_{n-1}$,
		$A=F_{n-1}$ and $B=F_{n-2}$, we have $[B,u]=0$ and
		\[
		[F_n,u]=\beta[A,u]-[uA+Au,u]
		=\beta(Au-uA)-(Au^2-u^2A)=0.
		\]
		It remains to commute with $v=u_{n-2}$.  Put additionally
		$C=F_{n-3}$; then
		$A=\beta B-C-vB-Bv+\beta Cv$ and $[A,v]=0$.
		Substitution into
		\[
		[F_n,v]=-[B,v]-[u,v]A-A[u,v]+\beta(Buv-vBu)
		\]
		reduces its right-hand side to
		\[
		[u,vBv-\kappa Cv],
		\]
		which is zero by the already proved sandwich identity on the smaller
		level.  To check the reduction without commuting illegal factors, put
		$T=uv-vu$, $L=Bv-vB$.  Then $Tv+vT=\beta T$ and
		\[
		TvB+TBv+vBT+BvT
		=\beta(uvB-Bvu)+uvBv-vBvu+L;
		\]
		these cancel exactly the corresponding terms in the displayed
		commutator.  This completes the centrality induction.
		
		Multiplying \eqref{tl:Frec} by $v=u_{n-1}$ on the right and using
		the sandwich identity gives
		\[
		F_nv=(\beta^2-1)F_{n-2}v-vF_{n-1}v=F_{n-2}v.
		\]
		By centrality the same relation holds on the left.  Consequently,
		for every polynomial $Q$,
		\begin{equation}\label{tl:polycap}
			Q(F_n)u_{n-1}=Q(F_{n-2})u_{n-1}.
		\end{equation}
		
		Finally, the base case of \eqref{tl:Fquad} is the scalar identity
		$(\beta^2-2)^2-\beta^2(\beta^2-2)+\beta^2+\beta^2-4=0$.
		For the induction set
		\[
		A=F_{n-1},\quad B=F_{n-2},\quad v=u_{n-1},\qquad
		S=vA+Av-\beta Bv,
		\]
		so $F_n=\beta A-B-S$, $[A,B]=[B,v]=0$,
		$vAv=\kappa Bv$, $v^2=\beta v$.  Let
		$K=A^2-\beta AB+B^2+\alpha1=0$ by induction.  Expansion yields
		\[
		F_n^2-\beta AF_n+A^2+\alpha1=K+S^2+2BS-\beta SA,
		\]
		where
		\[
		S^2=vA^2v+\beta AvA-2B(vA+Av)-\beta\alpha B^2v.
		\]
		Therefore
		\[
		S^2+2BS-\beta SA
		=vA^2v-\beta vA^2+\beta^2BvA-\beta\kappa B^2v
		=vKv-\beta vK=0.
		\]
		This proves the last identity and the lemma.
	\end{proof}
	
	\subsection{Central characters and identification with the standard element}
	In a generic semisimple specialization, $F_n$ acts on $V_{n,0}$
	by the scalar recurrence with initial values $\beta$, $\beta^2-2$,
	hence by $\Theta_{n+1}(\beta)$.  For $p>0$, the nonzero subspace
	$u_{n-1}V_{n,p}$ is identified, after removing the fixed final arc,
	with $V_{n-2,p-1}$ as an $A_{n-2}$-module.  Equation
	\eqref{tl:Fcap} gives the same central scalar on these spaces.
	Induction yields
	\begin{equation}\label{tl:central-character}
		F_n|_{V_{n,p}}=\Theta_{n-2p+1}(\beta)\operatorname{id}
		=(q^{n-2p+1}+q^{-(n-2p+1)})\operatorname{id}.
	\end{equation}
	This is Proposition~A.2 of \cite{RSA}.  In a generic specialization
	the full list of simple standard modules detects equality of central
	elements.  Their diagrammatic coefficients are polynomials in the
	parameter, so the identification with the element in that appendix
	and formula \eqref{tl:central-character} extend to all
	specializations.  This does not assert scalar action on every
	indecomposable projective; a nilpotent part is precisely what will
	detect nonsemisimplicity.
	
	\section{Partial closure and the Chebyshev transfer identity}\label{tl:closure}
	Let $\tau_n:A_n\to A_{n-1}$ close the last upper and lower endpoints
	outside the rectangle and remove every resulting circle with weight
	$\beta$.  It is the unnormalized partial closure, not a scalar trace.
	It satisfies
	\begin{equation}\label{tl:tau-basic}
		\tau_n(1)=\beta1,\qquad
		\tau_n(axb)=a\tau_n(x)b\quad(a,b\in A_{n-1}),\qquad
		u_nxu_n=\tau_n(x)u_n.
	\end{equation}
	These identities are obtained by stacking the same diagrams in the
	two orders.  For instance $\tau_3(u_2)=1$ while
	$\tau_3(u_1)=\beta u_1$.  Partial closure is not an algebra
	homomorphism, and one cannot replace $\tau_n(Q(F_n))$ by
	$Q(\tau_n(F_n))$.
	
	The map $a\mapsto au_n$ from $A_{n-1}$ into $A_{n+1}$ is injective:
	it adds fixed disjoint arcs on the last two boundary positions of
	each basis diagram.  Comparing \eqref{tl:Fsandwich} and
	\eqref{tl:tau-basic} therefore gives
	\begin{equation}\label{tl:tau-F}
		\tau_n(F_n)=(\beta^2-2)F_{n-1}.
	\end{equation}
	
	\begin{lemma}[Chebyshev transfer under closure]\label{tl:cheb-transfer}
		For all $n\ge1$ and $r\ge0$,
		\begin{equation}\label{tl:tau-P}
			\boxed{\tau_n(P_r(F_n))=\Theta_{r+1}(\beta)P_r(F_{n-1}).}
		\end{equation}
	\end{lemma}
	\begin{proof}
		Put $X=F_n$, $A=F_{n-1}$ and $X'=\beta A-X$.  They commute and
		\eqref{tl:Fquad} gives
		$XX'=A^2+\beta^2-4$, $X+X'=\beta A$.
		Thus
		\[
		\begin{aligned}
			Q(t)&=(1-Xt+t^2)(1-X't+t^2)\\
			&=1-\beta At+(A^2+\beta^2-2)t^2-\beta At^3+t^4.
		\end{aligned}
		\]
		All coefficients of $Q(t)^{-1}$ belong to $A_{n-1}$, so they can
		be moved through $\tau_n$.  Moreover
		$\tau_n(X')=\beta^2 A-(\beta^2-2)A=2A$.
		Using the formal generating series in \eqref{tl:cheb-extra},
		\begin{equation}\label{tl:trace-generating}
			\tau_n((1-Xt+t^2)^{-1})
			=\frac{\beta-2At+\beta t^2}
			{1-\beta At+(A^2+\beta^2-2)t^2-\beta At^3+t^4}.
		\end{equation}
		On the other hand, with $\beta=q+q^{-1}$,
		\[
		\begin{aligned}
			\sum_{r\ge0}\Theta_{r+1}(\beta)P_r(A)t^r
			&=\frac q{1-qAt+q^2t^2}
			+\frac{q^{-1}}{1-q^{-1}At+q^{-2}t^2}\\
			&=\frac{\beta-2At+\beta t^2}{Q(t)}.
		\end{aligned}
		\]
		Comparison of coefficients proves the identity.  It is a polynomial
		identity in the loop parameter, so no semisimplicity or nondegeneracy
		of the Markov form is used.
	\end{proof}
	
	\begin{example}
		For $r=2$, \eqref{tl:Fquad} and \eqref{tl:tau-F} give
		\[
		\tau_n(F_n^2-1)=\beta(\beta^2-3)(F_{n-1}^2-1)
		=\Theta_3(\beta)P_2(F_{n-1}).
		\]
		The factor is not obtained by applying $P_2$ to the right-hand side
		of \eqref{tl:tau-F}; the generating-function argument is what
		produces the correct basis of polynomial expressions.
	\end{example}
	
	\section{A compatible derivation of the diagram tower}\label{tl:derivation}
	\subsection{The central-difference criterion}
	The small example $D(u_1)=0$, $D(u_2)=(1-u_1)u_2$ in $TL_3(1)$
	suggests searching for values $D(u_i)=c_iu_i$, with $c_i\in A_i$.
	The square relation then asks for $u_ic_iu_i=0$.  Distant
	commutation suggests taking $c_i$ central relative to $A_{i-1}$,
	while adjacent relations suggest a telescoping difference.
	These requirements are met by $c_i=H_i-H_{i-1}$.
	
	\begin{lemma}[Central differences]\label{tl:central-differences}
		Suppose $H_n\in Z(A_n)$ satisfy
		\begin{equation}\label{tl:H-conditions}
			\tau_n(H_n)=\beta H_{n-1},\qquad
			(H_n-H_{n-2})u_{n-1}=0.
		\end{equation}
		Then $D_N(u_i)=(H_i-H_{i-1})u_i$ defines a derivation of $A_N$
		for every $N$, and $D_{N+1}|_{A_N}=D_N$.
	\end{lemma}
	\begin{proof}
		With $c_i=H_i-H_{i-1}$, closure gives
		$u_ic_iu_i=\tau_i(c_i)u_i=0$, because
		$\tau_i(H_{i-1})=\beta H_{i-1}$.  Hence
		$D(u_i^2-\beta u_i)=c_iu_i^2+u_ic_iu_i-\beta c_iu_i=0$.
		If $i<j-1$, $D(u_i)\in A_{i+1}$ commutes with $u_j$, and $c_j$
		commutes with $A_{j-1}$, giving preservation of the distant relation.
		
		For adjacent generators put $x=u_i$, $y=u_{i+1}$,
		$A=H_{i-1}$, $B=H_i$, $C=H_{i+1}$.  The legal commutations are:
		$A$ with $x,y$, $B$ with $y$, and $C$ with $x$.  We do not commute
		$B$ with $x$ or $C$ with $y$.  Direct expansion gives
		\[
		\begin{aligned}
			D(xyx-x)
			&=x(C-B)yx+xy(B-A)x\\
			&=Cxyx-xByx+xyBx-xyAx=(C-A)x=0.
		\end{aligned}
		\]
		Similarly,
		\[
		\begin{aligned}
			D(yxy-y)&=y\bigl((B-A)x+x(C-B)\bigr)y\\
			&=y[B,x]y=Byxy-yxyB=0,
		\end{aligned}
		\]
		where $(C-A)x=0$ is used in the middle line.  All defining relations
		are preserved.  The derivation on the free algebra descends to TL;
		its values do not depend on the final value of $N$.
	\end{proof}
	
	\subsection{The root-of-unity choice}
	Assume $\beta=q+q^{-1}\ne0$, $3\le e=\operatorname{ord}(q^2)<\infty$,
	and $\sigma=q^e$. Define
	\begin{equation}\label{tl:H-definition}
		\boxed{H_n=\frac{1+\sigma^{n+1}P_{e-2}(F_n)}e.}
	\end{equation}
	The identity $\Theta_{e-1}(\beta)=\sigma\beta$ and
	Lemma~\ref{tl:cheb-transfer} imply the first condition in
	\eqref{tl:H-conditions}; \eqref{tl:polycap} and $\sigma^2=1$
	imply the second.  We have therefore proved the explicit formula
	\begin{equation}\label{tl:D-general}
		\boxed{D_N(u_i)=\frac{\sigma^i}{e}
			\bigl(\sigma P_{e-2}(F_i)-P_{e-2}(F_{i-1})\bigr)u_i.}
	\end{equation}
	For $\sigma=-1$ it becomes
	\begin{equation}\label{tl:D-original}
		D_N(u_i)=\frac{(-1)^{i+1}}e
		\bigl(P_{e-2}(F_i)+P_{e-2}(F_{i-1})\bigr)u_i.
	\end{equation}
	The expression $H_n$ is central but need not be an idempotent.
	The scalars in the constant parts of successive $H_n$ cancel in
	the formula for $D$.
	
	The polynomial degree $e-2$ is dictated by closure: it makes its
	scalar multiplier equal to $\sigma\beta$, and the factor
	$\sigma^{n+1}$ corrects the sign between successive levels.  The
	normalization by $e$ fixes the first nonzero value, as follows.
	
	\subsection{The first nonzero level}
	For $n<e$ the algebra $A_n$ is semisimple.  On $V_{n,p}$ write
	$h=n-2p+1$.  If $e\nmid h$, \eqref{tl:cheb-laurent} gives
	\begin{equation}\label{tl:P-at-height}
		P_{e-2}(q^h+q^{-h})=-\sigma^h.
	\end{equation}
	Since $h\equiv n+1\pmod2$, this makes $H_n$ zero on all simple
	modules when $n<e-1$.  At $n=e-1$ the exceptional height is $h=e$,
	where $F_n=2\sigma$ and
	$P_{e-2}(2\sigma)=(e-1)\sigma^{e-2}$; there $H_n=1$.
	It follows that
	\begin{equation}\label{tl:H-initial}
		H_0=\cdots=H_{e-2}=0,\qquad H_{e-1}=f_{e-1},
	\end{equation}
	where $f_{e-1}$ is the Jones--Wenzl projector onto the trivial module
	on $e-1$ strands.  The value $H_0=0$ can also be checked directly
	from $P_{e-2}(\beta)=-\sigma$.
	
	For reference, in our positive loop convention the projectors satisfy
	\begin{equation}\label{tl:JW}
		f_1=1,\qquad
		f_{r+1}=f_r-\frac{[r]_q}{[r+1]_q}f_ru_rf_r,
	\end{equation}
	with $f_r^2=f_r$ and $u_jf_r=f_ru_j=0$ for $j<r$.
	Only $f_{e-1}$ is used here, so all required denominators are
	nonzero.  The singular step to $f_e$ is not used; compare
	\cite{Morrison}, taking account of its loop-sign convention.
	Thus
	\[
	D_e(u_i)=0\ (i<e-1),\qquad
	D_e(u_{e-1})=f_{e-1}u_{e-1}.
	\]
	
	\subsection{An outer class and its persistence along the tower}
	Let $p=f_{e-1}$, $r=f_{e-2}$ and $w=u_{e-1}$.  By
	\eqref{tl:root-values}, the ratio $[e-2]_q/[e-1]_q$ is $\beta$,
	so $p=r-\beta ru_{e-2}r$.  Since $r$ commutes with $w$,
	\[
	wpw=\beta rw-\beta r(wu_{e-2}w)r=0.
	\]
	The element $z=pwp$ is central in $A_e$, has square zero, and is
	nonzero.  To see the last point, its coefficient at the diagram $w$
	is one: a nonidentity term from the left $p$ preserves an upper arc
	among the first $e-1$ endpoints, and a nonidentity term from the right
	$p$ preserves such a lower arc; neither can yield $w$.
	Since $D_e(p)=0$,
	\begin{equation}\label{tl:initial-witness}
		D_e(z)=pD_e(w)p=z\ne0.
	\end{equation}
	An inner derivation annihilates the centre, so $D_e$ is outer.
	
	\begin{lemma}[Restriction through a bimodule retraction]\label{tl:outer-retraction}
		Suppose $\beta\ne0$, $D_N(A_m)\subseteq A_m$, and the restriction
		$D_m$ is outer.  Then $D_N$ is outer.
	\end{lemma}
	\begin{proof}
		Successive partial closures give the $A_m$-bimodule retraction
		\[
		E_{N,m}=\beta^{-(N-m)}\tau_{m+1}\circ\cdots\circ\tau_N:A_N\to A_m.
		\]
		If $D_N=[h,-]$, applying this retraction to $a\in A_m$ yields
		$D_m(a)=E_{N,m}(h)a-aE_{N,m}(h)$, a contradiction.
		No nondegeneracy of a trace form is required.
	\end{proof}
	
	Consequently $D_N$ is zero for $N<e$ and outer for every $N\ge e$.
	This proves global nonvanishing.  Its stronger blockwise version will
	be proved in Section~\ref{tl:blocknonzero}.
	
	\subsection{Words and interval generators}
	On a word the derivation is the explicit Leibniz sum
	\begin{equation}\label{tl:D-word}
		D(u_{i_1}\cdots u_{i_s})=
		\sum_{j=1}^s u_{i_1}\cdots u_{i_{j-1}}D(u_{i_j})
		u_{i_{j+1}}\cdots u_{i_s}.
	\end{equation}
	Preservation of \eqref{tl:relations} proves independence of the chosen
	word for the diagram.  The loop parameter is fixed under $D$:
	$D(\beta^L d)=\beta^LD(d)$.
	For the descending interval word $U_{[i,j]}=u_i u_{i-1}\cdots u_j$
	with $1\le j\le i\le N-1$,
	\begin{equation}\label{tl:D-interval}
		\boxed{D(U_{[i,j]})=(H_i-H_{j-1})U_{[i,j]}.}
	\end{equation}
	Indeed, in the inductive step $u_i$ commutes with $H_{i-1}\in A_{i-1}$
	and with $H_{j-1}$, since these use generators of index at most
	$i-2$.  The two adjacent differences therefore telescope.
	This formula is useful in the interval alphabet of the direct TL
	Gr\"obner--Shirshov notes; it does not identify resolution cells with
	their products.

	\section{Relative full twists and a compact braid formula}\label{tl:braids}
	Define
	\begin{equation}\label{tl:braid-generators}
		t_i=1-qu_i,\qquad t_i^{-1}=1-q^{-1}u_i.
	\end{equation}
	The two expressions are inverse by $u_i^2=(q+q^{-1})u_i$.
	For adjacent $u,v$, the difference between the two braid words is
	\[
	(1-qu)(1-qv)(1-qu)-(1-qv)(1-qu)(1-qv)
	=(-q+q^2\beta-q^3)(u-v)=0.
	\]
	Distant commutation is immediate.  Thus a braid word represents an
	element of TL, with every crossing standing for the explicit linear
	combination \eqref{tl:braid-generators}.
	
	The full twist of \cite[Appendix~A]{RSA} and its relative form are
	\begin{equation}\label{tl:full-relative}
		C_0=C_1=1,\qquad C_n=(t_1\cdots t_{n-1})^n,\qquad
		J_n=C_nC_{n-1}^{-1}\quad(n\ge1).
	\end{equation}
	The previous full twist is embedded on the first $n-1$ strands.
	The braid relations give
	\begin{equation}\label{tl:Jn-palindrome}
		J_n=t_{n-1}\cdots t_2t_1^2t_2\cdots t_{n-1},\qquad
		J_n=t_{n-1}J_{n-1}t_{n-1}.
	\end{equation}
	It is a relative full twist in which the last strand goes around the
	preceding ones and returns to its original endpoint.  Although $C_n$
	is central in $A_n$, $J_n$ need not be; it centralizes $A_{n-1}$.
	
	Set
	\begin{equation}\label{tl:Rn}
		R_n=q^{-n-1}J_n\quad(n\ge1),\qquad R_0=q^{-1}1.
	\end{equation}
	The last expression is a scalar initial convention, not a full twist
	on a negative number of strands.
	
	\begin{proposition}\label{tl:F-R}
		For every $q\ne0$ and every $n\ge0$,
		\begin{equation}\label{tl:F-as-R}
			\boxed{F_n=R_n+R_n^{-1}.}
		\end{equation}
		Consequently, for every $r\ge0$,
		\begin{equation}\label{tl:P-R}
			\boxed{P_r(F_n)=\sum_{k=0}^{r}R_n^{r-2k}.}
		\end{equation}
	\end{proposition}
	\begin{proof}
		First take generic $q$.  The full twist acts on $V_{n,p}$ as
		$q^{2p(n+1-p)}$, whereas $F_n$ acts as
		$q^{n-2p+1}+q^{-(n-2p+1)}$ \cite[Appendix~A]{RSA}.
		The restriction to $A_{n-1}$ is the direct sum of the admissible
		$V_{n-1,p}$ and $V_{n-1,p-1}$
		\cite[Proposition~4.1, Corollary~4.2]{RSA}.
		Therefore $J_n$ acts on these respective summands as
		\[
		q^{2p(n+1-p)-2p(n-p)}=q^{2p},\qquad
		q^{2p(n+1-p)-2(p-1)(n+1-p)}=q^{2(n+1-p)}.
		\]
		If $h=n-2p+1$, the eigenvalues of $R_n$ are $q^{-h}$ and $q^h$.
		Their symmetric sums equal the scalar of $F_n$ on all of $V_{n,p}$.
		All simple standard modules together detect equality in the generic
		semisimple algebra, proving \eqref{tl:F-as-R} there.
		
		In the fixed diagram basis both sides have Laurent-polynomial
		coefficients in $q$.  This follows for the right-hand side from
		\eqref{tl:braid-generators} and for the left-hand side from
		\eqref{tl:Frec}.  Their equality is thus an identity of coefficients,
		which specializes at every nonzero $q$, including roots of unity.
		This step retains the nilpotent parts: equality on simple modules
		only at a singular parameter would not suffice.  The case $n=0$
		is the definition of $F_0$.  Finally, \eqref{tl:P-R} follows from
		the Laurent identity \eqref{tl:cheb-laurent}, since an invertible
		$R_n$ commutes with its inverse.
	\end{proof}
	
	\begin{example}[Two strands and the normalization]
		Here
		\[
		J_2=1+q^2(q-q^{-1})u_1,\qquad
		J_2^{-1}=1-q^{-2}(q-q^{-1})u_1.
		\]
		Hence
		\[
		q^{-3}J_2+q^3J_2^{-1}
		=(q^3+q^{-3})1-(q-q^{-1})^2u_1=F_2,
		\]
		as in \eqref{tl:F2}.  This check fixes the power $q^{-n-1}$.
	\end{example}
	
	\begin{remark}[Why one full twist loses information]
		The assertion is not that $F_n$ is a polynomial in $C_n$ alone.
		For $n=4$ and $q=\exp(\sqrt{-1}\pi/6)$, $C_4$ has the same
		scalar $1$ on $V_{4,0}$ and $V_{4,2}$, whereas $F_4$ has values
		$-\sqrt3$ and $\sqrt3$.  This example is noted in
		\cite[Appendix~A]{RSA}.  A polynomial in $C_4$ cannot distinguish
		these values; the neighboring full twist in $J_n$ supplies the
		missing information.
	\end{remark}
	
	For the root-of-unity derivation define
	$S_n=\sum_{k=0}^{e-2}R_n^{e-2-2k}$.
	Then \eqref{tl:D-general} becomes
	\begin{equation}\label{tl:D-braid}
		\boxed{D_N(u_i)=\frac{\sigma^i}{e}
			(\sigma S_i-S_{i-1})u_i.}
	\end{equation}
	For $\sigma=-1$ this is $(-1)^{i+1}(S_i+S_{i-1})u_i/e$.
	There are at most $2(e-1)$ braid expressions before combining equal
	terms, independent of $N$; their lengths and their expanded planar
	supports may still grow with $N$.  The finite count is a compact
	notation, not a constant-time complexity claim.
	
	For example, at $e=4$, $\beta=\sqrt2$,
	\begin{equation}\label{tl:D-e4-braid}
		D(u_i)=\frac{(-1)^{i+1}}4
		(R_i^2+R_i^{-2}+R_{i-1}^2+R_{i-1}^{-2}+2)u_i.
	\end{equation}
	At $e=5$, the single-level sum is
	$S_n=R_n^3+R_n+R_n^{-1}+R_n^{-3}$; one must still use the correct
	$\sigma$ for the chosen root.

	\begin{figure}[htbp]\centering\begin{tikzpicture}
			\begin{scope}[xshift=0cm,x=.72cm,y=.53cm]
				\draw[line width=.7pt] (0,4) -- (0,3);
				\draw[line width=.8pt] (2,4) .. controls (2,3.6) and (1,3.4) .. (1,3);
				\draw[white,line width=3.5pt] (1,4) .. controls (1,3.6) and (2,3.4) .. (2,3);
				\draw[line width=.8pt] (1,4) .. controls (1,3.6) and (2,3.4) .. (2,3);
				\node[anchor=west,font=\scriptsize] at (2.4,3.5) {$t_2$};
				\draw[line width=.7pt] (2,3) -- (2,2);
				\draw[line width=.8pt] (1,3) .. controls (1,2.6) and (0,2.4) .. (0,2);
				\draw[white,line width=3.5pt] (0,3) .. controls (0,2.6) and (1,2.4) .. (1,2);
				\draw[line width=.8pt] (0,3) .. controls (0,2.6) and (1,2.4) .. (1,2);
				\node[anchor=west,font=\scriptsize] at (2.4,2.5) {$t_1$};
				\draw[line width=.7pt] (2,2) -- (2,1);
				\draw[line width=.8pt] (1,2) .. controls (1,1.6) and (0,1.4) .. (0,1);
				\draw[white,line width=3.5pt] (0,2) .. controls (0,1.6) and (1,1.4) .. (1,1);
				\draw[line width=.8pt] (0,2) .. controls (0,1.6) and (1,1.4) .. (1,1);
				\node[anchor=west,font=\scriptsize] at (2.4,1.5) {$t_1$};
				\draw[line width=.7pt] (0,1) -- (0,0);
				\draw[line width=.8pt] (2,1) .. controls (2,0.6) and (1,0.4) .. (1,0);
				\draw[white,line width=3.5pt] (1,1) .. controls (1,0.6) and (2,0.4) .. (2,0);
				\draw[line width=.8pt] (1,1) .. controls (1,0.6) and (2,0.4) .. (2,0);
				\node[anchor=west,font=\scriptsize] at (2.4,0.5) {$t_2$};
				\node[font=\scriptsize] at (0,4.4) {$1$};
				\node[font=\scriptsize] at (1,4.4) {$2$};
				\node[font=\scriptsize] at (2,4.4) {$3$};
				\node at (1.0,-.55) {$J_3$};
			\end{scope}
			\begin{scope}[xshift=5.3cm,x=.72cm,y=.53cm]
				\draw[line width=.7pt] (0,6) -- (0,5);
				\draw[line width=.7pt] (1,6) -- (1,5);
				\draw[line width=.8pt] (3,6) .. controls (3,5.6) and (2,5.4) .. (2,5);
				\draw[white,line width=3.5pt] (2,6) .. controls (2,5.6) and (3,5.4) .. (3,5);
				\draw[line width=.8pt] (2,6) .. controls (2,5.6) and (3,5.4) .. (3,5);
				\node[anchor=west,font=\scriptsize] at (3.4,5.5) {$t_3$};
				\draw[line width=.7pt] (0,5) -- (0,4);
				\draw[line width=.7pt] (3,5) -- (3,4);
				\draw[line width=.8pt] (2,5) .. controls (2,4.6) and (1,4.4) .. (1,4);
				\draw[white,line width=3.5pt] (1,5) .. controls (1,4.6) and (2,4.4) .. (2,4);
				\draw[line width=.8pt] (1,5) .. controls (1,4.6) and (2,4.4) .. (2,4);
				\node[anchor=west,font=\scriptsize] at (3.4,4.5) {$t_2$};
				\draw[line width=.7pt] (2,4) -- (2,3);
				\draw[line width=.7pt] (3,4) -- (3,3);
				\draw[line width=.8pt] (1,4) .. controls (1,3.6) and (0,3.4) .. (0,3);
				\draw[white,line width=3.5pt] (0,4) .. controls (0,3.6) and (1,3.4) .. (1,3);
				\draw[line width=.8pt] (0,4) .. controls (0,3.6) and (1,3.4) .. (1,3);
				\node[anchor=west,font=\scriptsize] at (3.4,3.5) {$t_1$};
				\draw[line width=.7pt] (2,3) -- (2,2);
				\draw[line width=.7pt] (3,3) -- (3,2);
				\draw[line width=.8pt] (1,3) .. controls (1,2.6) and (0,2.4) .. (0,2);
				\draw[white,line width=3.5pt] (0,3) .. controls (0,2.6) and (1,2.4) .. (1,2);
				\draw[line width=.8pt] (0,3) .. controls (0,2.6) and (1,2.4) .. (1,2);
				\node[anchor=west,font=\scriptsize] at (3.4,2.5) {$t_1$};
				\draw[line width=.7pt] (0,2) -- (0,1);
				\draw[line width=.7pt] (3,2) -- (3,1);
				\draw[line width=.8pt] (2,2) .. controls (2,1.6) and (1,1.4) .. (1,1);
				\draw[white,line width=3.5pt] (1,2) .. controls (1,1.6) and (2,1.4) .. (2,1);
				\draw[line width=.8pt] (1,2) .. controls (1,1.6) and (2,1.4) .. (2,1);
				\node[anchor=west,font=\scriptsize] at (3.4,1.5) {$t_2$};
				\draw[line width=.7pt] (0,1) -- (0,0);
				\draw[line width=.7pt] (1,1) -- (1,0);
				\draw[line width=.8pt] (3,1) .. controls (3,0.6) and (2,0.4) .. (2,0);
				\draw[white,line width=3.5pt] (2,1) .. controls (2,0.6) and (3,0.4) .. (3,0);
				\draw[line width=.8pt] (2,1) .. controls (2,0.6) and (3,0.4) .. (3,0);
				\node[anchor=west,font=\scriptsize] at (3.4,0.5) {$t_3$};
				\node[font=\scriptsize] at (0,6.4) {$1$};
				\node[font=\scriptsize] at (1,6.4) {$2$};
				\node[font=\scriptsize] at (2,6.4) {$3$};
				\node[font=\scriptsize] at (3,6.4) {$4$};
				\node at (1.5,-.55) {$J_4$};
			\end{scope}
		\end{tikzpicture}
		\caption{The relative twists $J_3=t_2t_1^2t_2$ and $J_4=t_3t_2t_1^2t_2t_3$. Each crossing abbreviates $t_i=1-qu_i$; it is not an extra generator of the planar algebra.}
		\label{fig:relative-twists}
	\end{figure}
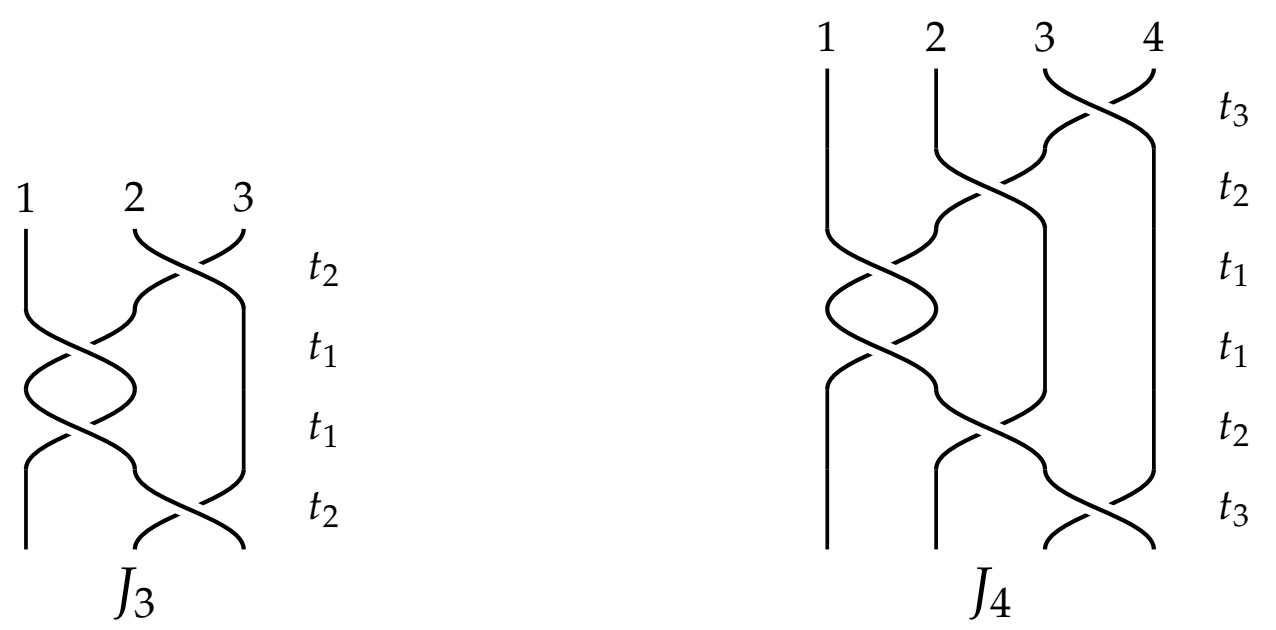
	
	\section{The marked-circle two-cocycle}\label{tl:phi}
	On the vector space of loop-free diagrams, vary the loop weight but
	keep the basis fixed.  If stacking $d_1,d_2$ produces $L$ circles and
	the remaining diagram $d_0$, then
	\[
	\mu_{\beta+t}(d_1,d_2)=(\beta+t)^L d_0.
	\]
	Its linear coefficient is the ordinary bilinear cochain
	\begin{equation}\label{tl:Phi-definition}
		\boxed{\Phi_\beta(d_1,d_2)=
			\begin{cases}L\beta^{L-1}d_0,&L>0,\\0,&L=0.\end{cases}}
	\end{equation}
	The value is the sum over choices of one marked internal circle: the
	marked circle has weight one and the others retain weight $\beta$.
	This is a bilinear operation, not a single diagram assigned to a class.
	
	\begin{proposition}\label{tl:Phi-cocycle}
		The cochain $\Phi_\beta$ is a normalized Hochschild two-cocycle for
		all parameters.  For a central block idempotent $c_B$, its restriction
		\begin{equation}\label{tl:Phi-block}
			\Phi_B(a,b)=c_B\Phi_\beta(a,b)c_B\qquad(a,b\in B)
		\end{equation}
		is a two-cocycle of $B$.
	\end{proposition}
	\begin{proof}
		Associativity of $\mu_{\beta+t}$ holds identically on diagrams.
		The coefficient of $t$ is
		\[
		a\Phi(b,c)-\Phi(ab,c)+\Phi(a,bc)-\Phi(a,b)c=0.
		\]
		Stacking with the identity creates no circles, giving normalization.
		Because $B=c_BA$ is an ideal with identity $c_B$, the Hochschild
	cochain equation restricts to $B$, and central projection of its
	values preserves that equation. This proves the assertion for $B$.
	\end{proof}
	
	Already $\Phi(u_i,u_i)=u_i$ is nonzero as a cochain; this alone says
	nothing about its class.  The blockwise nonvanishing proof below is
	therefore necessary.  Higher coefficients
	\[
	\Phi_r(d_1,d_2)=\binom Lr\beta^{L-r}d_0\quad(L\ge r)
	\]
	remain bilinear, and are zero when $L<r$.  Their index is the power
	of the deformation parameter, not their cohomological degree.
	For instance associativity to second order gives
	$\delta\Phi_2=\Phi\circ\Phi$, not a claim that $\Phi_2$ is an
	independent four-cocycle.  Higher even Hochschild classes will be
	represented instead by cup powers with $2r$ input slots.
	
	\section{Nonvanishing on every nonsemisimple block}\label{tl:blocknonzero}
	Throughout this section $\beta=q+q^{-1}\ne0$,
	$e=\operatorname{ord}(q^2)\ge3$, and the field is $\mathbb C$.
	We use the block structure recalled in Section~\ref{tl:prelim} and
	the centre and cohomology of $\Lambda(m)$ proved in Part~I.
	The proof of nonvanishing is organized to keep the derivation and
	the deformation arguments logically separate.
	
	\subsection{Colliding central characters have different slopes}
	On a standard module with height $h_p=N-2p+1$, define the polynomial
	\[
	f_p(b)=\Theta_{h_p}(b).
	\]
	This is the scalar of the universal central element $F_N(b)$ on the
	universal standard module.  It does not require a choice of a local
	branch $q(b)$.  Differentiating by \eqref{tl:cheb-extra} gives
	\begin{equation}\label{tl:character-slope}
		f_p'(\beta)=h_pP_{h_p-1}(\beta)=h_p[h_p]_q.
	\end{equation}
	Two adjacent labels in a noncritical reflection block have positive
	heights $h,h'$ with $h+h'=2ke$, $k\ge1$, and $e\nmid h$.
	At the special parameter their central characters coincide, but
	$[h']_q=-[h]_q$, so
	\begin{equation}\label{tl:slope-gap}
		\boxed{\Theta_h'(\beta)-\Theta_{h'}'(\beta)
			=2ke[h]_q\ne0.}
	\end{equation}
	The equality of characters and the inequality of slopes are the
	spectral input in both arguments.  An assertion merely that the
	spectrum changes would not by itself exclude a trivial deformation;
	the trace argument below supplies that additional step.
	
	\subsection{The central nilpotent witness on a block}
	Fix a nonsemisimple block $B=c_BA_N$ with common central character
	$\lambda_B=q^{h_p}+q^{-h_p}$.  Put
	\begin{equation}\label{tl:zB}
		z_B=c_B(F_N-\lambda_B1).
	\end{equation}
	It belongs to the radical of $Z(B)$: it acts as zero on all simples
	in the block.  Since the radical of $Z(\Lambda(m))$ has square zero
	and centres are Morita invariant, $z_B^2=0$.
	
	\begin{lemma}\label{tl:zB-nonzero}
		For every such block $B$, $z_B\ne0$.
	\end{lemma}
	\begin{proof}
		The noncritical value $\lambda_B\ne\pm2$ singles out the block
		among the generalized eigenspaces of $F_N$.  Choose a polynomial
		$h$ with $h(F_N(\beta))=c_B$ by the Chinese remainder theorem for
		its minimal polynomial.  In particular $h(\lambda_B)=1$.
		Suppose $z_B=0$, and form the universal central family
		\[
		K(b)=h(F_N(b))(F_N(b)-\lambda_B1).
		\]
		At $b=\beta+t$ it has the form $tG+O(t^2)$ in the fixed diagram
		basis.  Centrality implies $G\in Z(A_N(\beta))$; the first-order
		change of multiplication contributes nothing because $K(\beta)=0$.
		
		On $V_{N,p}$ for $p\in B$, the element $G$ acts as the scalar
		\[
		\left.\frac{d}{dt}\right|_{0}
		h(f_p(\beta+t))(f_p(\beta+t)-\lambda_B)=f_p'(\beta).
		\]
		Passing to the simple heads, these scalars would have to coincide:
		the centre of a single finite-dimensional complex block is local and
		has one scalar character on its simples.  For adjacent labels this
		contradicts \eqref{tl:slope-gap}.  Thus $z_B\ne0$.
	\end{proof}
	
	\subsection{A formula for the action on the central element}
	\begin{proposition}\label{tl:DF}
		For every $n\ge0$,
		\begin{equation}\label{tl:DF-formula}
			\boxed{D_n(F_n)=g_n(F_n),\qquad
				g_n(x)=\frac{x-\sigma^{n+1}\Theta_{e-1}(x)}e.}
		\end{equation}
	\end{proposition}
	\begin{proof}
		We record the full algebraic induction.  Set
		$Q(x)=\Theta_{e-1}(x)$, $U(x)=P_{e-2}(x)$,
		$V(x)=P_{e-1}(x)$.
		If commuting $X,Y$ satisfy
		$X^2-\beta XY+Y^2+\beta^2-4=0$, then
		\begin{equation}\label{tl:transport}
			Q(X)=\sigma\bigl(\beta V(Y)-XU(Y)\bigr).
		\end{equation}
		To prove this as a polynomial identity, in the corresponding quotient
		ring set
		\[
		s=\frac{X-q^{-1}Y}{q-q^{-1}},\qquad
		s^{-1}=\frac{qY-X}{q-q^{-1}}.
		\]
		Their product is one, and $Y=s+s^{-1}$,
		$X=qs+q^{-1}s^{-1}$.  Formula \eqref{tl:cheb-laurent} and
		$q^e=\sigma$ identify both sides of \eqref{tl:transport} with
		$\sigma(q^{-1}s^{e-1}+qs^{-(e-1)})$.
		This is an identity in the quotient ring, so it remains valid on
		nondiagonalizable algebra elements.
		
		The base cases $n=0,1$ follow from
		$\Theta_{e-1}(\beta)=\sigma\beta$ and
		$\Theta_{e-1}(\beta^2-2)=\beta^2-2$.
		For the induction put
		\[
		A=F_{n-1},\quad B=F_{n-2},\quad C=F_n,\quad u=u_{n-1},
		\qquad T=C+B-\beta A,\quad s_0=\sigma^{n+1}.
		\]
		All of $A,B,C$ commute, and $B$ commutes with $u$.
		The central recurrence is equivalent to
		\begin{equation}\label{tl:anti-A}
			uA+Au=-T+\beta Bu.
		\end{equation}
		Apply \eqref{tl:transport} to $A,B$ and use
		$2V(B)-BU(B)=Q(B)$.  It gives
		\begin{equation}\label{tl:anti-Q}
			uQ(A)+Q(A)u=\sigma\beta Q(B)u+\sigma U(B)T.
		\end{equation}
		Apply the same identity to the pairs $(B,A)$ and $(C,A)$:
		\begin{equation}\label{tl:sum-Q}
			Q(C)+Q(B)=\sigma\beta Q(A)-\sigma U(A)T.
		\end{equation}
		
		With $K=H_{n-1}$ and $L=H_{n-2}$, differentiating the recurrence
		\eqref{tl:Frec} yields
		\[
		\begin{aligned}
			D(C)={}&\beta D(A)-D(B)-uD(A)-D(A)u+\beta D(B)u\\
			&+(K-L)T.
		\end{aligned}
		\]
		By induction,
		\[
		eD(A)=A-\sigma s_0Q(A),\quad eD(B)=B-s_0Q(B),\quad
		e(K-L)=\sigma s_0U(A)-s_0U(B).
		\]
		Substituting these and using \eqref{tl:anti-Q} cancels all terms
		containing $Q(B)u$ and $U(B)T$.  What remains is
		\[
		eD(C)=C+s_0\bigl(-\sigma\beta Q(A)+Q(B)+\sigma U(A)T\bigr)
		=C-s_0Q(C)
		\]
		by \eqref{tl:sum-Q}, proving the formula.
	\end{proof}
	
	\begin{theorem}[Blockwise non-innerness and normalization]\label{tl:block-D}
		For each nonsemisimple block $B$,
		\begin{equation}\label{tl:DzB}
			D_N(z_B)=z_B\ne0.
		\end{equation}
		In particular $D_B$ is outer.  Under the Morita identification with
		Part~I, normalized by the degree-one weight derivation there,
		$[D_B]=Y_B$.
	\end{theorem}
	\begin{proof}
		For a noncritical height $h$ we have
		\[
		Q(\lambda_B)=\sigma^h\lambda_B=\sigma^{N+1}\lambda_B,
		\qquad U(\lambda_B)=-\sigma^h=-\sigma^{N+1}.
		\]
		Since $Q'=(e-1)U$, it follows that
		\[
		g_N(\lambda_B)=0,\qquad g_N'(\lambda_B)=1.
		\]
		In the block, $c_BF_N=\lambda_Bc_B+z_B$ and $z_B^2=0$, whence
		$c_Bg_N(F_N)=z_B$.  A derivation annihilates a central idempotent:
		from $D(c)=2cD(c)$ and $(1-2c)^2=1$ we get $D(c)=0$.
		Therefore Proposition~\ref{tl:DF} gives $D(z_B)=z_B$.
		Lemma~\ref{tl:zB-nonzero} now excludes innerness.
		
		The class $Y$ of Part~I acts as the identity on the radical of the
		centre, while inner derivations act trivially there.
		Since $HH^1(B)=\mathbb CY_B$, the equality on the nonzero $z_B$
		fixes the scalar coefficient at one.  In particular $D_B$ acts as
		the identity on the entire radical of $Z(B)$, not merely on $z_B$.
	\end{proof}
	
	\subsection{Central nilpotents as sums of commutators}
	In $\Lambda(m)$ put $\boldsymbol{\omega}_m=0$.  The arrow products
	\eqref{eq:arrow-products} give
	\[
	[\mathbf{b}_j,\mathbf{a}_j]=\boldsymbol{\omega}_j-\boldsymbol{\omega}_{j+1},\qquad
	\boxed{\boldsymbol{\omega}_i=\sum_{j=i}^{m-1}[\mathbf{b}_j,\mathbf{a}_j].}
	\]
	Consequently, in any representation of the constant algebra
	$\Lambda(m)\otimes R$ on a finite free module $M$ over a commutative
	$\mathbb C$-algebra $R$, every $\boldsymbol{\omega}_i$ has zero matrix trace.
	If a central element has the form
	For $R=\mathbb C[t]/(t^2)$, comparing the coefficients of
$1$ and $t$ in the centrality equations gives
$Z(\Lambda(m)\otimes R)=Z(\Lambda(m))\otimes R$. Therefore
\begin{equation}\label{tl:central-over-R}
		K=\alpha1+\sum_i\nu_i \boldsymbol{\omega}_i,
	\end{equation}
	then
	\begin{equation}\label{tl:trace-test}
		\operatorname{tr}_M(K)=\alpha\operatorname{rank}_R M.
	\end{equation}
	If $K$ acts by a scalar $f$ and the rank is positive, characteristic
	zero gives $f=\alpha$.  Thus the scalar of one central element of a
	constant $\Lambda(m)\otimes R$ must be the same on any such module
	on which it acts scalarly.
	
	The commutator expression is essential here.  A nilpotent matrix over
	$\mathbb C[t]/(t^2)$ need not have zero trace: the one-by-one matrix
	$(t)$ is already a counterexample.  We use actual sums of
	commutators, not that incorrect inference from nilpotence.
	
	\subsection{Lifting the block to first order}
	Work over $R=\mathbb C[t]/(t^2)$ in $TL_N(\beta+t)$.  The central
	idempotent $c_B$ lifts to a central idempotent $c_B(t)$, producing a
	finite free block deformation $B_R$.
	For a concrete construction, the polynomial identity
	\begin{equation}\label{tl:F-annihilator}
		\prod_{p=0}^{\lfloor N/2\rfloor}
		\bigl(F_N(b)-\Theta_{N-2p+1}(b)1\bigr)=0
	\end{equation}
	holds generically on all simple standard modules, hence identically
	in the diagram basis.  At $b=\beta+t$, group the factors according
	to their constant characters.  The group for $\lambda_B$ is
	coprime to the other groups in $R[x]$, since their constant
	characters differ.  The Chinese remainder theorem produces
	$c_B(t)$.  On the standard modules belonging to $B$ it acts as one,
	and on the others as zero.  This separates the block from other
	blocks; it does not attempt to split the colliding characters
	inside the same block.
	
	The standard link-state bases are independent of $t$, and their action
	matrices are polynomial in the loop parameter. Hence
	$V_{N,p}(\beta+t)$ for $p\in B$ are finite free $R$-modules
	supported on $B_R$.

	We now identify the first-order class of the lifted block explicitly.
	Use the fixed diagram basis to regard the deformed multiplication on
	$A\otimes R$ as $\mu_t=\mu+t\Phi_\beta$. Choose an $R$-linear
	isomorphism
	\[
	 T_t:B\otimes_{\mathbb C}R\longrightarrow B_R,
	 \qquad T_t(a)=a+t k(a),
	\]
	reducing to the natural inclusion of $B$ in $A$ modulo $t$ and
	chosen with $T_t(c_B)=c_B(t)$; here $k:B\to A$ is a
	$\mathbb C$-linear map, not necessarily valued in $B$.
	Such an isomorphism exists because $B_R$ is a direct summand of the
	finite free $R$-module $A\otimes R$ and has reduction $B$.
	Let $\mu_t^B$ be the multiplication on $B\otimes R$ transported by
	$T_t$. Expanding
	$T_t^{-1}(T_t(a)*_t T_t(b))$ to first order gives
	\[
	 \mu_t^B(a,b)=ab+t\Psi(a,b),\qquad
	 \Psi(a,b)=\Phi_\beta(a,b)+a k(b)+k(a)b-k(ab).
	\]
	The expression for $\Psi$ takes values in $B$, since the left-hand
	side does. Put $\xi(a)=c_Bk(a)c_B$. Multiplying the equality by
	$c_B$ on both sides and using $a=c_Ba=ac_B$ for $a\in B$ yields
	\begin{equation}\label{tl:block-deformation-comparison}
	 \boxed{
	 \Psi(a,b)=\Phi_B(a,b)+a\xi(b)-\xi(ab)+\xi(a)b
	 =\Phi_B(a,b)+(\delta_{\mathrm{bar}}\xi)(a,b).}
	\end{equation}
	Consequently $[\Psi]=[\Phi_B]\in HH^2(B)$, independently
	of the chosen $R$-linear identification. This also agrees with
	Morita invariance for infinitesimal deformations \cite{Transfer}.
	
	\begin{theorem}[Blockwise nontriviality of the marked-circle class]
		\label{tl:block-Phi}
		For every nonsemisimple block $B$ under the present nonzero-parameter
		hypotheses, $[\Phi_B]\ne0$ in $HH^2(B)$.
	\end{theorem}
	\begin{proof}
		Suppose $[\Phi_B]=0$.  The first-order deformation is then
		$R$-algebra isomorphic to the constant algebra $B\otimes R$, by an
		isomorphism reducing to the identity modulo $t$.
		Choose a full basic idempotent $f\in B$ with $fBf\cong\Lambda(m)$;
		only its existence is needed.  Transfer the standard $B_R$-modules
		through the supposed trivialization and apply $f$.  The resulting
		modules $M_p$ over $\Lambda(m)\otimes R$ are free of positive ranks.
		Indeed, the image of an idempotent is a projective direct summand of
		a finite free module over the local ring $R$, and fullness ensures
		that its reduction is not zero.
		
		The transferred image of $c_B(t)F_N(\beta+t)$ in this constant
		corner is central, so it has the form \eqref{tl:central-over-R}
		with coefficients in $R$.  On $M_p$ it still acts as the scalar
		$f_p(\beta+t)$, because trivialization and passage to the corner
		preserve scalar action.  Equation \eqref{tl:trace-test} forces
		\[
		f_p(\beta+t)=\alpha(t)\pmod{t^2}
		\quad\text{for every }p\in B.
		\]
		All derivatives $f_p'(\beta)$ would therefore be equal inside the
		block.  This contradicts \eqref{tl:slope-gap} for any adjacent pair.
		Thus the deformation class cannot be zero.
	\end{proof}
	
	\begin{remark}[Normalization versus nonvanishing]
		Write $[\Phi_B]=b_BX_B$.  Theorem~\ref{tl:block-Phi} proves
		$b_B\ne0$, not a numerical formula for $b_B$ relative to a
		particular arrow normalization.  The difference of the slopes in
		\eqref{tl:slope-gap} is a nonvanishing witness, not automatically
		the coordinate $b_B$.  Numerical normalization is unnecessary for a
		choice of ring generators, but is needed to match a previously fixed
		representative $\eta_\ell$ exactly.
	\end{remark}
	
	\section{The diagrammatic cohomology ring and bracket}\label{tl:diagramring}
	For a nonsemisimple block $B$ of basic size $m$, put
	$x_B=[\Phi_B]$ and $y_B=[D_B]$.  The preceding two theorems,
	Theorem~\ref{cup:thm:ring}, and Morita invariance give
	\begin{equation}\label{tl:block-ring}
		\boxed{HH^*(B)\cong
			\frac{\mathbb C[Z_1,\ldots,Z_{m-1},x_B,y_B]}
			{(Z_iZ_j,\ Z_ix_B,\ Z_iy_B,\ x_B^m,\ y_B^2,\ x_B^{m-1}y_B)}.}
	\end{equation}
	Here $|Z_i|=0$, $|x_B|=2$, $|y_B|=1$.  Rescaling the even
	variable by the nonzero $b_B$ does not change these relations.
	Explicit nonzero classes are represented by
	\begin{equation}\label{tl:all-representatives}
		\begin{gathered}
			\Phi_B^{\smile\ell}\quad(1\le\ell\le m-1),\\
			\Phi_B^{\smile\ell}\smile D_B\quad(0\le\ell\le m-2).
		\end{gathered}
	\end{equation}
	For comparison with the original Anick normalization,
	\[
	[\eta_\ell]=(-1)^{\binom\ell2}b_B^{-\ell}
	[\Phi_B^{\smile\ell}],\qquad
	[\theta_\ell]=(-1)^{\binom\ell2}b_B^{-\ell}
	[\Phi_B^{\smile\ell}\smile D_B].
	\]
	These are equalities of cohomology classes; no assertion is made
	that the Anick cochain and an ordinary bar cochain are the same
	function on the same domain.
	
	For ordinary bar cochains the actual operations are
	\[
	(\Phi_B^{\smile\ell})(a_1,\ldots,a_{2\ell})
	=\Phi_B(a_1,a_2)\cdots\Phi_B(a_{2\ell-1},a_{2\ell}),
	\]
	\[
	(\Phi_B^{\smile\ell}\smile D_B)(a_1,\ldots,a_{2\ell+1})
	=\Phi_B(a_1,a_2)\cdots\Phi_B(a_{2\ell-1},a_{2\ell})D_B(a_{2\ell+1}).
	\]
	All products on the right are diagram stacking.  Formula
	\eqref{tl:D-braid} supplies the derivation values without first
	expanding the central elements completely.
	
	Since the normalization of $D_B$ is fixed by its action on the
	radical of the centre, Theorem~\ref{gerst:complete-bracket} gives
	\begin{equation}\label{tl:bracket}
		[y_B,x_B]=-x_B,\qquad [y_B,Z_i]=Z_i,
		\qquad [x_B,x_B]=[x_B,Z_i]=[y_B,y_B]=[Z_i,Z_j]=0.
	\end{equation}
	All further brackets follow by the graded Leibniz identity.  In
	particular,
	\[
	[x_B^ry_B,x_B^s]=-s x_B^{r+s},\qquad
	[x_B^ry_B,x_B^sy_B]=(r-s)x_B^{r+s}y_B.
	\]
	These are cohomology identities, not a claim that the chosen diagram
	cochains satisfy all of them strictly before quotienting by boundaries.
	
	For the full algebra the result is the product over all blocks.
	A singleton block contributes one copy of $\mathbb C$ in degree zero
	and nothing in positive degree.  Distinct block components have zero
	cup products and brackets with one another.  Central multiplication
	selects each component, giving
	\begin{equation}\label{tl:centre-generation}
		HH^1(A)=Z(A)[D],\qquad HH^2(A)=Z(A)[\Phi].
	\end{equation}
	Thus the positive part is generated by these classes over the centre,
	not necessarily by two elements as a complex vector space.
	The module in \eqref{tl:centre-generation} is generally not free over
	$Z(A)$: central nilpotents annihilate its positive-degree generators.
	
	If $\operatorname{Bl}(A)$ is the set of all blocks and $m_B$ is the size of the
	basic algebra of $B$, with $m_B=1$ for a simple block, then
	\begin{equation}\label{tl:full-hilbert}
		\sum_{r\ge0}\dim HH^r(A)t^r
		=\sum_{B\in\operatorname{Bl}(A)}
		\left(m_B+\sum_{r=1}^{2m_B-2}t^r\right).
	\end{equation}
	The positive representatives do not furnish a planar basis of the
	whole centre.  Its abstract structure is included in
	\eqref{tl:block-ring}, while an explicit corner-to-block realization
	is given by \eqref{tl:center-transfer} once Morita data are chosen.
	One witness $z_B$ does not span the radical of the centre if $m_B>2$.

	\section{Worked diagrammatic examples}\label{tl:examples}
	\subsection{The complete five-dimensional example: \texorpdfstring{$TL_3(1)$}{TL3(1)}}
	Put $u=u_1$, $v=u_2$.  The diagram basis is
	$1,u,v,uv,vu$, with $u^2=u$, $v^2=v$, $uvu=u$, $vuv=v$.
	The central sequence starts with
	\[
	F_0=1,\quad F_1=-1,\quad F_2=-2+3u,
	\qquad H_0=H_1=0,\quad H_2=1-u.
	\]
	Thus
	\begin{equation}\label{tl:TL3D}
		D(u)=0,\qquad D(v)=v-uv.
	\end{equation}
	The values on the full basis are
	\[
	\begin{array}{c|ccccc}
		a&1&u&v&uv&vu\\\hline
		D(a)&0&0&v-uv&0&vu-u.
	\end{array}
	\]
	For instance $D(v^2)=(v-uv)v+v(v-uv)=v-uv$,
	$D(uvu)=u(v-uv)u=0$, and
	$D(vuv)=(v-uv)uv+vu(v-uv)=v-uv$.
	This explicitly checks the nontrivial defining relations.
	
	Choose
	\[
	\varepsilon_1=1-u,\quad\varepsilon_2=u,
	\qquad \mathbf{a}=uv-u,\quad \mathbf{b}=vu-u,
	\qquad \boldsymbol{\omega}=\mathbf{b}\mathbf{a}=u+v-uv-vu.
	\]
	They satisfy
	\[
	\mathbf{a}=\varepsilon_2\mathbf{a}\varepsilon_1,\quad
	\mathbf{b}=\varepsilon_1\mathbf{b}\varepsilon_2,\quad
	\mathbf{a}\mathbf{b}=0,\quad \mathbf{b}\mathbf{a}=\boldsymbol{\omega}\ne0,\quad \boldsymbol{\omega}^2=0.
	\]
	The five elements $\varepsilon_1,\varepsilon_2,\mathbf{a},\mathbf{b},\boldsymbol{\omega}$ form a
	basis: $uv=\mathbf{a}+\varepsilon_2$, $vu=\mathbf{b}+\varepsilon_2$ and
	$v=\boldsymbol{\omega}+\mathbf{a}+\mathbf{b}+\varepsilon_2$.  This gives an explicit isomorphism
	$TL_3(1)\cong\Lambda(2)$, not merely a Morita equivalence.
	Furthermore
	\[
	D(\varepsilon_1)=D(\varepsilon_2)=D(\mathbf{a})=0,
	\qquad D(\mathbf{b})=\mathbf{b},\quad D(\boldsymbol{\omega})=\boldsymbol{\omega}.
	\]
	The derivation is exactly the weight representative $Y$ in these
	coordinates.
	
	For the marked-circle cochain keep the coefficients in $\mathbf{a},\mathbf{b}$ fixed
	while varying the loop weight.  Direct multiplication gives
	\[
	\mu_{1+t}(\mathbf{a},\mathbf{b})=2t\,u,
	\qquad \varepsilon_2\Phi(\mathbf{a},\mathbf{b})\varepsilon_2=2\varepsilon_2.
	\]
	If $\Phi=\delta h$, then $\mathbf{a}\mathbf{b}=0$ and the one-dimensional adjacent
	corner spaces would imply
	$\varepsilon_2(\delta h)(\mathbf{a},\mathbf{b})\varepsilon_2=0$.
	The displayed value therefore also proves nonvanishing directly.
	Under the chosen quiver normalization,
	\[
	X=\tfrac12[\Phi],\qquad Y=[D].
	\]
	The ring is
	\[
	HH^*(TL_3(1))\cong
	\mathbb C[Z,X,Y]/(Z^2,ZX,ZY,X^2,Y^2,XY),
	\]
	where $Z$ denotes the degree-zero class of $\boldsymbol{\omega}$.  We have $[Y,Z]=Z$ and $[Y,X]=-X$.  Nonzero brackets persist although
	all cup products of two positive-degree classes vanish.
	
	The operation $\Phi\smile\Phi$ is nevertheless a nonzero cochain:
	\[
	(\Phi\smile\Phi)(u,u,u,u)=u.
	\]
	Its class is zero because $HH^4(TL_3(1))=0$.  This is an explicit
	warning against equating cochain-level nonvanishing with a nonzero
	cohomology class.
	
	\begin{figure}[htbp]
		\centering
		\begin{tikzpicture}
			\begin{scope}[shift={(0.0,0)}]
				\draw[line width=.8pt] (0.000,1.850) .. controls (0.000,1.200) and (0.000,0.650) .. (0.000,0.000);
				\draw[line width=.8pt] (0.480,1.850) .. controls (0.480,1.200) and (0.480,0.650) .. (0.480,0.000);
				\draw[line width=.8pt] (0.960,1.850) .. controls (0.960,1.200) and (0.960,0.650) .. (0.960,0.000);
				\fill (0.000,1.850) circle (1.2pt);
				\fill (0.480,1.850) circle (1.2pt);
				\fill (0.960,1.850) circle (1.2pt);
				\fill (0.000,0.000) circle (1.2pt);
				\fill (0.480,0.000) circle (1.2pt);
				\fill (0.960,0.000) circle (1.2pt);
				\node at (0.480,-.40) {$1$};
			\end{scope}
			\begin{scope}[shift={(2.5,0)}]
				\draw[line width=.8pt] (0.000,1.850) .. controls (0.000,1.400) and (0.480,1.400) .. (0.480,1.850);
				\draw[line width=.8pt] (0.960,1.850) .. controls (0.960,1.200) and (0.960,0.650) .. (0.960,0.000);
				\draw[line width=.8pt] (0.000,0.000) .. controls (0.000,0.450) and (0.480,0.450) .. (0.480,0.000);
				\fill (0.000,1.850) circle (1.2pt);
				\fill (0.480,1.850) circle (1.2pt);
				\fill (0.960,1.850) circle (1.2pt);
				\fill (0.000,0.000) circle (1.2pt);
				\fill (0.480,0.000) circle (1.2pt);
				\fill (0.960,0.000) circle (1.2pt);
				\node at (0.480,-.40) {$u$};
			\end{scope}
			\begin{scope}[shift={(5.0,0)}]
				\draw[line width=.8pt] (0.000,1.850) .. controls (0.000,1.200) and (0.000,0.650) .. (0.000,0.000);
				\draw[line width=.8pt] (0.480,1.850) .. controls (0.480,1.400) and (0.960,1.400) .. (0.960,1.850);
				\draw[line width=.8pt] (0.480,0.000) .. controls (0.480,0.450) and (0.960,0.450) .. (0.960,0.000);
				\fill (0.000,1.850) circle (1.2pt);
				\fill (0.480,1.850) circle (1.2pt);
				\fill (0.960,1.850) circle (1.2pt);
				\fill (0.000,0.000) circle (1.2pt);
				\fill (0.480,0.000) circle (1.2pt);
				\fill (0.960,0.000) circle (1.2pt);
				\node at (0.480,-.40) {$v$};
			\end{scope}
			\begin{scope}[shift={(7.5,0)}]
				\draw[line width=.8pt] (0.000,1.850) .. controls (0.000,1.400) and (0.480,1.400) .. (0.480,1.850);
				\draw[line width=.8pt] (0.960,1.850) .. controls (0.960,1.200) and (0.000,0.650) .. (0.000,0.000);
				\draw[line width=.8pt] (0.480,0.000) .. controls (0.480,0.450) and (0.960,0.450) .. (0.960,0.000);
				\fill (0.000,1.850) circle (1.2pt);
				\fill (0.480,1.850) circle (1.2pt);
				\fill (0.960,1.850) circle (1.2pt);
				\fill (0.000,0.000) circle (1.2pt);
				\fill (0.480,0.000) circle (1.2pt);
				\fill (0.960,0.000) circle (1.2pt);
				\node at (0.480,-.40) {$uv$};
			\end{scope}
			\begin{scope}[shift={(10.0,0)}]
				\draw[line width=.8pt] (0.000,1.850) .. controls (0.000,1.200) and (0.960,0.650) .. (0.960,0.000);
				\draw[line width=.8pt] (0.480,1.850) .. controls (0.480,1.400) and (0.960,1.400) .. (0.960,1.850);
				\draw[line width=.8pt] (0.000,0.000) .. controls (0.000,0.450) and (0.480,0.450) .. (0.480,0.000);
				\fill (0.000,1.850) circle (1.2pt);
				\fill (0.480,1.850) circle (1.2pt);
				\fill (0.960,1.850) circle (1.2pt);
				\fill (0.000,0.000) circle (1.2pt);
				\fill (0.480,0.000) circle (1.2pt);
				\fill (0.960,0.000) circle (1.2pt);
				\node at (0.480,-.40) {$vu$};
			\end{scope}
		\end{tikzpicture}
		\caption{The five basis diagrams of $TL_3$. Products are stacked with the first factor above the second. Cochains act on these basis inputs and extend linearly.}
		\label{fig:TL3basis}
	\end{figure}
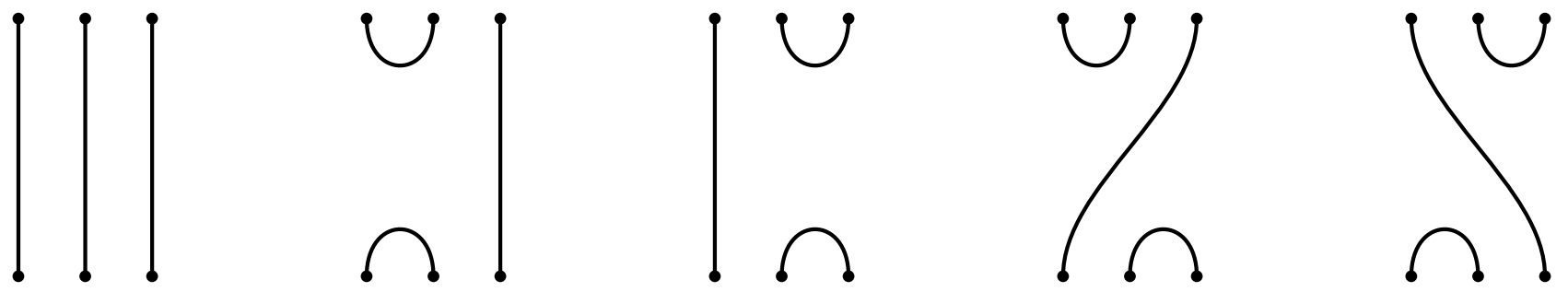
	
	\subsection{A full corner inside \texorpdfstring{$TL_4(\sqrt2)$}{TL4(sqrt2)}}
	Let $s=\sqrt2$ and $u=u_1$, $v=u_2$, $w=u_3$.
	Set
	\begin{equation}\label{tl:TL4projectors}
		p=1-s(u+v)+uv+vu,
		\qquad r=\frac v s-\frac12vuwv.
	\end{equation}
	The letter $r$ avoids using $q$ both for a root parameter and for a
	projector.  Direct diagram multiplication gives
	\[
	p^2=p,\quad r^2=r,\quad pr=rp=0.
	\]
	Put
	\[
	\mathbf{a}=rwp,\qquad \mathbf{b}=pwr,\qquad \boldsymbol{\omega}=\mathbf{b}\mathbf{a}.
	\]
	Then
	\begin{equation}\label{tl:TL4arrows}
		\mathbf{a}\mathbf{b}=0,\qquad \mathbf{b}\mathbf{a}=\frac1s pwp\ne0,\qquad \boldsymbol{\omega}^2=0.
	\end{equation}
	The full corner for the nonsemisimple block has unit $f=p+r$ and
	basis $p,r,\mathbf{a},\mathbf{b},\boldsymbol{\omega}$, hence is isomorphic to $\Lambda(2)$.
	
	For clarity, the other block is not absent.  With $h=uw$, its
	central identity is
	\begin{equation}\label{tl:Hmatrix}
		H=h+vhv-s^{-1}(hv+vh),\qquad c_B=1-H.
	\end{equation}
	It is a four-dimensional simple matrix block; the nonsemisimple
	block $B$ has dimension ten, whereas $fBf$ has dimension five.
	In particular $p+r\ne c_B$.  The element $f$ is full in $B$,
	not the identity of the whole block.
	
	The derivation is
	\[
	D(u)=D(v)=0,\qquad
	D(w)=pw=w-suw-svw+uvw+vuw.
	\]
	Although $p$ involves only $u,v$, the projector $r$ also contains
	$w$, so its derivative must be checked.  Using $D(w)=pw$ and $up=vp=0$,
	the $vuwv$ term in $r$ has zero derivative, so $D(p)=D(r)=0$.
	It follows that $D(\mathbf{a})=0$, $D(\mathbf{b})=\mathbf{b}$, and $D(\boldsymbol{\omega})=\boldsymbol{\omega}$.
	The direct coefficient computation is
	\begin{equation}\label{tl:TL4Phi}
		r\Phi(\mathbf{a},\mathbf{b})r=2s\,r.
	\end{equation}
	For every ordinary one-cochain $h$, the same corner of
	$(\delta h)(\mathbf{a},\mathbf{b})$ is zero: $\mathbf{a}\mathbf{b}=0$ and both surviving corner
	factors are multiples of $\mathbf{a}$ or $\mathbf{b}$.  Therefore
	\[
	X_B=\frac1{2\sqrt2}[\Phi_B],\qquad Y_B=[D_B].
	\]
	The full cohomology ring is the product of the $\Lambda(2)$ ring
	with a copy of $\mathbb C$ in degree zero.  Equations
	\eqref{tl:TL4projectors}--\eqref{tl:TL4Phi} are fixed-normalization
	identities; changing the arrow scales changes the coordinate of $\Phi$.
	
	\subsection{The first nonzero square: \texorpdfstring{$TL_6(1)$}{TL6(1)}}
	The height list is $7,5,3,1$.  The nonsemisimple block has labels
	$\{0,1,3\}$ and $m=3$; $\{2\}$ is critical and simple.
	Thus
	\[
	\sum_j\dim HH^j(TL_6(1))t^j=4+t+t^2+t^3+t^4.
	\]
	On the nonsemisimple block $[\Phi_B]^2\ne0$.  In the normalized
	quiver coordinates, the sign $X^2=-[\eta_2]$ comes from
	\[
	\Delta_{2,2}(\mathrm C_{1;2,2})
	=\mathrm C_{1;2,0}\otimes_\Lambda\mathrm C_{1;0,2}
	-\mathrm D_2^{(1)}\otimes_\Lambda\mathrm D_2^{(1)}.
	\]
	The first summand evaluates to zero under $\eta_1\otimes\eta_1$;
	the second gives $-e_3$.  A single cut in the middle of the word
	would miss the nonzero class.
	
	\subsection{Two independent block directions: \texorpdfstring{$TL_6(\sqrt2)$}{TL6(sqrt2)}}
	There are two nonsemisimple blocks,
	\[
	B_+=\{0,3\},\quad\lambda_+=\sqrt2,
	\qquad B_-=\{1,2\},\quad\lambda_-=-\sqrt2,
	\]
	both of basic size two.  Their standard-character derivatives are
	\[
	\begin{array}{c|c|c|c}
		\text{block}&p&h_p&f_p'(\sqrt2)\\\hline
		B_+&0&7&-7\\B_+&3&1&1\\
		B_-&1&5&-5\\B_-&2&3&3.
	\end{array}
	\]
	Each pair has distinct slopes, proving nonvanishing of both
	restrictions of $\Phi$.  If $F=F_6$, then
	\[
	c_+=\tfrac12+\frac{3F}{4\sqrt2}-\frac{F^3}{8\sqrt2},\qquad
	c_-=1-c_+,
	\]
	\[
	z_+=c_+(F-\sqrt2),\qquad z_-=c_-(F+\sqrt2),\qquad
	D(F)=\frac{F^3-2F}{4}.
	\]
	The cubic projectors impose both value and derivative conditions at
	$\pm\sqrt2$; linear interpolation of values alone would not remove
	nilpotent parts.  The theorem gives $z_\pm\ne0$, $z_\pm^2=0$,
	$D(z_\pm)=z_\pm$.  Here
	\[
	\sum_j\dim HH^j(TL_6(\sqrt2))t^j=4+2t+2t^2.
	\]
	One global derivation and one global two-cocycle generate the two
	separate directions only after multiplication by $c_+$ and $c_-$.
	
	\subsection{Eight and ten strands: block sizes versus algebra dimensions}
	The following data follow from the reflection rule and
	the standard-module dimensions; no additional matrix computation is needed.
	\[
	\begin{array}{c|c|c|c}
		N&e&\text{blocks of labels }p&\text{basic sizes}\\\hline
		4&4&\{0,1\},\{2\}&2,1\\
		6&3&\{0,1,3\},\{2\}&3,1\\
		6&4&\{0,3\},\{1,2\}&2,2\\
		8&4&\{0,1,4\},\{2,3\}&3,2\\
		10&3&\{0,2,3,5\},\{1\},\{4\}&4,1,1.
	\end{array}
	\]
	At $N=8,e=4$ the Hilbert series is
	$5+2t+2t^2+t^3+t^4$.  The square of $[\Phi]$ is nonzero on the
	size-three block and zero on the size-two block.
	For the block $\{0,1,4\}$,
	\[
	\begin{array}{c|ccc}
		p&0&1&4\\\hline
		\dim V_{8,p}&1&7&14\\
		\dim L_{8,p}&1&6&8\\
		\dim P_{8,p}&8&21&14
	\end{array}
	\]
	so the block dimension is $246$, the direct sum of one projective of
	each type has dimension $43$, and the basic corner has dimension $9$.
	These three numbers describe different objects.
	
	For $TL_{10}(1)$ the nonsemisimple block has $m=4$ and
	$[\Phi_B]^3\ne0$ in degree six.  Its dimension data are
	\[
	\begin{array}{c|rrrr}
		p&0&2&3&5\\\hline
		\dim V_{10,p}&1&35&75&42\\
		\dim L_{10,p}&1&34&41&1\\
		\dim P_{10,p}&36&110&117&42.
	\end{array}
	\]
	Consequently
	\[
	\dim B=8615,\qquad
	\dim(P_{10,0}\oplus P_{10,2}\oplus P_{10,3}\oplus P_{10,5})=305,
	\qquad \dim\Lambda(4)=13.
	\]
	The full algebra has dimension $16796$, not $305$.
	The two critical singleton blocks contribute the remaining
	$9^2+90^2=8181$ dimensions.  Its Hochschild series is
	\[
	6+t+t^2+t^3+t^4+t^5+t^6.
	\]
	At this parameter the central projector of the nonsemisimple block is
	\[
	c_B=1-\frac{(F_{10}-1)^2}{9}.
	\]
	The root $1$ is treated to order two and the simple critical value
	$-2$ to order one.  This supplies block selection without the full
	collection of primitive projectors.
	
	\section{A worked visual calculus for the diagrammatic generators}
\label{vis:calculus}
The figures in this section describe operations on the fixed Kauffman
basis. They are not pictures of the abstract isomorphism type of the
cohomology ring. Every input is specified, and all coefficients are
retained. A plain equality denotes an equality of cochains or algebra
elements; square brackets are used when equality holds only in
cohomology. Crossings in the earlier relative-twist figures are the
abbreviations $t_i^{\pm1}=1-q^{\pm1}u_i$, not new TL basis elements.

\subsection{The derivation on three and four strands}
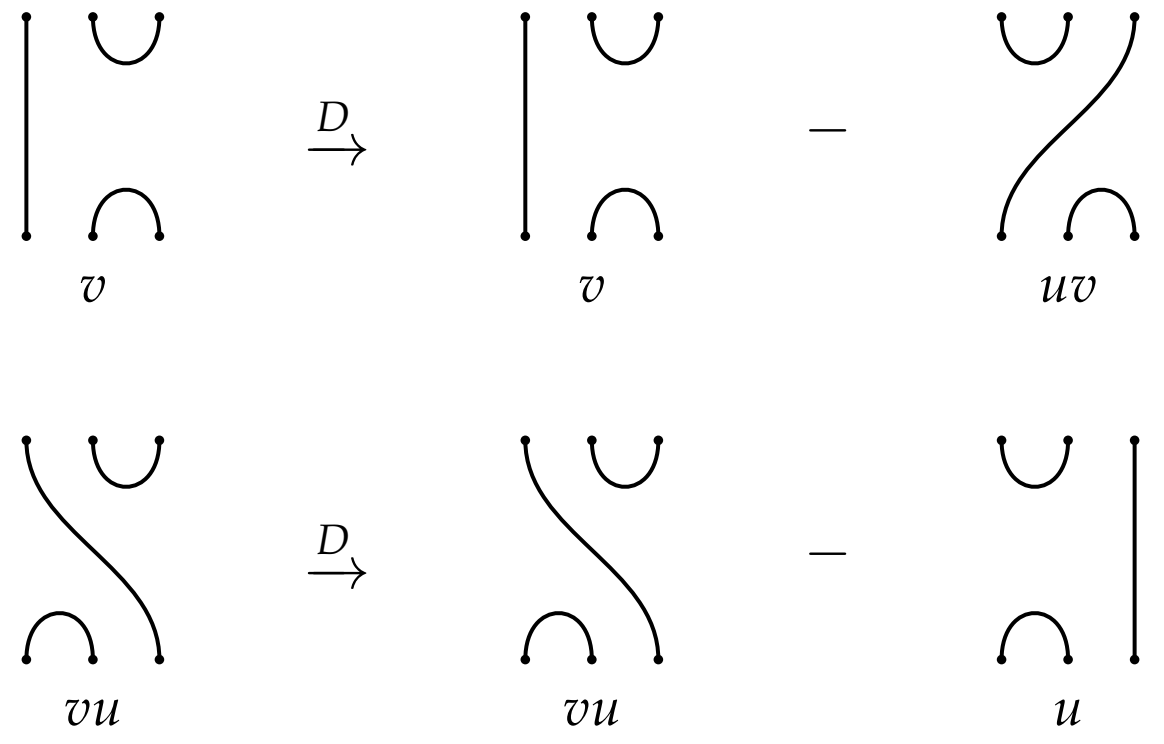
\begin{figure}[htbp]
\centering
\begin{tikzpicture}[scale=1]
\begin{scope}[shift={(0,2.8)}]
\draw[line width=.75pt] (0.000,1.450) .. controls (0.000,0.841) and (0.000,0.609) .. (0.000,0);
\draw[line width=.75pt] (0.440,1.450) .. controls (0.440,1.040) and (0.880,1.040) .. (0.880,1.450);
\draw[line width=.75pt] (0.440,0) .. controls (0.440,0.410) and (0.880,0.410) .. (0.880,0);
\fill (0.000,0) circle (.9pt);\fill (0.000,1.450) circle (.9pt);
\fill (0.440,0) circle (.9pt);\fill (0.440,1.450) circle (.9pt);
\fill (0.880,0) circle (.9pt);\fill (0.880,1.450) circle (.9pt);
\node[font=\small] at (0.440,-.35) {$v$};
\end{scope}
\node at (2.05,3.5) {$\xrightarrow{D}$};
\begin{scope}[shift={(3.3,2.8)}]
\draw[line width=.75pt] (0.000,1.450) .. controls (0.000,0.841) and (0.000,0.609) .. (0.000,0);
\draw[line width=.75pt] (0.440,1.450) .. controls (0.440,1.040) and (0.880,1.040) .. (0.880,1.450);
\draw[line width=.75pt] (0.440,0) .. controls (0.440,0.410) and (0.880,0.410) .. (0.880,0);
\fill (0.000,0) circle (.9pt);\fill (0.000,1.450) circle (.9pt);
\fill (0.440,0) circle (.9pt);\fill (0.440,1.450) circle (.9pt);
\fill (0.880,0) circle (.9pt);\fill (0.880,1.450) circle (.9pt);
\node[font=\small] at (0.440,-.35) {$v$};
\end{scope}
\node at (5.3,3.5) {$-$};
\begin{scope}[shift={(6.45,2.8)}]
\draw[line width=.75pt] (0.000,1.450) .. controls (0.000,1.040) and (0.440,1.040) .. (0.440,1.450);
\draw[line width=.75pt] (0.880,1.450) .. controls (0.880,0.841) and (0.000,0.609) .. (0.000,0);
\draw[line width=.75pt] (0.440,0) .. controls (0.440,0.410) and (0.880,0.410) .. (0.880,0);
\fill (0.000,0) circle (.9pt);\fill (0.000,1.450) circle (.9pt);
\fill (0.440,0) circle (.9pt);\fill (0.440,1.450) circle (.9pt);
\fill (0.880,0) circle (.9pt);\fill (0.880,1.450) circle (.9pt);
\node[font=\small] at (0.440,-.35) {$uv$};
\end{scope}
\begin{scope}[shift={(0,0)}]
\draw[line width=.75pt] (0.000,1.450) .. controls (0.000,0.841) and (0.880,0.609) .. (0.880,0);
\draw[line width=.75pt] (0.440,1.450) .. controls (0.440,1.040) and (0.880,1.040) .. (0.880,1.450);
\draw[line width=.75pt] (0.000,0) .. controls (0.000,0.410) and (0.440,0.410) .. (0.440,0);
\fill (0.000,0) circle (.9pt);\fill (0.000,1.450) circle (.9pt);
\fill (0.440,0) circle (.9pt);\fill (0.440,1.450) circle (.9pt);
\fill (0.880,0) circle (.9pt);\fill (0.880,1.450) circle (.9pt);
\node[font=\small] at (0.440,-.35) {$vu$};
\end{scope}
\node at (2.05,0.7) {$\xrightarrow{D}$};
\begin{scope}[shift={(3.3,0)}]
\draw[line width=.75pt] (0.000,1.450) .. controls (0.000,0.841) and (0.880,0.609) .. (0.880,0);
\draw[line width=.75pt] (0.440,1.450) .. controls (0.440,1.040) and (0.880,1.040) .. (0.880,1.450);
\draw[line width=.75pt] (0.000,0) .. controls (0.000,0.410) and (0.440,0.410) .. (0.440,0);
\fill (0.000,0) circle (.9pt);\fill (0.000,1.450) circle (.9pt);
\fill (0.440,0) circle (.9pt);\fill (0.440,1.450) circle (.9pt);
\fill (0.880,0) circle (.9pt);\fill (0.880,1.450) circle (.9pt);
\node[font=\small] at (0.440,-.35) {$vu$};
\end{scope}
\node at (5.3,0.7) {$-$};
\begin{scope}[shift={(6.45,0)}]
\draw[line width=.75pt] (0.000,1.450) .. controls (0.000,1.040) and (0.440,1.040) .. (0.440,1.450);
\draw[line width=.75pt] (0.880,1.450) .. controls (0.880,0.841) and (0.880,0.609) .. (0.880,0);
\draw[line width=.75pt] (0.000,0) .. controls (0.000,0.410) and (0.440,0.410) .. (0.440,0);
\fill (0.000,0) circle (.9pt);\fill (0.000,1.450) circle (.9pt);
\fill (0.440,0) circle (.9pt);\fill (0.440,1.450) circle (.9pt);
\fill (0.880,0) circle (.9pt);\fill (0.880,1.450) circle (.9pt);
\node[font=\small] at (0.440,-.35) {$u$};
\end{scope}
\end{tikzpicture}
\caption{The derivation on $TL_3(1)$: $D(v)=v-uv$ and $D(vu)=vu-u$. The input labels name the algebra elements; arrows labelled $D$ are operations, not strands. The first factor of a word is stacked above the second.}
\label{fig:vis-D3}
\end{figure}

In $TL_3(1)$ use $u=u_1$, $v=u_2$. From \eqref{tl:TL3D},
\[
 D(u)=0,\quad D(v)=v-uv,\quad D(uv)=0,\quad D(vu)=vu-u.
\]
The last two values can be checked before making any picture:
$u(v-uv)=0$ and $(v-uv)u=vu-u$, using $u^2=u$ and $uvu=u$.
Figure~\ref{fig:vis-D3} displays the nonzero outputs in the planar basis.

\begin{figure}[htbp]
\centering
\begin{tikzpicture}[scale=1]
\begin{scope}[shift={(0.0,0)}]
\draw[line width=.75pt] (0.000,1.450) .. controls (0.000,1.040) and (0.440,1.040) .. (0.440,1.450);
\draw[line width=.75pt] (0.880,1.450) .. controls (0.880,0.841) and (0.880,0.609) .. (0.880,0);
\draw[line width=.75pt] (0.000,0) .. controls (0.000,0.410) and (0.440,0.410) .. (0.440,0);
\fill (0.000,0) circle (.9pt);\fill (0.000,1.450) circle (.9pt);
\fill (0.440,0) circle (.9pt);\fill (0.440,1.450) circle (.9pt);
\fill (0.880,0) circle (.9pt);\fill (0.880,1.450) circle (.9pt);
\node[font=\small] at (0.440,-.35) {$u$};
\end{scope}
\node at (2.2,.7) {$+$};
\begin{scope}[shift={(2.9,0)}]
\draw[line width=.75pt] (0.000,1.450) .. controls (0.000,0.841) and (0.000,0.609) .. (0.000,0);
\draw[line width=.75pt] (0.440,1.450) .. controls (0.440,1.040) and (0.880,1.040) .. (0.880,1.450);
\draw[line width=.75pt] (0.440,0) .. controls (0.440,0.410) and (0.880,0.410) .. (0.880,0);
\fill (0.000,0) circle (.9pt);\fill (0.000,1.450) circle (.9pt);
\fill (0.440,0) circle (.9pt);\fill (0.440,1.450) circle (.9pt);
\fill (0.880,0) circle (.9pt);\fill (0.880,1.450) circle (.9pt);
\node[font=\small] at (0.440,-.35) {$v$};
\end{scope}
\node at (5.1,.7) {$-$};
\begin{scope}[shift={(5.8,0)}]
\draw[line width=.75pt] (0.000,1.450) .. controls (0.000,1.040) and (0.440,1.040) .. (0.440,1.450);
\draw[line width=.75pt] (0.880,1.450) .. controls (0.880,0.841) and (0.000,0.609) .. (0.000,0);
\draw[line width=.75pt] (0.440,0) .. controls (0.440,0.410) and (0.880,0.410) .. (0.880,0);
\fill (0.000,0) circle (.9pt);\fill (0.000,1.450) circle (.9pt);
\fill (0.440,0) circle (.9pt);\fill (0.440,1.450) circle (.9pt);
\fill (0.880,0) circle (.9pt);\fill (0.880,1.450) circle (.9pt);
\node[font=\small] at (0.440,-.35) {$uv$};
\end{scope}
\node at (7.999999999999999,.7) {$-$};
\begin{scope}[shift={(8.7,0)}]
\draw[line width=.75pt] (0.000,1.450) .. controls (0.000,0.841) and (0.880,0.609) .. (0.880,0);
\draw[line width=.75pt] (0.440,1.450) .. controls (0.440,1.040) and (0.880,1.040) .. (0.880,1.450);
\draw[line width=.75pt] (0.000,0) .. controls (0.000,0.410) and (0.440,0.410) .. (0.440,0);
\fill (0.000,0) circle (.9pt);\fill (0.000,1.450) circle (.9pt);
\fill (0.440,0) circle (.9pt);\fill (0.440,1.450) circle (.9pt);
\fill (0.880,0) circle (.9pt);\fill (0.880,1.450) circle (.9pt);
\node[font=\small] at (0.440,-.35) {$vu$};
\end{scope}
\node at (4.8,-1.05) {$z=u+v-uv-vu,\qquad z^2=0,\qquad D(z)=z$};
\end{tikzpicture}
\caption{A geometric witness of outerness. This linear combination is central and nonzero; it is not a detached circle. Any inner derivation kills it, whereas $D$ fixes it.}
\label{fig:vis-z3}
\end{figure}
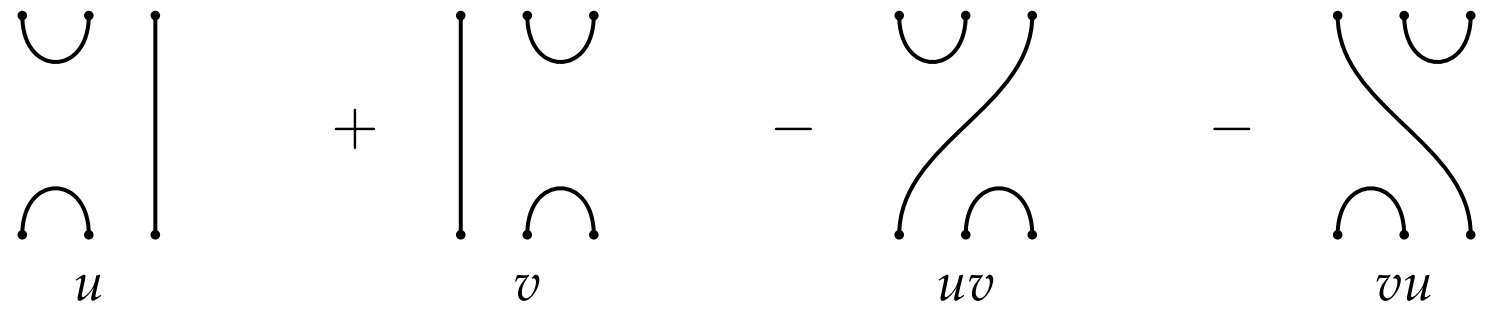

The element $z=u+v-uv-vu$ satisfies
\[
 uz=zu=vz=zv=0,\qquad z^2=0,\qquad D(z)=z.
\]
For example $uz=u+uv-uv-u=0$ and
$D(z)=(v-uv)-(vu-u)=z$. Its four distinct basis terms show $z\ne0$.
Thus Figure~\ref{fig:vis-z3} is a witness of outerness, not a depiction
of a loop of weight $\beta$: a detached ordinary TL circle is scalar,
whereas $z$ is a nonzero central nilpotent.

\begin{figure}[htbp]
\centering
\begin{tikzpicture}[scale=0.98]
\begin{scope}[shift={(0.0,0)}]
\draw[line width=.75pt] (0.000,1.450) .. controls (0.000,0.841) and (0.000,0.609) .. (0.000,0);
\draw[line width=.75pt] (0.350,1.450) .. controls (0.350,0.841) and (0.350,0.609) .. (0.350,0);
\draw[line width=.75pt] (0.700,1.450) .. controls (0.700,1.040) and (1.050,1.040) .. (1.050,1.450);
\draw[line width=.75pt] (0.700,0) .. controls (0.700,0.410) and (1.050,0.410) .. (1.050,0);
\fill (0.000,0) circle (.9pt);\fill (0.000,1.450) circle (.9pt);
\fill (0.350,0) circle (.9pt);\fill (0.350,1.450) circle (.9pt);
\fill (0.700,0) circle (.9pt);\fill (0.700,1.450) circle (.9pt);
\fill (1.050,0) circle (.9pt);\fill (1.050,1.450) circle (.9pt);
\node[font=\small] at (0.525,-.35) {$u_3$};
\end{scope}
\node at (0.525,1.95) {$1$};
\begin{scope}[shift={(2.65,0)}]
\draw[line width=.75pt] (0.000,1.450) .. controls (0.000,1.040) and (0.350,1.040) .. (0.350,1.450);
\draw[line width=.75pt] (0.700,1.450) .. controls (0.700,1.040) and (1.050,1.040) .. (1.050,1.450);
\draw[line width=.75pt] (0.000,0) .. controls (0.000,0.410) and (0.350,0.410) .. (0.350,0);
\draw[line width=.75pt] (0.700,0) .. controls (0.700,0.410) and (1.050,0.410) .. (1.050,0);
\fill (0.000,0) circle (.9pt);\fill (0.000,1.450) circle (.9pt);
\fill (0.350,0) circle (.9pt);\fill (0.350,1.450) circle (.9pt);
\fill (0.700,0) circle (.9pt);\fill (0.700,1.450) circle (.9pt);
\fill (1.050,0) circle (.9pt);\fill (1.050,1.450) circle (.9pt);
\node[font=\small] at (0.525,-.35) {$u_1u_3$};
\end{scope}
\node at (3.175,1.95) {$-\sqrt2$};
\begin{scope}[shift={(5.3,0)}]
\draw[line width=.75pt] (0.000,1.450) .. controls (0.000,0.841) and (0.000,0.609) .. (0.000,0);
\draw[line width=.75pt] (0.350,1.450) .. controls (0.350,1.040) and (0.700,1.040) .. (0.700,1.450);
\draw[line width=.75pt] (1.050,1.450) .. controls (1.050,0.841) and (0.350,0.609) .. (0.350,0);
\draw[line width=.75pt] (0.700,0) .. controls (0.700,0.410) and (1.050,0.410) .. (1.050,0);
\fill (0.000,0) circle (.9pt);\fill (0.000,1.450) circle (.9pt);
\fill (0.350,0) circle (.9pt);\fill (0.350,1.450) circle (.9pt);
\fill (0.700,0) circle (.9pt);\fill (0.700,1.450) circle (.9pt);
\fill (1.050,0) circle (.9pt);\fill (1.050,1.450) circle (.9pt);
\node[font=\small] at (0.525,-.35) {$u_2u_3$};
\end{scope}
\node at (5.825,1.95) {$-\sqrt2$};
\begin{scope}[shift={(7.949999999999999,0)}]
\draw[line width=.75pt] (0.000,1.450) .. controls (0.000,1.040) and (0.350,1.040) .. (0.350,1.450);
\draw[line width=.75pt] (0.700,1.450) .. controls (0.700,0.841) and (0.000,0.609) .. (0.000,0);
\draw[line width=.75pt] (1.050,1.450) .. controls (1.050,0.841) and (0.350,0.609) .. (0.350,0);
\draw[line width=.75pt] (0.700,0) .. controls (0.700,0.410) and (1.050,0.410) .. (1.050,0);
\fill (0.000,0) circle (.9pt);\fill (0.000,1.450) circle (.9pt);
\fill (0.350,0) circle (.9pt);\fill (0.350,1.450) circle (.9pt);
\fill (0.700,0) circle (.9pt);\fill (0.700,1.450) circle (.9pt);
\fill (1.050,0) circle (.9pt);\fill (1.050,1.450) circle (.9pt);
\node[font=\small] at (0.525,-.35) {$u_1u_2u_3$};
\end{scope}
\node at (8.475,1.95) {$1$};
\begin{scope}[shift={(10.6,0)}]
\draw[line width=.75pt] (0.000,1.450) .. controls (0.000,0.680) and (1.050,0.680) .. (1.050,1.450);
\draw[line width=.75pt] (0.350,1.450) .. controls (0.350,1.040) and (0.700,1.040) .. (0.700,1.450);
\draw[line width=.75pt] (0.000,0) .. controls (0.000,0.410) and (0.350,0.410) .. (0.350,0);
\draw[line width=.75pt] (0.700,0) .. controls (0.700,0.410) and (1.050,0.410) .. (1.050,0);
\fill (0.000,0) circle (.9pt);\fill (0.000,1.450) circle (.9pt);
\fill (0.350,0) circle (.9pt);\fill (0.350,1.450) circle (.9pt);
\fill (0.700,0) circle (.9pt);\fill (0.700,1.450) circle (.9pt);
\fill (1.050,0) circle (.9pt);\fill (1.050,1.450) circle (.9pt);
\node[font=\small] at (0.525,-.35) {$u_2u_1u_3$};
\end{scope}
\node at (11.125,1.95) {$1$};
\end{tikzpicture}
\caption{All five planar terms of $D(u_3)$ in $TL_4(\sqrt2)$. The numbers above the diagrams are coefficients in a sum, not labels of separate classes.}
\label{fig:vis-D4}
\end{figure}
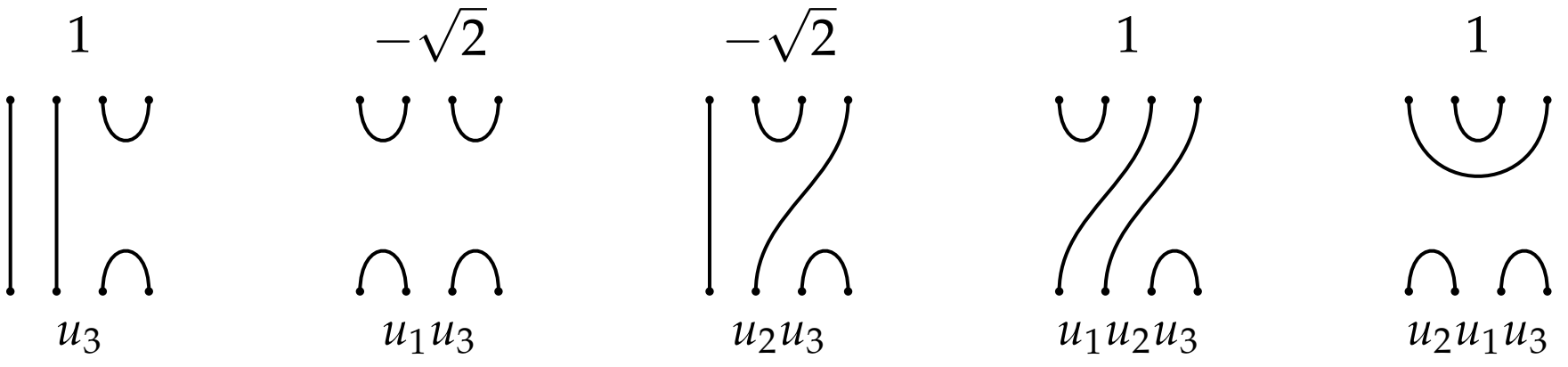

In $TL_4(\sqrt2)$ put $p=1-\sqrt2(u_1+u_2)+u_1u_2+u_2u_1$.
Then $D(u_1)=D(u_2)=0$ and
\begin{equation}\label{vis:D4}
 D(u_3)=pu_3=u_3-\sqrt2u_1u_3-\sqrt2u_2u_3
            +u_1u_2u_3+u_2u_1u_3.
\end{equation}
Figure~\ref{fig:vis-D4} makes the five terms explicit. These are outputs
on one ordinary generator, not five different cohomology classes.
Using $u_1p=u_2p=0$ and $u_3pu_3=0$ gives the checks
\[
 D(u_1u_3)=0,\qquad
 D(u_3^2)=\sqrt2pu_3,\qquad
 D(u_3u_2u_3)=pu_3=D(u_3).
\]
The compact relative-twist expression and its planar expansion therefore
agree on concrete local relations.

\subsection{One marked circle, and the cocycle equation}
For basis diagrams, $\Phi_\beta$ glues two inputs and chooses one of the
internal circles to mark. A marked circle is evaluated at $1$, all
others at $\beta$. This accounts for the coefficient $L\beta^{L-1}$,
including the factor $L$ counting the distinct choices. In particular,
\[
 \Phi_\beta(1,u_i)=0,\qquad\Phi_\beta(u_i,u_i)=u_i.
\]
For the four-strand diagram $h=u_1u_3$, distant commutation gives
\begin{equation}\label{vis:two-loops}
 h^2=\beta^2h,\qquad \Phi_\beta(h,h)=2\beta h.
\end{equation}
There are two marked configurations, not one configuration assigned an
unexplained factor of two (Figures~\ref{fig:vis-mark01} and
\ref{fig:vis-mark}). Marking is part of the multilinear operation; it
is not an additional unmarked element of $TL_N(\beta)$.
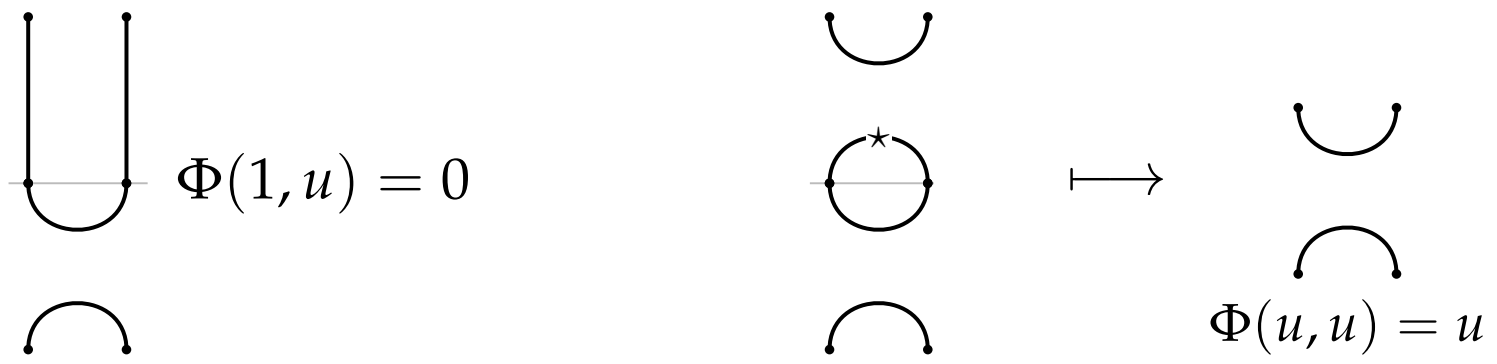
\begin{figure}[htbp]
\centering
\begin{tikzpicture}[scale=1]
\draw[gray!53,line width=.36pt] (-.13,1.10) -- (.79,1.10);
\draw[gray!53,line width=.36pt] (5.17,1.10) -- (5.99,1.10);
\begin{scope}[shift={(0,1.1)}]
\draw[line width=.75pt] (0.000,1.100) .. controls (0.000,0.638) and (0.000,0.462) .. (0.000,0);
\draw[line width=.75pt] (0.650,1.100) .. controls (0.650,0.638) and (0.650,0.462) .. (0.650,0);
\fill (0.000,0) circle (.9pt);\fill (0.000,1.100) circle (.9pt);
\fill (0.650,0) circle (.9pt);\fill (0.650,1.100) circle (.9pt);
\end{scope}
\begin{scope}[shift={(0,0)}]
\draw[line width=.75pt] (0.000,1.100) .. controls (0.000,0.690) and (0.650,0.690) .. (0.650,1.100);
\draw[line width=.75pt] (0.000,0) .. controls (0.000,0.410) and (0.650,0.410) .. (0.650,0);
\fill (0.000,0) circle (.9pt);\fill (0.000,1.100) circle (.9pt);
\fill (0.650,0) circle (.9pt);\fill (0.650,1.100) circle (.9pt);
\end{scope}
\node at (1.95,1.1) {$\Phi(1,u)=0$};
\begin{scope}[shift={(5.3,1.1)}]
\draw[line width=.75pt] (0.000,1.100) .. controls (0.000,0.690) and (0.650,0.690) .. (0.650,1.100);
\draw[line width=.75pt] (0.000,0) .. controls (0.000,0.410) and (0.650,0.410) .. (0.650,0);
\fill (0.000,0) circle (.9pt);\fill (0.000,1.100) circle (.9pt);
\fill (0.650,0) circle (.9pt);\fill (0.650,1.100) circle (.9pt);
\end{scope}
\begin{scope}[shift={(5.3,0)}]
\draw[line width=.75pt] (0.000,1.100) .. controls (0.000,0.690) and (0.650,0.690) .. (0.650,1.100);
\draw[line width=.75pt] (0.000,0) .. controls (0.000,0.410) and (0.650,0.410) .. (0.650,0);
\fill (0.000,0) circle (.9pt);\fill (0.000,1.100) circle (.9pt);
\fill (0.650,0) circle (.9pt);\fill (0.650,1.100) circle (.9pt);
\end{scope}
\node[fill=white,inner sep=.2pt] at (5.625,1.4) {$\scriptstyle\star$};\node at (7.2,1.1) {$\longmapsto$};
\begin{scope}[shift={(8.4,0.5)}]
\draw[line width=.75pt] (0.000,1.100) .. controls (0.000,0.690) and (0.650,0.690) .. (0.650,1.100);
\draw[line width=.75pt] (0.000,0) .. controls (0.000,0.410) and (0.650,0.410) .. (0.650,0);
\fill (0.000,0) circle (.9pt);\fill (0.000,1.100) circle (.9pt);
\fill (0.650,0) circle (.9pt);\fill (0.650,1.100) circle (.9pt);
\node[font=\small] at (0.325,-.35) {$\Phi(u,u)=u$};
\end{scope}
\end{tikzpicture}
\caption{The thin gray separators distinguish the two input diagrams. No internal circle contributes zero; one marked internal circle contributes the resulting diagram once. These are local two-strand patches, with any spectator strands left unchanged.}
\label{fig:vis-mark01}
\end{figure}
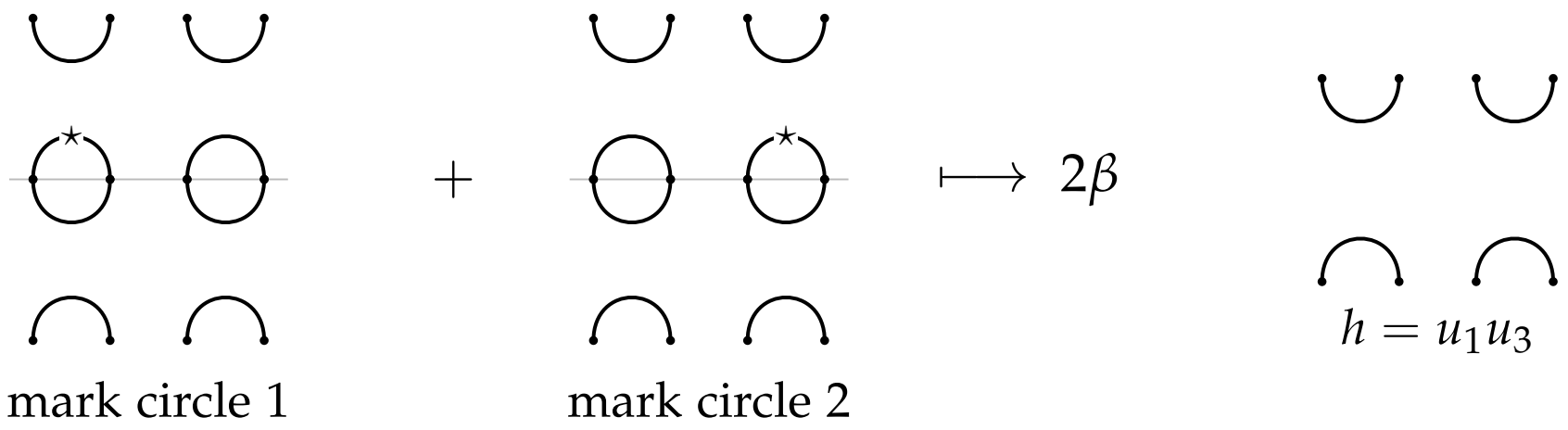
\begin{figure}[htbp]
\centering
\begin{tikzpicture}[scale=1]
\draw[gray!53,line width=.36pt] (-.17,1.15) -- (1.82,1.15);
\draw[gray!53,line width=.36pt] (3.83,1.15) -- (5.82,1.15);
\begin{scope}[shift={(0.0,1.15)}]
\draw[line width=.75pt] (0.000,1.150) .. controls (0.000,0.740) and (0.550,0.740) .. (0.550,1.150);
\draw[line width=.75pt] (1.100,1.150) .. controls (1.100,0.740) and (1.650,0.740) .. (1.650,1.150);
\draw[line width=.75pt] (0.000,0) .. controls (0.000,0.410) and (0.550,0.410) .. (0.550,0);
\draw[line width=.75pt] (1.100,0) .. controls (1.100,0.410) and (1.650,0.410) .. (1.650,0);
\fill (0.000,0) circle (.9pt);\fill (0.000,1.150) circle (.9pt);
\fill (0.550,0) circle (.9pt);\fill (0.550,1.150) circle (.9pt);
\fill (1.100,0) circle (.9pt);\fill (1.100,1.150) circle (.9pt);
\fill (1.650,0) circle (.9pt);\fill (1.650,1.150) circle (.9pt);
\end{scope}
\begin{scope}[shift={(0.0,0)}]
\draw[line width=.75pt] (0.000,1.150) .. controls (0.000,0.740) and (0.550,0.740) .. (0.550,1.150);
\draw[line width=.75pt] (1.100,1.150) .. controls (1.100,0.740) and (1.650,0.740) .. (1.650,1.150);
\draw[line width=.75pt] (0.000,0) .. controls (0.000,0.410) and (0.550,0.410) .. (0.550,0);
\draw[line width=.75pt] (1.100,0) .. controls (1.100,0.410) and (1.650,0.410) .. (1.650,0);
\fill (0.000,0) circle (.9pt);\fill (0.000,1.150) circle (.9pt);
\fill (0.550,0) circle (.9pt);\fill (0.550,1.150) circle (.9pt);
\fill (1.100,0) circle (.9pt);\fill (1.100,1.150) circle (.9pt);
\fill (1.650,0) circle (.9pt);\fill (1.650,1.150) circle (.9pt);
\end{scope}
\node[fill=white,inner sep=.2pt] at (0.275,1.46) {$\scriptstyle\star$};
\node[font=\small] at (0.825,-.42) {mark circle 1};
\begin{scope}[shift={(4.0,1.15)}]
\draw[line width=.75pt] (0.000,1.150) .. controls (0.000,0.740) and (0.550,0.740) .. (0.550,1.150);
\draw[line width=.75pt] (1.100,1.150) .. controls (1.100,0.740) and (1.650,0.740) .. (1.650,1.150);
\draw[line width=.75pt] (0.000,0) .. controls (0.000,0.410) and (0.550,0.410) .. (0.550,0);
\draw[line width=.75pt] (1.100,0) .. controls (1.100,0.410) and (1.650,0.410) .. (1.650,0);
\fill (0.000,0) circle (.9pt);\fill (0.000,1.150) circle (.9pt);
\fill (0.550,0) circle (.9pt);\fill (0.550,1.150) circle (.9pt);
\fill (1.100,0) circle (.9pt);\fill (1.100,1.150) circle (.9pt);
\fill (1.650,0) circle (.9pt);\fill (1.650,1.150) circle (.9pt);
\end{scope}
\begin{scope}[shift={(4.0,0)}]
\draw[line width=.75pt] (0.000,1.150) .. controls (0.000,0.740) and (0.550,0.740) .. (0.550,1.150);
\draw[line width=.75pt] (1.100,1.150) .. controls (1.100,0.740) and (1.650,0.740) .. (1.650,1.150);
\draw[line width=.75pt] (0.000,0) .. controls (0.000,0.410) and (0.550,0.410) .. (0.550,0);
\draw[line width=.75pt] (1.100,0) .. controls (1.100,0.410) and (1.650,0.410) .. (1.650,0);
\fill (0.000,0) circle (.9pt);\fill (0.000,1.150) circle (.9pt);
\fill (0.550,0) circle (.9pt);\fill (0.550,1.150) circle (.9pt);
\fill (1.100,0) circle (.9pt);\fill (1.100,1.150) circle (.9pt);
\fill (1.650,0) circle (.9pt);\fill (1.650,1.150) circle (.9pt);
\end{scope}
\node[fill=white,inner sep=.2pt] at (5.375,1.46) {$\scriptstyle\star$};
\node[font=\small] at (4.825,-.42) {mark circle 2};
\node at (3,1.15) {$+$};\node at (7.1,1.15) {$\longmapsto\ 2\beta$};
\begin{scope}[shift={(9.2,0.42)}]
\draw[line width=.75pt] (0.000,1.450) .. controls (0.000,1.040) and (0.550,1.040) .. (0.550,1.450);
\draw[line width=.75pt] (1.100,1.450) .. controls (1.100,1.040) and (1.650,1.040) .. (1.650,1.450);
\draw[line width=.75pt] (0.000,0) .. controls (0.000,0.410) and (0.550,0.410) .. (0.550,0);
\draw[line width=.75pt] (1.100,0) .. controls (1.100,0.410) and (1.650,0.410) .. (1.650,0);
\fill (0.000,0) circle (.9pt);\fill (0.000,1.450) circle (.9pt);
\fill (0.550,0) circle (.9pt);\fill (0.550,1.450) circle (.9pt);
\fill (1.100,0) circle (.9pt);\fill (1.100,1.450) circle (.9pt);
\fill (1.650,0) circle (.9pt);\fill (1.650,1.450) circle (.9pt);
\node[font=\small] at (0.825,-.35) {$h=u_1u_3$};
\end{scope}
\end{tikzpicture}
\caption{The two contributions to $\Phi_\beta(h,h)$. Stack two copies of $h$ (separated by the thin gray horizontal rules); exactly two internal circles occur. The star records a chosen circle of weight $1$, while the unmarked one retains weight $\beta$. Thus the two terms each evaluate to $\beta h$.}
\label{fig:vis-mark}
\end{figure}

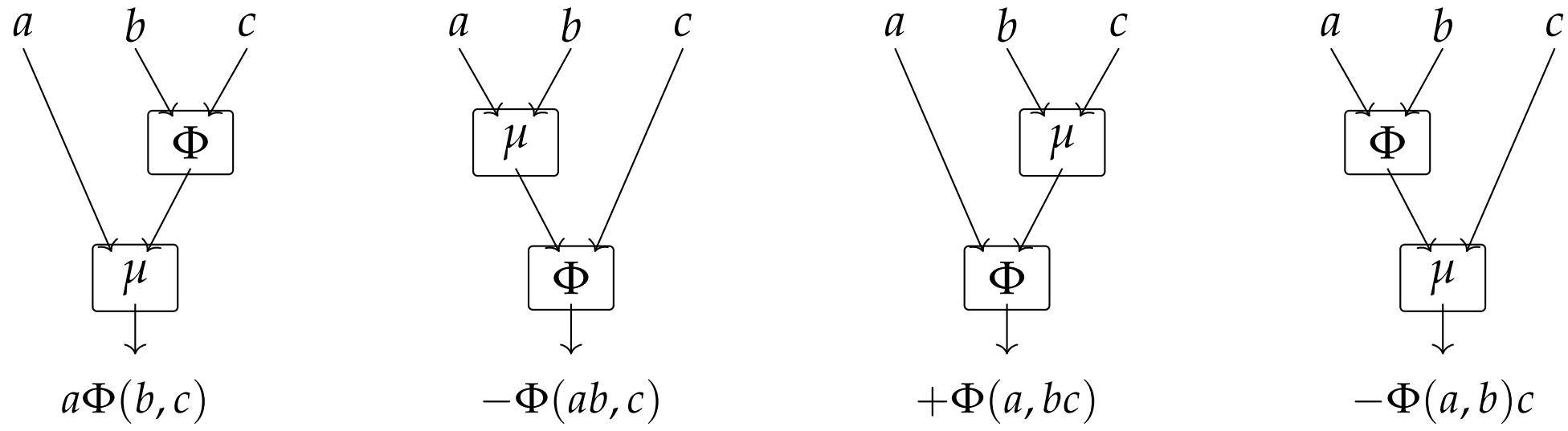
\begin{figure}[htbp]
\centering
\begin{tikzpicture}[scale=0.97]
\node (a0) at (-0.95,2.8) {$a$};
\node (b0) at (0,2.8) {$b$};
\node (c0) at (0.95,2.8) {$c$};
\node[draw,rounded corners=1pt,minimum width=.7cm,minimum height=.4cm] at (0.47,1.8) {$\Phi$};
\node[draw,rounded corners=1pt,minimum width=.7cm,minimum height=.4cm] at (0,0.65) {$\mu$};
\draw[->] (0,2.5999999999999996) -- (0.31999999999999995,2.02);
\draw[->] (0.95,2.5999999999999996) -- (0.62,2.02);
\draw[->] (-0.95,2.5999999999999996) -- (-0.2,0.87);
\draw[->] (0.47,1.58) -- (0.1,0.87);
\draw[->] (0,0.43000000000000005) -- (0,0);\node[font=\small] at (0,-0.4) {$a\Phi(b,c)$};
\node (a3x7) at (2.75,2.8) {$a$};
\node (b3x7) at (3.7,2.8) {$b$};
\node (c3x7) at (4.65,2.8) {$c$};
\node[draw,rounded corners=1pt,minimum width=.7cm,minimum height=.4cm] at (3.2300000000000004,1.8) {$\mu$};
\node[draw,rounded corners=1pt,minimum width=.7cm,minimum height=.4cm] at (3.7,0.65) {$\Phi$};
\draw[->] (2.75,2.5999999999999996) -- (3.0800000000000005,2.02);
\draw[->] (3.7,2.5999999999999996) -- (3.3800000000000003,2.02);
\draw[->] (4.65,2.5999999999999996) -- (3.9000000000000004,0.87);
\draw[->] (3.2300000000000004,1.58) -- (3.6,0.87);
\draw[->] (3.7,0.43000000000000005) -- (3.7,0);\node[font=\small] at (3.7,-0.4) {$-\Phi(ab,c)$};
\node (a7x4) at (6.45,2.8) {$a$};
\node (b7x4) at (7.4,2.8) {$b$};
\node (c7x4) at (8.35,2.8) {$c$};
\node[draw,rounded corners=1pt,minimum width=.7cm,minimum height=.4cm] at (7.87,1.8) {$\mu$};
\node[draw,rounded corners=1pt,minimum width=.7cm,minimum height=.4cm] at (7.4,0.65) {$\Phi$};
\draw[->] (7.4,2.5999999999999996) -- (7.72,2.02);
\draw[->] (8.35,2.5999999999999996) -- (8.02,2.02);
\draw[->] (6.45,2.5999999999999996) -- (7.2,0.87);
\draw[->] (7.87,1.58) -- (7.5,0.87);
\draw[->] (7.4,0.43000000000000005) -- (7.4,0);\node[font=\small] at (7.4,-0.4) {$+\Phi(a,bc)$};
\node (a11x1) at (10.15,2.8) {$a$};
\node (b11x1) at (11.1,2.8) {$b$};
\node (c11x1) at (12.049999999999999,2.8) {$c$};
\node[draw,rounded corners=1pt,minimum width=.7cm,minimum height=.4cm] at (10.629999999999999,1.8) {$\Phi$};
\node[draw,rounded corners=1pt,minimum width=.7cm,minimum height=.4cm] at (11.1,0.65) {$\mu$};
\draw[->] (10.15,2.5999999999999996) -- (10.479999999999999,2.02);
\draw[->] (11.1,2.5999999999999996) -- (10.78,2.02);
\draw[->] (12.049999999999999,2.5999999999999996) -- (11.299999999999999,0.87);
\draw[->] (10.629999999999999,1.58) -- (11.0,0.87);
\draw[->] (11.1,0.43000000000000005) -- (11.1,0);\node[font=\small] at (11.1,-0.4) {$-\Phi(a,b)c$};
\end{tikzpicture}
\caption{The four terms of $\delta\Phi=0$. Each box retains its number of input slots. The diagram expresses the derivative of associativity; it does not identify the four individual terms.}
\label{fig:vis-cocycle}
\end{figure}

With three inputs, differentiating the associativity relation for
$\mu_{\beta+t}$ gives the exact cochain identity
\begin{equation}\label{vis:cocycle}
 a\Phi(b,c)-\Phi(ab,c)+\Phi(a,bc)-\Phi(a,b)c=0.
\end{equation}
Figure~\ref{fig:vis-cocycle} displays its four terms. Which gluing first
creates a particular circle can depend on the parenthesization, but
the two first-order expansions of the same complete diagram agree.
This is why marking produces a cocycle rather than an arbitrary rule.

\subsection{Cup products preserve the input slots}
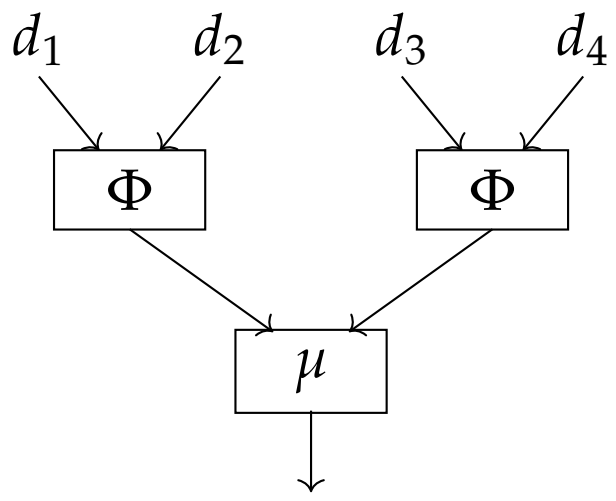
\begin{figure}[htbp]
\centering
\begin{tikzpicture}[scale=1]
\node at (0.0,3) {$d_1$};
\node at (1.2,3) {$d_2$};
\node at (2.4,3) {$d_3$};
\node at (3.5999999999999996,3) {$d_4$};
\node[draw,minimum width=1.0cm,minimum height=.5cm] at (0.6,2) {$\Phi$};
\draw[->] (0.0,2.75)--(0.39999999999999997,2.26);
\draw[->] (1.2,2.75)--(0.8,2.26);
\node[draw,minimum width=1.0cm,minimum height=.5cm] at (3.0,2) {$\Phi$};
\draw[->] (2.4,2.75)--(2.8,2.26);
\draw[->] (3.5999999999999996,2.75)--(3.2,2.26);
\node[draw,minimum width=1cm,minimum height=.5cm] at (1.8,.8) {$\mu$};\draw[->] (.6,1.74)--(1.55,1.06);\draw[->] (3,1.74)--(2.05,1.06);\draw[->] (1.8,.54)--(1.8,0);
\node[anchor=west,align=left] at (5,1.7) {$\Phi(d_1,d_2)\,\Phi(d_3,d_4)$\\[3mm] $\text{four inputs, degree }4$\\[3mm] $\text{output product: ordinary gluing}$};
\end{tikzpicture}
\caption{Cup multiplication: two independent marked-circle operations, followed by ordinary gluing of their outputs. Circles first created at the final gluing are not marked. This is not the second derivative of a bilinear product.}
\label{fig:vis-cup}
\end{figure}

The ordinary bar formula has no hidden choice of a middle cut:
\begin{equation}\label{vis:cup}
 (\Phi\smile\Phi)(d_1,d_2,d_3,d_4)
 =\Phi(d_1,d_2)\Phi(d_3,d_4).
\end{equation}
Each of the two operations marks its own internal gluing circle. Any
circle produced only when their outputs are multiplied has the usual
unmarked weight. This is a four-cochain. By contrast, the coefficient
of $t^2$ in $\mu_{\beta+t}(d_1,d_2)$ still has two inputs and is not a
four-cochain. Figure~\ref{fig:vis-cup} keeps these arities visible.

For example, on four strands,
\[
 (\Phi\smile\Phi)(u_1,u_1,u_3,u_3)=u_1u_3\ne0.
\]
At $\beta=\sqrt2$ the class of this cochain is nevertheless zero,
because the only nonsemisimple basic block is $\Lambda(2)$ and has
$HH^4=0$. At $\beta=1$ on six strands the nonsemisimple block is
$\Lambda(3)$, and its class $[\Phi_B]^2$ is nonzero by the already
proved block theorem. Nonzero values on single tuples do not decide
whether a cochain is a boundary.

The same distinction has a small formal-chain picture. In $\Lambda(3)$
let $\mathrm C=\mathrm C_{1;2,2}$; the degree $(2,2)$ part of its
shuffle diagonal is
\[
 \Delta_{2,2}(\mathrm C)=
 \mathrm C_{1;2,0}\otimes_\Lambda\mathrm C_{1;0,2}
 -\mathrm C_{2;1,1}\otimes_\Lambda\mathrm C_{2;1,1}.
\]
The first tensor evaluates to zero under $\eta_1\otimes\eta_1$, while
$\eta_1(\mathrm C_{2;1,1})=-e_3$. Hence
\[
 (\eta_1\smile\eta_1)(\mathrm C)=-e_3=-\eta_2(\mathrm C).
\]
The relevant paths are of the type depicted in the chain drawing of
Section~\ref{sec:differential}; none may be replaced by its algebra product. This calculation explains the
nonzero higher class after the Morita and bar comparisons, without
pretending that a cohomology class is a single unmarked staircase.

\subsection{A bracket picture and an explicit coboundary}
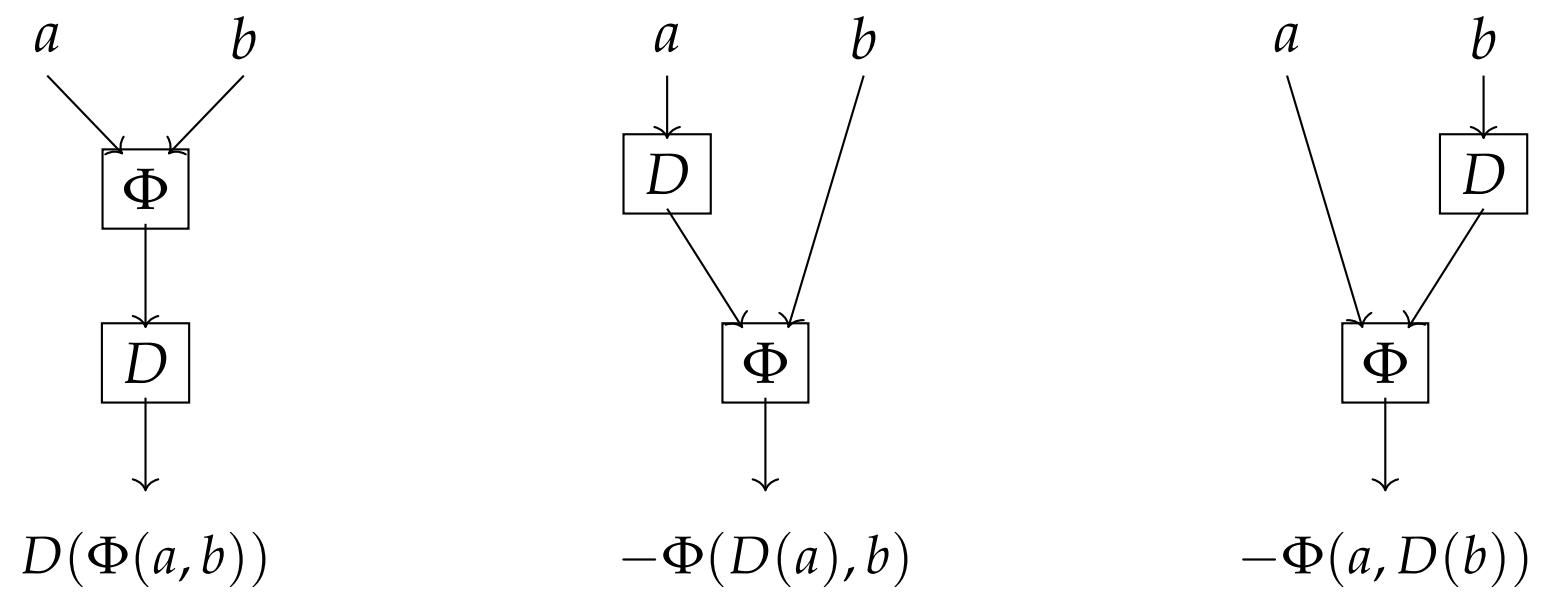
\begin{figure}[htbp]
\centering
\begin{tikzpicture}[scale=1]
\node at (-0.65,3) {$a$};\node at (0.65,3) {$b$};
\node[draw] at (0.0,2.0) {$\Phi$};\node[draw] at (0.0,.85) {$D$};
\draw[->] (-0.65,2.75)--(-0.15,2.23);\draw[->] (0.65,2.75)--(0.15,2.23);\draw[->] (0.0,1.77)--(0.0,1.08);
\draw[->] (0.0,.62)--(0.0,0);\node[font=\small] at (0.0,-.45) {$D(\Phi(a,b))$};
\node at (3.4499999999999997,3) {$a$};\node at (4.75,3) {$b$};
\node[draw] at (3.4499999999999997,2.1) {$D$};\node[draw] at (4.1,.85) {$\Phi$};
\draw[->] (3.4499999999999997,2.75)--(3.4499999999999997,2.33);\draw[->] (3.4499999999999997,1.87)--(3.9499999999999997,1.08);\draw[->] (4.75,2.75)--(4.25,1.08);
\draw[->] (4.1,.62)--(4.1,0);\node[font=\small] at (4.1,-.45) {$-\Phi(D(a),b)$};
\node at (7.549999999999999,3) {$a$};\node at (8.85,3) {$b$};
\node[draw] at (8.85,2.1) {$D$};\node[draw] at (8.2,.85) {$\Phi$};
\draw[->] (8.85,2.75)--(8.85,2.33);\draw[->] (8.85,1.87)--(8.35,1.08);\draw[->] (7.549999999999999,2.75)--(8.049999999999999,1.08);
\draw[->] (8.2,.62)--(8.2,0);\node[font=\small] at (8.2,-.45) {$-\Phi(a,D(b))$};
\end{tikzpicture}
\caption{The exact cochain formula for $[D,\Phi]$. In general its value is not pointwise $-\Phi$. In $TL_3(1)$ the correction is explicitly $[D,\Phi]+\Phi=\delta K$, with $K$ given in \eqref{vis:K}.}
\label{fig:vis-bracket}
\end{figure}

For a one-cochain $D$ and a two-cochain $\Phi$ the Gerstenhaber bracket
is the exact operation
\begin{equation}\label{vis:bracket}
 [D,\Phi](a,b)=D(\Phi(a,b))-\Phi(D(a),b)-\Phi(a,D(b)).
\end{equation}
The three terms are nested operations (Figure~\ref{fig:vis-bracket}),
unlike the independent outputs in a cup product. Already in $TL_3(1)$,
$[D,\Phi](u,u)=0$ whereas $-\Phi(u,u)=-u$. Thus the class relation
$[[D],[\Phi]]=-[\Phi]$ must not be drawn as pointwise equality.

In this example the discrepancy can be given explicitly. Define a
normalized linear cochain $K$ by
\begin{equation}\label{vis:K}
 \begin{array}{c|ccccc}
 a&1&u&v&uv&vu\\\hline
 K(a)&0&u&uv-v&0&u.
 \end{array}
\end{equation}
Then, as ordinary bilinear cochains,
\begin{equation}\label{vis:Kidentity}
 \boxed{[D,\Phi]+\Phi=\delta K.}
\end{equation}
To verify this, all pairs with a unit vanish. In the order
$u,v,uv,vu$ the common table on the remaining pairs is
\[
 \begin{array}{c|cccc}
 &u&v&uv&vu\\\hline
 u&u&uv&uv&u\\
 v&0&0&0&0\\
 uv&0&0&0&0\\
 vu&vu&v&v&vu
 \end{array}
\]
where the table is obtained by substituting \eqref{vis:K} into
$aK(b)-K(ab)+K(a)b$. The same substitution of the known $D$ and
$\Phi$ into the left side gives exactly this table.
Consequently the equality of classes is directly checked, and the
correcting one-cochain is itself expressed in the planar diagram basis.

The displayed cochain identities can also be checked directly in the
finite diagram basis; the general blockwise results are proved above.

\section{The zero-parameter boundary}\label{tl:zero}
	The nonzero-parameter block theorem must not be specialized blindly
	at $e=2$.  We state the separate zero-parameter conclusions and work out the
	smallest complete example.
	
	\subsection{An alternating-weight derivation}
	For every $N$ the formula
	\begin{equation}\label{tl:Dzero}
		D_N^{(0)}(u_i)=(-1)^{i+1}u_i
	\end{equation}
	defines a derivation of $TL_N(0)$.  The square relation is preserved
	because $u_i^2=0$.  If $\varepsilon_i=(-1)^{i+1}$, an adjacent
	triple has weight $2\varepsilon_i+\varepsilon_{i+1}=\varepsilon_i$,
	the same as its reduction $u_i$.  Distant relations also have equal
	weights on both sides.  Thus the algebra is graded by the signed
	weight of the generators, and $D^{(0)}$ is its Euler derivation.
	
	For odd $N$, $TL_N(0)$ is semisimple over $\mathbb C$
	\cite[Proposition~8.2]{RSA}, so positive cohomology vanishes and
	this derivation is inner.  For $N=3$ one can see it explicitly:
	\[
	h=\tfrac12(u_1u_2-u_2u_1),\qquad
	[h,u_1]=u_1,\quad [h,u_2]=-u_2.
	\]
	For even $N=2m$, the derivation is outer.  The central recursion
	\eqref{tl:Frec}, valid polynomially at zero, becomes
	\[
	F_0=0,\quad F_1=-2,\quad
	F_n=-F_{n-2}-u_{n-1}F_{n-1}-F_{n-1}u_{n-1}.
	\]
	By induction the even $F_{2m}$ have weight one and the odd ones
	weight zero.  Moreover
	\[
	\operatorname{coeff}_{u_1}(F_{2m})=4(-1)^{m-1}\ne0:
	\]
	Here $\operatorname{coeff}_{u_1}$ denotes the coefficient of the
	Catalan basis diagram $u_1$. Since $F_2=4u_1$, each subsequent even
	step changes this coefficient's
	sign; the terms involving the new last generator cannot produce the
	diagram $u_1$ because they retain a last upper or lower arc.
	Hence
	\[
	F_{2m}\in Z(TL_{2m}(0)),\qquad
	D^{(0)}(F_{2m})=F_{2m}\ne0,
	\]
	proving outer nonvanishing.  This does not contradict innerness on
	the next odd level: the retraction in
	Lemma~\ref{tl:outer-retraction} divides by $\beta$ and does not
	exist in that form at zero.
	
	\subsection{The marked-circle operation at zero}
	There is no need to interpret a power $0^0$ in the formula.  Take
	the coefficient of $t$ in $t^Ld_0$:
	\begin{equation}\label{tl:Phizero}
		\Phi_0(d_1,d_2)=
		\begin{cases}d_0,&L=1,\\0,&L=0\text{ or }L\ge2.\end{cases}
	\end{equation}
	It is a two-cocycle, but its class is not uniformly the generator
	of the second cohomology in all even-strand cases.
	
	\subsection{The dual numbers: the full smallest exceptional ring}
	For $A=TL_2(0)=\mathbb C[u]/(u^2)$, the derivation $D(u)=u$ is
	outer because the algebra is commutative.  However
	$\Phi_0(u,u)=u$ is a coboundary: for $h(1)=0$, $h(u)=1/2$,
	\[
	(\delta h)(u,u)=uh(u)+h(u)u=u.
	\]
	A genuine second-degree generator is the normalized cochain
	$\Psi(u,u)=1$, not the loop-parameter derivative.
	The normalized bar complex has one generator $u^{\otimes n}$ in
	each degree, so the cochain complex is
	\[
	A\xrightarrow{0}A\xrightarrow{2u}A\xrightarrow{0}A
	\xrightarrow{2u}A\longrightarrow\cdots.
	\]
	It follows that $HH^{2r}(A)=\mathbb C[\Psi^{\smile r}]$ for
	$r\ge1$ and $HH^{2r+1}(A)=\mathbb C[\Psi^{\smile r}\smile D]$
	for $r\ge0$.  Writing $X=[\Psi]$, $Y=[D]$ gives
	\begin{equation}\label{tl:dualnumbers-ring}
		HH^*(A)=\mathbb C[u,X,Y]/(u^2,uX,uY,Y^2),
		\quad |u|=0,\ |X|=2,\ |Y|=1.
	\end{equation}
	The even powers do not terminate, in contrast to
	\eqref{cup:eq:ring}.  In particular the finite resolution of the
	nonzero-parameter basic algebra cannot be transferred to the entire
	even-strand zero-parameter family by changing only a scalar.
	
	\begin{remark}[Boundary of the claim]
		This manuscript proves the alternating outer derivation on all even
		zero-parameter levels and gives the complete $N=2$ example.  It does
		not claim a uniform ring theorem for all even $TL_N(0)$ or a
		uniform nonvanishing theorem for $[\Phi_0]$.  The odd case is
		complete through semisimplicity.  This limitation does not affect the
		nonzero-parameter block theorem.
	\end{remark}

\section{Loop-weight response and an ice-type realization}
\label{phys:section}
This section has two distinct ingredients. The loop/Potts and TL/six-vertex
dictionaries are established constructions \cite{TemperleyLieb,BKW,BJJZ}.
The statements about $D,\Phi$, their cup product, and their bracket are
applications of the cochains proved above. We work with finite planar
geometries or open strips. Affine or periodic TL and topological weights
for noncontractible loops are not included in the cohomology theorem.
In particular, none of the following arguments asserts a thermodynamic
limit, a critical exponent, or a new solution of an ice model.

\subsection{Marking as differentiation of a finite loop sum}
Let $\Omega$ be a finite set of completely specified loop configurations
with a fixed boundary closure. Suppose all circles are contractible and
all local weights $w(\omega)$ are independent of $\beta$. Set
\begin{equation}\label{phys:Z}
 Z(\beta)=\sum_{\omega\in\Omega}w(\omega)\beta^{\ell(\omega)},
\end{equation}
where $\ell(\omega)$ is the number of circles. This is the standard
loop-fugacity formulation, specialized here to a finite planar setting;
see \cite[Sections~2.1 and 4]{BJJZ} for loop transfer matrices.
A termwise derivative gives
\begin{equation}\label{phys:mark-derivative}
 Z'(\beta)=\sum_{\omega\in\Omega}\ \sum_{c\in\mathrm{Loops}(\omega)}
 w(\omega)\beta^{\ell(\omega)-1}.
\end{equation}
Thus a derivative selects one marked circle. At one multiplication step,
this is precisely \eqref{vis:two-loops} and the operation $\Phi$.

When $\beta>0$, $w(\omega)\ge0$ and $Z(\beta)>0$, normalize the weights
to a probability distribution. Direct differentiation, without a
large-volume limit, yields
\begin{equation}\label{phys:expectation}
 \beta\partial_\beta\log Z=\langle\ell\rangle,
 \qquad (\beta\partial_\beta)^2\log Z=
 \langle\ell^2\rangle-\langle\ell\rangle^2.
\end{equation}
The second equality follows by differentiating the normalized weights:
$\beta\partial_\beta p_\omega=(\ell(\omega)-\langle\ell\rangle)p_\omega$.
For complex or signed weights, \eqref{phys:mark-derivative} remains an
algebraic identity, but the probabilistic language is not justified.
The ordinary second derivative counts two distinct marked circles and
weights them by $\ell(\ell-1)$; it is not a cup square of a bilinear
operation, whose four inputs must be preserved.

\paragraph{Internal gluing and boundary closure.}
A physical partition function may be written schematically as
$\operatorname{cl}_\beta(T_\beta)$, where $\operatorname{cl}_\beta$ is a
specified closing functional and $T_\beta$ is a product of local faces.
Its full derivative is
\[
 \partial_\beta\operatorname{cl}_\beta(T_\beta)
 =(\partial_\beta\operatorname{cl}_\beta)(T_\beta)
   +\operatorname{cl}_\beta(\partial_\beta T_\beta).
\]
The first term marks circles created only by closure; $\Phi$ accounts
for circles at internal multiplication steps. If
$T_r=\mu_\beta(T_{r-1},X_r)$, then
\[
 \dot T_r=\Phi(T_{r-1},X_r)+\mu_\beta(\dot T_{r-1},X_r)
                  +\mu_\beta(T_{r-1},\dot X_r).
\]
Local Boltzmann weights that themselves depend on $\beta$ contribute
the final term. Omitting these terms would incorrectly identify the
full derivative of a transfer matrix with one application of $\Phi$.

\subsection{The planar Potts normalization}
For a finite connected plane graph $G=(V,E)$, the random-cluster form is
\begin{equation}\label{phys:FK}
 Z_G(Q,v)=\sum_{A\subseteq E}v^{|A|}Q^{k(A)},
\end{equation}
where $k(A)$ counts connected components, including isolated vertices.
For positive integer $Q$, this follows by expanding
$\prod_{ij\in E}(1+v\delta_{\sigma_i,\sigma_j})$ and summing over spins.
The expression then defines a polynomial for general $Q$; this is the
Fortuin--Kasteleyn representation \cite{FK}.

Take a regular neighbourhood in the sphere of the occupied subgraph
$(V,A)$. Let $\ell(A)$ be the number of boundary components. These are
the medial loops. Every component of this neighbourhood has genus zero,
and it deformation retracts to the corresponding graph component. Hence
\[
 |V|-|A|=2k(A)-\ell(A),\qquad
 \ell(A)=2k(A)+|A|-|V|.
\]
Writing $Q=\beta^2$ with $\beta\ne0$ and $x=v/\beta$, we obtain the
fully normalized identity
\begin{equation}\label{phys:Potts-loop}
 \boxed{Z_G(\beta^2,v)=
 \beta^{|V|}\sum_{A\subseteq E}x^{|A|}\beta^{\ell(A)}.}
\end{equation}
This is the planar loop correspondence of \cite{BKW}; the elementary
Euler calculation specifies exactly the normalization used here.
On a torus the regular neighbourhood may have positive genus, so the
same formula without a topological correction is not asserted.

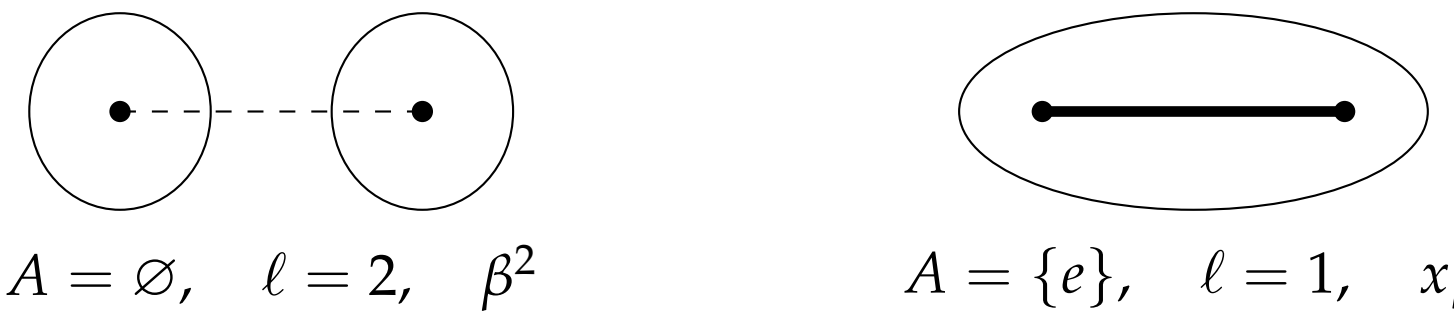
\begin{figure}[htbp]
\centering
\begin{tikzpicture}[scale=1]
\draw[dashed] (0.0,0)--(2.0,0);
\fill (0.0,0) circle (2pt);\fill (2.0,0) circle (2pt);
\draw (0.0,0) ellipse (.6 and .65);\draw (2.0,0) ellipse (.6 and .65);
\node at (1.0,-1.1) {$A=\varnothing,\quad \ell=2,\quad\beta^2$};
\draw[line width=2pt] (6.1,0)--(8.1,0);
\fill (6.1,0) circle (2pt);\fill (8.1,0) circle (2pt);
\draw (7.1,0) ellipse (1.55 and .65);
\node at (7.1,-1.1) {$A=\{e\},\quad\ell=1,\quad x\beta$};
\end{tikzpicture}
\caption{A one-edge planar random-cluster example. Boundary components of a regular neighbourhood of the occupied subgraph are the medial loops. After the factor $\beta^{|V|}$ is removed, the two contributions are $\beta^2$ and $x\beta$, where $x=v/\beta$.}
\label{fig:phys-FK}
\end{figure}

For a single edge and its two vertices,
\[
 Z_G(Q,v)=Q^2+Qv,\qquad
 \beta^{-2}Z_G(\beta^2,\beta x)=\beta^2+x\beta.
\]
The empty occupied set has two boundary loops, and the occupied edge
has one (Figure~\ref{fig:phys-FK}). At fixed $x$, differentiation gives
$2\beta+x$, exactly the sum of configurations with one marked circle.
The mean loop count in this normalized two-state model is
$(2\beta+x)/(\beta+x)$ when $\beta,x>0$.

The qualification ``fixed $x$'' is important. At fixed $v$, $x=v/\beta$
also varies; marked loops alone do not give $\partial_\beta Z_G$.
Similarly $Q\partial_Q\log Z_G$ at fixed $v$ counts FK clusters,
not the same observable as \eqref{phys:expectation}.
For the family $\beta=2\cos(\pi/e)$, one has
$Q=4\cos^2(\pi/e)$; $e=3,4,6$ give $Q=1,2,3$.
The blockwise nontriviality of $\Phi$ is an obstruction to trivializing
an algebra deformation, not a proof of thermodynamic singularity at
these numbers.

\subsection{A local face and six ice-rule entries}
\begin{figure}[htbp]
\centering
\begin{tikzpicture}[scale=1]
\node[draw,minimum width=1.8cm,minimum height=1.6cm] at (.8,.8) {$X_i(s)$};\node at (2.4,.8) {$=$};\node at (3.5,.8) {$a(s)$};
\begin{scope}[shift={(4.3,0)}]
\draw[line width=.75pt] (0.000,1.600) .. controls (0.000,0.928) and (0.000,0.672) .. (0.000,0);
\draw[line width=.75pt] (0.700,1.600) .. controls (0.700,0.928) and (0.700,0.672) .. (0.700,0);
\fill (0.000,0) circle (.9pt);\fill (0.000,1.600) circle (.9pt);
\fill (0.700,0) circle (.9pt);\fill (0.700,1.600) circle (.9pt);
\end{scope}
\node at (6.0,.8) {$+$};\node at (7.1,.8) {$b(s)$};
\begin{scope}[shift={(8,0)}]
\draw[line width=.75pt] (0.000,1.600) .. controls (0.000,1.190) and (0.700,1.190) .. (0.700,1.600);
\draw[line width=.75pt] (0.000,0) .. controls (0.000,0.410) and (0.700,0.410) .. (0.700,0);
\fill (0.000,0) circle (.9pt);\fill (0.000,1.600) circle (.9pt);
\fill (0.700,0) circle (.9pt);\fill (0.700,1.600) circle (.9pt);
\end{scope}
\end{tikzpicture}
\caption{A local face is a weighted sum of two connectivity states. Here $a(s)=\sin(\lambda-s)$, $b(s)=\sin s$, and $\beta=2\cos\lambda$. Spectator strands have been omitted.}
\label{fig:phys-face}
\end{figure}
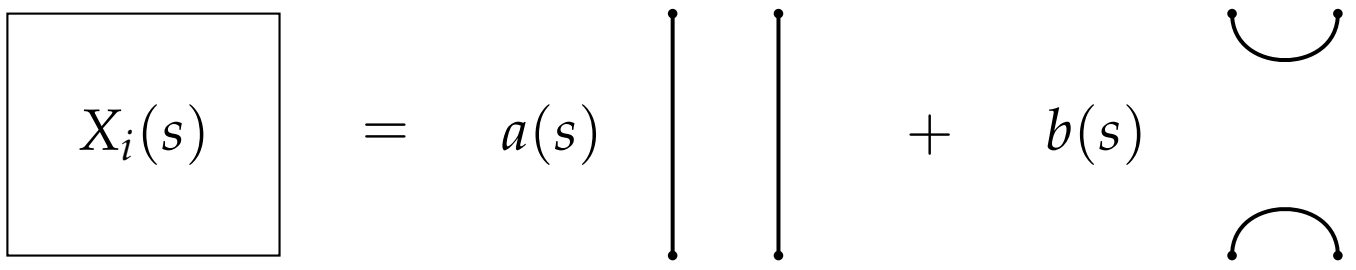

Let $q=e^{\mathrm i\lambda}$, $\beta=2\cos\lambda$, and define the
local algebra element
\begin{equation}\label{phys:face}
 X_i(s)=a(s)1+b(s)u_i,\qquad
 a(s)=\sin(\lambda-s),\quad b(s)=\sin s.
\end{equation}
This is the standard TL form of the local checked $R$ operator
\cite[Section~3, equation (3.2)]{BJJZ}, with our positive loop convention.
Figure~\ref{fig:phys-face} shows its two connectivity states. We shall
verify the representation used here rather than leaving a sign convention
implicit.

On $V=(\mathbb C^2)^{\otimes N}$, with two-site basis ordered as
$\uparrow\uparrow,\uparrow\downarrow,\downarrow\uparrow,
\downarrow\downarrow$, put
\begin{equation}\label{phys:E}
 E(q)=\begin{pmatrix}
 0&0&0&0\\0&q^{-1}&-1&0\\0&-1&q&0\\0&0&0&0
 \end{pmatrix},\qquad
 \rho_q(u_i)=1^{\otimes(i-1)}\otimes E(q)\otimes
 1^{\otimes(N-i-1)}.
\end{equation}
The $2\times2$ middle block immediately gives $E^2=(q+q^{-1})E$.
For completeness let
\[
 v=q^{-1}|\uparrow\downarrow\rangle-|\downarrow\uparrow\rangle,
 \qquad w=\langle\uparrow\downarrow|-q\langle\downarrow\uparrow|.
\]
Then $E=vw$ and $wv=\beta$. The image of $E_{12}$ on three sites is
spanned by $v\otimes|\uparrow\rangle$ and $v\otimes|\downarrow\rangle$;
direct substitution gives $E_{12}E_{23}$ equal to the identity on these
two vectors. Thus $E_{12}E_{23}E_{12}=E_{12}$; exchanging the tensor
positions gives the other adjacent relation. Distant generators act on
disjoint tensor factors. Therefore \eqref{phys:E} defines a representation
of the open TL algebra. It is a conventional six-vertex gauge; signs and
root conventions in \cite{BJJZ} must be translated before comparing matrices.

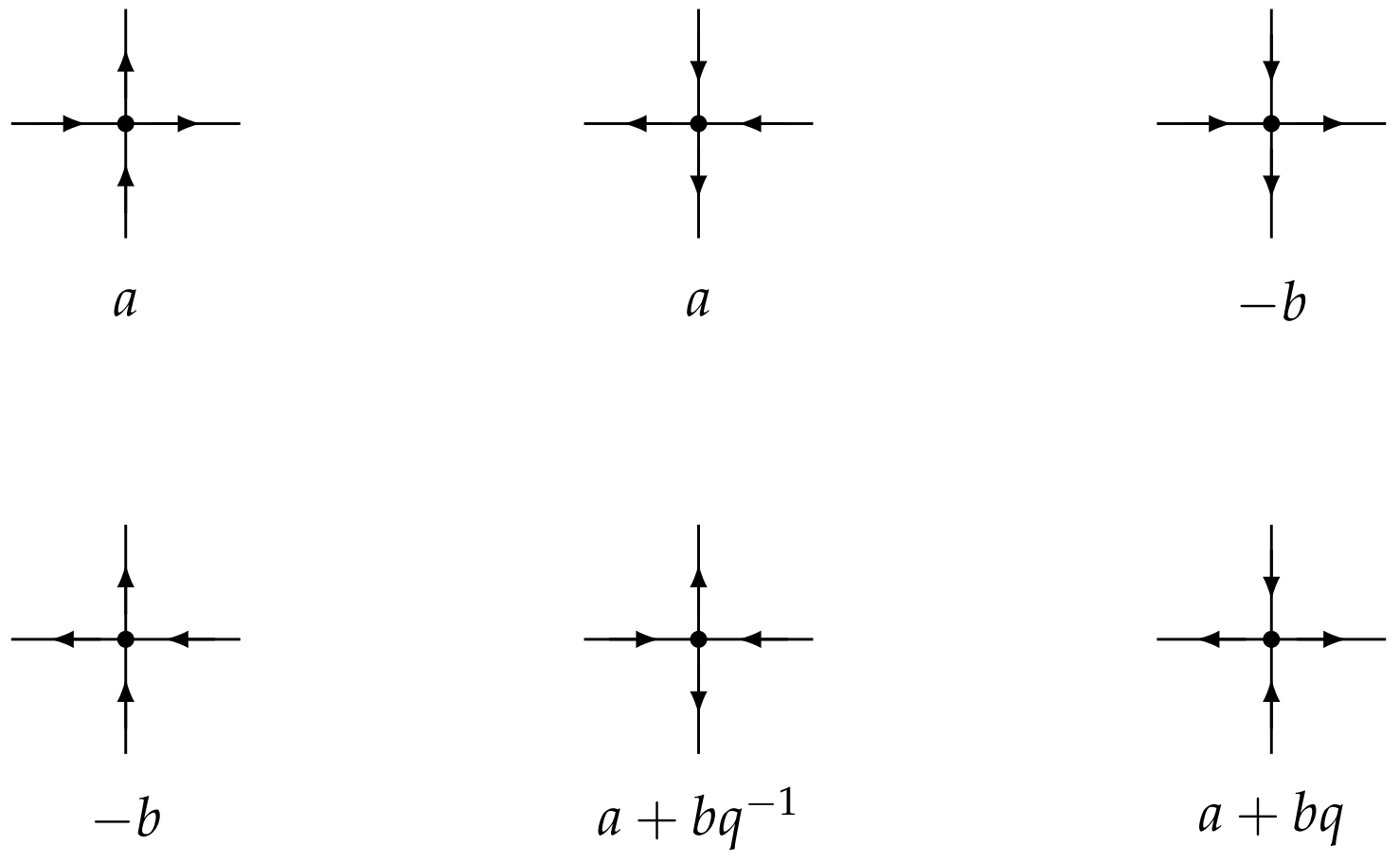
\begin{figure}[htbp]
\centering
\begin{tikzpicture}[scale=1]
\draw[line width=.55pt] (-0.8,3.6)--(0.0,3.6);
\draw[-{Latex[length=1.65mm]},line width=.55pt] (-0.6240000000000001,3.6000000000000005)--(-0.272,3.6000000000000005);
\draw[line width=.55pt] (0.0,2.8)--(0.0,3.6);
\draw[-{Latex[length=1.65mm]},line width=.55pt] (0.0,2.976)--(0.0,3.3280000000000003);
\draw[line width=.55pt] (0.0,3.6)--(0.8,3.6);
\draw[-{Latex[length=1.65mm]},line width=.55pt] (0.17600000000000002,3.6000000000000005)--(0.528,3.6000000000000005);
\draw[line width=.55pt] (0.0,3.6)--(0.0,4.4);
\draw[-{Latex[length=1.65mm]},line width=.55pt] (0.0,3.7760000000000002)--(0.0,4.128);
\fill (0.0,3.6) circle (1.7pt);\node at (0.0,2.35) {$a$};
\draw[line width=.55pt] (3.2,3.6)--(4.0,3.6);
\draw[-{Latex[length=1.65mm]},line width=.55pt] (3.8240000000000003,3.6000000000000005)--(3.4720000000000004,3.6000000000000005);
\draw[line width=.55pt] (4.0,2.8)--(4.0,3.6);
\draw[-{Latex[length=1.65mm]},line width=.55pt] (4.0,3.4240000000000004)--(4.0,3.072);
\draw[line width=.55pt] (4.0,3.6)--(4.8,3.6);
\draw[-{Latex[length=1.65mm]},line width=.55pt] (4.624,3.6000000000000005)--(4.272,3.6000000000000005);
\draw[line width=.55pt] (4.0,3.6)--(4.0,4.4);
\draw[-{Latex[length=1.65mm]},line width=.55pt] (4.0,4.224)--(4.0,3.8720000000000008);
\fill (4.0,3.6) circle (1.7pt);\node at (4.0,2.35) {$a$};
\draw[line width=.55pt] (7.2,3.6)--(8.0,3.6);
\draw[-{Latex[length=1.65mm]},line width=.55pt] (7.376,3.6000000000000005)--(7.728000000000001,3.6000000000000005);
\draw[line width=.55pt] (8.0,2.8)--(8.0,3.6);
\draw[-{Latex[length=1.65mm]},line width=.55pt] (8.0,3.4240000000000004)--(8.0,3.072);
\draw[line width=.55pt] (8.0,3.6)--(8.8,3.6);
\draw[-{Latex[length=1.65mm]},line width=.55pt] (8.176,3.6000000000000005)--(8.528,3.6000000000000005);
\draw[line width=.55pt] (8.0,3.6)--(8.0,4.4);
\draw[-{Latex[length=1.65mm]},line width=.55pt] (8.0,4.224)--(8.0,3.8720000000000008);
\fill (8.0,3.6) circle (1.7pt);\node at (8.0,2.35) {$-b$};
\draw[line width=.55pt] (-0.8,0)--(0.0,0);
\draw[-{Latex[length=1.65mm]},line width=.55pt] (-0.17600000000000002,0.0)--(-0.528,0.0);
\draw[line width=.55pt] (0.0,-0.8)--(0.0,0);
\draw[-{Latex[length=1.65mm]},line width=.55pt] (0.0,-0.6240000000000001)--(0.0,-0.272);
\draw[line width=.55pt] (0.0,0)--(0.8,0);
\draw[-{Latex[length=1.65mm]},line width=.55pt] (0.6240000000000001,0.0)--(0.272,0.0);
\draw[line width=.55pt] (0.0,0)--(0.0,0.8);
\draw[-{Latex[length=1.65mm]},line width=.55pt] (0.0,0.17600000000000002)--(0.0,0.528);
\fill (0.0,0) circle (1.7pt);\node at (0.0,-1.25) {$-b$};
\draw[line width=.55pt] (3.2,0)--(4.0,0);
\draw[-{Latex[length=1.65mm]},line width=.55pt] (3.3760000000000003,0.0)--(3.728,0.0);
\draw[line width=.55pt] (4.0,-0.8)--(4.0,0);
\draw[-{Latex[length=1.65mm]},line width=.55pt] (4.0,-0.17600000000000002)--(4.0,-0.528);
\draw[line width=.55pt] (4.0,0)--(4.8,0);
\draw[-{Latex[length=1.65mm]},line width=.55pt] (4.624,0.0)--(4.272,0.0);
\draw[line width=.55pt] (4.0,0)--(4.0,0.8);
\draw[-{Latex[length=1.65mm]},line width=.55pt] (4.0,0.17600000000000002)--(4.0,0.528);
\fill (4.0,0) circle (1.7pt);\node at (4.0,-1.25) {$a+bq^{-1}$};
\draw[line width=.55pt] (7.2,0)--(8.0,0);
\draw[-{Latex[length=1.65mm]},line width=.55pt] (7.824,0.0)--(7.472000000000001,0.0);
\draw[line width=.55pt] (8.0,-0.8)--(8.0,0);
\draw[-{Latex[length=1.65mm]},line width=.55pt] (8.0,-0.6240000000000001)--(8.0,-0.272);
\draw[line width=.55pt] (8.0,0)--(8.8,0);
\draw[-{Latex[length=1.65mm]},line width=.55pt] (8.176,0.0)--(8.528,0.0);
\draw[line width=.55pt] (8.0,0)--(8.0,0.8);
\draw[-{Latex[length=1.65mm]},line width=.55pt] (8.0,0.6240000000000001)--(8.0,0.272);
\fill (8.0,0) circle (1.7pt);\node at (8.0,-1.25) {$a+bq$};
\end{tikzpicture}
\caption{The six nonzero vertex weights of $R=P\check X$ in the stated gauge. Left and bottom are input channels, right and top are output channels; a $+$ spin points east or north. Each vertex has two arrows in and two out. The weights can be complex or signed: positivity is not asserted in this gauge.}
\label{fig:phys-ice}
\end{figure}

Set $\check X=aI+bE$ and let $P$ interchange the two tensor factors.
The unchecked vertex operator is
\begin{equation}\label{phys:ice-matrix}
 R=P\check X=
 \begin{pmatrix}
 a&0&0&0\\
 0&-b&a+bq&0\\
 0&a+bq^{-1}&-b&0\\
 0&0&0&a
 \end{pmatrix}.
\end{equation}
For generic $s$ there are precisely six nonzero entries; each preserves
the total number of up spins. With left and bottom channels as inputs,
this is the two-arrows-in, two-arrows-out rule drawn in
Figure~\ref{fig:phys-ice}. At the trigonometric values in \eqref{phys:face},
$a+bq^{\pm1}=\sin\lambda\,e^{\pm\mathrm i s}$.
The weights in this gauge need not be nonnegative. We do not identify
this generic operator with the equal-weight square-ice model, nor claim
that changing gauge leaves arbitrary boundaries unchanged.

\subsection{What the cochains do to local faces}
Hold the scalar coefficients of the local faces fixed when evaluating
a cochain on them. Bilinearity and normalization give
\begin{equation}\label{phys:Phi-face}
 \Phi_\beta(X_i(s),X_i(t))=b(s)b(t)u_i.
\end{equation}
Only the product of their two cup--cap terms creates a circle. If
$\lambda$, and hence $a,b$, are also varied, their derivatives are
additional terms, as explained above.
For the derivation,
\[
 D(X_i(s))=b(s)D(u_i),\qquad
 D(X_{i_1}(s_1)\cdots X_{i_r}(s_r))=
 \sum_jX_{i_1}(s_1)\cdots D(X_{i_j}(s_j))\cdots X_{i_r}(s_r).
\]
This is an explicit algebraic insertion rule for a face product. It does
not by itself show that $D$ preserves a chosen transfer matrix or its
trace, or that $D$ is a physical conserved charge.

Over dual numbers, $g_\varepsilon=\operatorname{id}+\varepsilon D$
is an automorphism of the fixed gluing algebra, with inverse
$\operatorname{id}-\varepsilon D$. Transporting a two-cochain gives
\begin{equation}\label{phys:transport-Phi}
 g_\varepsilon\Phi(g_\varepsilon^{-1}a,g_\varepsilon^{-1}b)
 =\Phi(a,b)+\varepsilon[D,\Phi](a,b).
\end{equation}
This explains the bracket diagram as the response of the marked-circle
rule to an infinitesimal algebra automorphism. On the blocks of the main
theorem, its cohomology class scales by $1-\varepsilon$, because
$[[D_B],[\Phi_B]]=-[\Phi_B]$. An exact cochain identity may contain a
coboundary, as \eqref{vis:Kidentity} demonstrates.

\subsection{A crucial limitation: changing coefficient representations}
\begin{proposition}\label{phys:rho-boundary}
Let a representation of the deformed algebra on a fixed vector space
have expansion $\rho_t=\rho+t\rho_1$ and
$\mu_t=\mu+t\Phi\pmod{t^2}$. Give $\operatorname{End}(V)$ its
$A$-bimodule structure through $\rho$. Then
\begin{equation}\label{phys:rhoPhi}
 \boxed{\rho\circ\Phi=\delta\rho_1
 \quad\text{in }C^2(A,\operatorname{End}(V)).}
\end{equation}
\end{proposition}
\begin{proof}
The coefficient of $t$ in
$\rho_t(\mu_t(a,b))=\rho_t(a)\rho_t(b)$ is
\[
 \rho(\Phi(a,b))+\rho_1(ab)
 =\rho(a)\rho_1(b)+\rho_1(a)\rho(b),
\]
which is exactly \eqref{phys:rhoPhi}.
\end{proof}
This proposition applies to \eqref{phys:E} at every parameter covered by
the nonzero-root theorem: $q\ne\pm1$ and a local choice $q(t)$ satisfies
$q(t)+q(t)^{-1}=\beta+t$ with $q'(0)=q^2/(q^2-1)$. For the local matrix,
\[
 E'=\begin{pmatrix}
 0&0&0&0\\0&-1/(q^2-1)&0&0\\
 0&0&q^2/(q^2-1)&0\\0&0&0&0
 \end{pmatrix}.
\]
Differentiating $E^2=\beta E$ gives
$EE'+E'E-\beta E'=E$, the matrix image of
$\Phi(u_i,u_i)=u_i$. Thus the first-order deformation is visible in
matrix entries while its class with these new coefficients is a boundary.

There is no contradiction with $[\Phi_B]\ne0$ with regular coefficients.
The coefficient map $A\to\operatorname{End}(V)$ need not induce an
injective map on Hochschild cohomology. The nontrivial algebra class and
a scalar response in a fixed statistical model are distinct assertions.
The coefficient complex of Appendix~\ref{app:coefficients} is the tool
for asking what survives in other representations; the main theorem is
not a classification of all physical observables.

\section{Coefficients in an arbitrary bimodule}\label{app:coefficients}
	The resolution itself is not restricted to regular coefficients.
	Let $M$ be a $\Lambda(m)$-bimodule.  Its cochain space is
	\[
	\mathsf C^n(\Lambda(m),M)=
	\bigoplus_{\substack{k\ge1,\ u,v\ge0\\u+v=n,\ k+\max(u,v)\le m}}
	e_{k+u}Me_{k+v}.
	\]
	The base index is denoted by $k$ here to distinguish it from the root
	parameter in the TL chapters. For a fixed degree-$(n-1)$ cochain
	$\phi$, put $M_{k;u,v}=\phi(\mathrm C_{k;u,v})$ when $u+v=n-1$.
	For $n\ge1$, the full coefficient
	differential evaluated on a degree-$n$ generator is
	\begin{equation}\label{app:general-coeff-d}
		\begin{aligned}
			(\mathsf d^{n-1}\phi)(\mathrm C_{k;u,v})={}&
			\mathbf{a}_{k+u-1}M_{k;u-1,v}
			+(-1)^{u+v}M_{k;u,v-1}\mathbf{b}_{k+v-1}\\
			&+(-1)^u \mathbf{b}_{k+u}M_{k+1;u,v-1}
			+(-1)^u M_{k+1;u-1,v}\mathbf{a}_{k+v}.
		\end{aligned}
	\end{equation}
	Inadmissible terms are zero.  Products in this formula are the left
	and right actions on $M$, not multiplication of two elements of
	$\Lambda(m)$ unless $M=\Lambda(m)$.
	
	The reduction to the three families $\mathrm D,\mathrm U,\mathrm V$ used
	$e_i\Lambda(m)e_j=0$ for $|i-j|\ge2$.  That vanishing need not hold
	for $e_iMe_j$.  For example, let $S_i$ be the vertex simple and take
	$M=S_i\otimes S_j^*$.  All arrows act by zero, so the differential
	\eqref{app:general-coeff-d} vanishes.  A nonzero coordinate requires
	$k+u=i$, $k+v=j$; therefore
	\begin{equation}\label{app:simple-bimodule}
		HH^n(\Lambda(m),S_i\otimes S_j^*)\cong
		\begin{cases}
			\Bbbk,&n=|i-j|,|i-j|+2,\ldots,i+j-2,\\
			0,&\text{otherwise}.
		\end{cases}
	\end{equation}
	The same argument is valid over the base field $\Bbbk$ with its split
	vertex simples.  In particular, for $i=3,j=1$ the class is supported
	on $\mathrm C_{1;2,0}$ in degree two, a generator whose regular-coefficient
	value was forced to zero.  Deleting such a generator from the
	resolution would give an incorrect answer for this bimodule.
	
	If $B$ is a TL block and $fBf\cong\Lambda(m)$ is a full basic corner,
	then for any $B$-bimodule $M$,
	\[
	HH^n(B,M)\cong HH^n(\Lambda(m),fMf).
	\]
	The coefficients, not just the resolution, must be transferred.
	For finite-dimensional left modules $V,W$ the familiar bridge is
	\[
	HH^n(B,\operatorname{Hom}_{\mathbb C}(V,W))
	\cong\operatorname{Ext}_B^n(V,W),
	\qquad (afb)(v)=a f(bv).
	\]
	This gives a useful coefficient extension of the explicit resolution;
	it is not a classification of all bimodules nor an assertion that
	the regular generators $D,\Phi$ generate these other cohomology groups.

\section*{Declaration of AI-assisted research and manuscript preparation}

The author used OpenAI's ChatGPT during the development of this work and the preparation of the manuscript. The assistance included mathematical discussions, exploration of proof strategies, symbolic computations, verification of examples, suggestions for arguments, and the preparation and revision of parts of the mathematical exposition, LaTeX source, and TikZ diagrams.

AI-generated suggestions were treated as research proposals requiring mathematical verification, rather than as independent proofs or authoritative sources. The author takes full responsibility for the mathematical correctness, originality, references, and final content of the manuscript.

\end{document}